\documentclass[dvipsnames]{article}
\usepackage{graphicx} 

\usepackage[citestyle=alphabetic-verb,bibstyle=alphabetic]{biblatex}
\usepackage{mathtools}
\usepackage{enumitem}
\usepackage{amssymb}
\usepackage{stmaryrd}
\usepackage{xcolor}
\usepackage{longtable}
\usepackage{hyperref}
\usepackage[mode=buildnew]{standalone}
\usepackage{tikz}
\usepackage{tikz-3dplot}
\usetikzlibrary{arrows,shapes,positioning,calc,3d,decorations.pathreplacing,decorations.markings,cd,patterns}
\usepackage{xargs}
\tikzstyle{->-}= [postaction={decorate, decoration={markings,mark=at position .55 with {\arrow{to}}}}]
\tikzstyle{-<-}= [postaction={decorate, decoration={markings,mark=at position .45 with {\arrowreversed{to};}}}]

\usetikzlibrary{calc}

\usepackage{geometry}

\usepackage{array}

\DeclareMathOperator{\CAT}{CAT}
\DeclareMathOperator{\Cay}{Cay}
\DeclareMathOperator{\lk}{lk}
\DeclareMathOperator{\st}{st}
\DeclareMathOperator{\Sym}{Sym}

\DeclareMathOperator{\Aut}{Aut}
\DeclareMathOperator{\Alt}{Alt}
\DeclareMathOperator{\PGL}{PGL}
\DeclareMathOperator{\SU}{SU}

\DeclareMathOperator{\PGU}{PGU}
\DeclareMathOperator{\GU}{GU}
\DeclareMathOperator{\id}{id}
\DeclareMathOperator{\score}{score}
\DeclareMathOperator{\FF}{\mathbb F}
\DeclareMathOperator{\QQ}{\mathbb Q}
\DeclareMathOperator{\RR}{\mathbb R}
\DeclareMathOperator{\ZZ}{\mathbb Z}
\DeclareMathOperator{\NN}{\mathbb N}
\DeclareMathOperator{\KK}{\mathbb K}
\DeclareMathOperator{\LL}{\mathbb L}

\DeclareMathOperator{\calO}{\mathcal O}
\DeclareMathOperator{\SL}{SL}
\DeclareMathOperator{\GL}{GL}
\DeclareMathOperator{\PSL}{PSL}

\DeclareMathOperator{\Span}{span}

\DeclareMathOperator{\Radu}{R}
\DeclareMathOperator{\JW}{JW}
\DeclareMathOperator{\type}{\textsc{typ}}
\DeclareMathOperator{\im}{im}

\DeclareMathOperator{\Tr}{Tr}
\DeclareMathOperator{\Cl}{Cl}

\DeclareFontFamily{U}{mathb}{\hyphenchar\font45}
\DeclareFontShape{U}{mathb}{m}{n}{
      <5> <6> <7> <8> <9> <10> gen * mathb
      <10.95> mathb10 <12> <14.4> <17.28> <20.74> <24.88> mathb12
      }{}
\DeclareSymbolFont{mathb}{U}{mathb}{m}{n}
\DeclareFontSubstitution{U}{mathb}{m}{n}
\DeclareMathSymbol{\righttoleftarrow}      {3}{mathb}{"FD}

\newcommand{\defeq}{\mathrel{\mathop{:}}=}

\newcommand{\abs}[1]{\lvert #1 \rvert}

\newcommand{\norm}[1]{\lVert #1 \rVert}
\newcommand{\gen}[1]{\langle #1 \rangle}

\newcommand{\op}{\text{op}}

\newcommandx{\Deltaopv}[1][1={}]{\Delta^\op_{*,#1}}

\newcommand{\Stab}{\operatorname{Stab}}

\renewcommand{\setminus}{\smallsetminus}

\newcommand{\finres}[1]{{#1}^{(\infty)}}
\newcommand{\sgind}[2]{[#1:#2]}
\newcommand{\sqcomp}{S}
\newcommand{\sradu}{\sqcomp_\text{R}}
\newcommand{\sjw}{\sqcomp_\text{JW}}
\newcommand{\tradu}{\tilde{\sqcomp}_\text{R}}
\newcommand{\tjw}{\tilde{\sqcomp}_\text{JW}}
\newcommand{\gradu}{\Gamma_\text{R}}
\newcommand{\gjw}{\Gamma_\text{JW}}

\newcommand{\bfG}{\mathbf{G}}

\usepackage{amsthm}
\newtheorem{introtheorem}{Theorem}

\newtheorem{introcorollary}[introtheorem]{Corollary}

\newtheorem{theorem}{Theorem}[section]
\newtheorem{lemma}[theorem]{Lemma}
\newtheorem{corollary}[theorem]{Corollary}
\newtheorem{proposition}[theorem]{Proposition}

\theoremstyle{definition}
\newtheorem{definition}[theorem]{Definition}
\newtheorem{example}[theorem]{Example}
\theoremstyle{remark}

\newtheorem{remark}[theorem]{Remark}

\title{Non-residually finite $\tilde{C}_2$-lattices}
\author{Thomas Titz Mite, Stefan Witzel}
\date{August 2026}

\newcommand{\newcomment}[4]{%
\newcounter{#2counter}
\expandafter\newcommand\csname #1\endcsname[1]{%
\refstepcounter{#2counter}%
{\color{#4}(#3\arabic{#2counter})}\marginpar{\scriptsize\raggedright\textbf{\color{#4}(#2 \arabic{#2counter}):} ##1}%
}}
\reversemarginpar

\definecolor{darkgreen}{rgb}{0,0.6,0}
\newcomment{sw}{Stefan}{S}{darkgreen}
\newcomment{ttm}{Thomas}{T}{blue}

\begin{document}

\maketitle

\begin{abstract}\noindent
  We provide the first known examples of non-residually finite lattices on irreducible buildings. They contain the first known simple $\CAT(0)$-groups with property (T), and the first known $\CAT(0)$-groups that are not quasi-isometric to a direct product.

  We also classify type-preserving vertex-regular lattices on buildings of type $\tilde{A}_2$ and thickness three, and discover an arithmetic example that is not commensurable to previously studied lattices.
\end{abstract}

\section{Introduction}

Classically, the lattices on Euclidean buildings that were most studied are $S$-arithmetic subgroups of reductive groups over global fields, which are residually finite. In contrast to these, there are lattices on two-dimensional exotic Euclidean buildings, many of them known quite explicitly \cite{Ronan84,Tits85,Kantor86,CartwrightManteroStegerZappa93,Barre00,Essert,Radu_NonDesarguesian,Witzel17}, and conjecturally none of them is residually finite (see \cite[Conjecture~1.5]{BaderCapraceLecureux}, \cite[Conjecture~B]{LecureuxWitzel2026}). However, until now, non-residual finiteness could not be verified for a single one. We provide the first such examples.

\begin{introtheorem}[Main theorem]\label{thm:main}
  There are five finite triangle complexes $Y_1^2, Y_1^3, Y_2^3, Y_3^3, Y_4^3$ whose universal covers $X_i^q \defeq \tilde Y_i^q$ are exotic buildings of type $\tilde C_2$ and whose fundamental groups $\Gamma_i^q \defeq \pi_1(Y_i^q)$ are not residually finite. More precisely, the finite residual $\check{\Gamma}_i^q \defeq (\Gamma_i^q)^{(\infty)}$ has finite index. The buildings $X_i^q$ are pairwise not isomorphic and their fundamental groups pairwise not quasi-isometric. Since the $\Gamma_i^q$ are uniform building lattices they are $\CAT(0)$-groups and have Kazhdan's property~(T).
\end{introtheorem}

Together with \cite[Theorem~A]{LecureuxWitzel2026} we deduce:

\begin{introcorollary}
  The groups $\check{\Gamma}_1^2, \check{\Gamma}_1^3, \check{\Gamma}_2^3, \check{\Gamma}_3^3, \check{\Gamma}_4^3$ are simple $\CAT(0)$ groups that have property (T) and are not quasi-isometric to a proper direct product.
\end{introcorollary}

The building $X_i^q$ has thickness $q+1$ and $\Gamma_i^q$ acts on it preserving types. The action is regular on vertices of each special type. Each $\Gamma_i^q$ admits an extension $\bar{\Gamma}^q_i = \Gamma^q_i \rtimes C_2$ by involutory building automorphism that is not type-preserving. The groups $\bar{\Gamma}^q_i$ act regularly on the special vertices of their buildings.

For $q = 2$ the group $\Gamma^2_1$ is its own finite residual $(\Gamma^2_1)^{(\infty)}$ and $\bar{\Gamma}^2_1$ is the full automorphism group $\Aut(X^2_1)$ of the building. For $q=3$ the extension $\bar{\Gamma}^3_i$ is not unique and the full automorphism group is strictly larger while the finite residual is strictly smaller.
More precisely, letting $\hat \Gamma_i^q \defeq \Aut(X_i^q)$ denote the automorphism group of the building and $\check \Gamma_i^q \defeq (\Gamma_i^q)^{(\infty)}$ the finite residual we have
\begin{align*}
  \sgind{\hat \Gamma_i^q}{\Gamma_i^q} &\defeq \begin{cases} 2 & q=2 \\    8 & q=3\;,  \end{cases} \;
                                     &
  \sgind{\Gamma_i^q}{\check\Gamma_i^q} &\defeq \begin{cases} 1 & q=2 \\ 4 & q=3, i =1,2\\ 8 & q=3, i =3,4\;.\\\end{cases}
\end{align*}

The group $\bar \Gamma_1^2$ admits a geometrically natural presentation with $15$ generators and $24$ relations of lengths $2$ and $4$, see Proposition~\ref{prop:presentation_bar_2}.
The lattices $\bar \Gamma_i^3$ admit analogous presentations with $40$ generators and $64$ relations, see Appendix~\ref{app:C2tilde_thickness4}.

The complexes $Y^q_i$ were found by a computationally expensive computer search. Now that they have been found, computer assistance is needed only for the following routine verifications: 1. to verify that $\check{\Gamma}^q_i$ is the finite residual by checking that a certain group presentation presents the trivial group; 2. to verify that $\hat{\Gamma}^q_i$ is the full automorphism group of the building by reconstructing certain balls and showing that they are rigid; 3. to show that the buildings for $q = 3$ are pairwise non-isomorphic by showing that the quotients $\check{\Gamma}^q_i \backslash X^q_i$ are non-isomorphic. Code carrying out these verifications is provided in the GitHub repository \cite{TitzMiteCode}.

Trees are Euclidean buildings, as are their products. Hence, the irreducible lattices of products of trees first studied by Wise, Burger and Mozes \cite{Wise96,BurgerMozes00} are the first known non-residually finite lattices on two-dimensional buildings, many of them simple. The decomposition into a product of trees is crucially used in proving that they are not residually finite. Note that from a geometric perspective, all uniform lattices on products of trees represent a single quasi-isometry class, while there are infinitely many quasi-isometry classes of non-arithmetic lattices on irreducible buildings.

In constructing the lattices in the Main Theorem we make use of non-residually finite lattices on products of trees. Indeed, in each case we started with a non-residually finite lattice $\Lambda$ on a product of trees $T_1 \times T_2$ and built the irreducible building $X$ around $T_1 \times T_2$ in such a way that the action of $\Lambda$ extends to a cocompact lattice $\Gamma$ on $X$. So the fact that $\Gamma$ is not residually finite is ensured by the construction. It remains an open problem how, given an arbitrary lattice $\Gamma < \Aut(X)$ on an exotic building $X$, one can identify a non-trivial element of the finite residual.

\medskip

As mentioned before, the complexes $Y^q_i$ were found by a greedy search through a certain parameter space. In smaller settings this parameter space can be completely enumerated. As an illustration of this we obtain a second result which is independent of the above: we classify all lattices that act regularly on the vertices of each type in a building of type $\tilde{A}_2$ of thickness $3$: there are $13$ such lattices. If we consider a lattice from \cite{CartwrightManteroStegerZappa93}, it acts on the set of types as $C_3$ and the kernel of this action is of the kind we consider, but only four of the $13$ lattices admit such an extension to a CMSZ-lattice. Three of the lattices are $S$-arithmetic in characteristic $2$ and three are $S$-arithmetic in characteristic $0$ while the remaining $7$ act on exotic buildings. All of the arithmetic ones are commensurable to previously studied groups except for one, which is therefore of particular interest:

\begin{introtheorem}\label{thm:intro_arithmetic}
    Let $K \defeq \QQ(\sqrt{-23})$ with ring of integers ${\cal O}_K = \ZZ[(1+\sqrt{-23})/2]$ and denote the Galois involution by $\overline{\rule[1.5ex]{.1em}{0pt}\cdot\rule[1ex]{.1em}{0pt}}$.
    Let $\omega \defeq (143 + 7\sqrt{-23})/2$ and consider the hermitian matrix
  \[
    Q = 
    \left(\begin{array}{rrr}
    83 & \omega  & \omega \\
    \bar \omega  & 83 & \omega \\
    \bar \omega & \bar \omega & 83
    \end{array}\right) \in M_3({\cal O}_K) \cap \GL_3(K).
  \]
  It defines a hermitian form $f(v,w) = v^T Q \overline{w}$ on ${\cal O}_{K}^3$ and an associated general unitary group
  \[
    \GU(f)(R) = \{M \in \GL_3(R) \mid M^TQ\overline{M} = \lambda Q \text{ for some } \lambda \in R^\times\}.
  \]
  The $S$-arithmetic subgroup $\GU(f)(\calO_K[1/2])$ acts as a lattice on the Bruhat--Tits building associated to $\PGL_3(\QQ_2)$. The action is transitive on vertices of each type and vertex stabilizers have order at least $3$.
\end{introtheorem}

The paper is organized as follows. Section~\ref{sec:buildings_lattices} collects basic results on two-dimensional buildings, their lattices, and residual finiteness. Section~\ref{sec:tree_products} collects results and examples on products of trees needed later on. Section~\ref{sec:non-rf_c2tilde} is the core of the article where we analyze the buildings $X^q_i$, the lattices $\Gamma^q_i$ and prove the Main Theorem. The complexes $Y^q_i$ were found via a computer search that has proven useful for other purposes; we describe it in Section~\ref{sec:searching}. Section~\ref{sec:a2tilde} contains the classification of vertex regular $\tilde{A}_2$-lattices including Theorem~\ref{thm:intro_arithmetic}. Appendix \ref{app:C2tilde_thickness4} contains descriptions of the complexes and lattices in Theorem~\ref{thm:main} for $q = 3$. Appendix \ref{app:a2tildeGABs} contains descriptions of the complexes that give rise to the $\tilde A_2$-lattices in Section \ref{sec:a2tilde}.

\paragraph{Acknowledgments.} The authors would like to thank Pierre-Emmanuel Caprace and Bernhard Mühlherr for helpful discussions. We thank Steffen Kionke and Tim Steger for valuable input on how to determine the arithmetic structure of the lattice in Theorem~\ref{thm:intro_arithmetic}. The second author was partially funded through the DFG Heisenberg project WI 4079/6, which supported both authors.


\section{Buildings and lattices}\label{sec:buildings_lattices}

\subsection{Non-positively curved complexes}\label{sec:non-positively_curved_complexes}

The spaces we will be interested in are non-positively curved chamber complexes in the following sense (see \cite[Chapter~I.7]{BridsonHaefliger} for context). Let $C$ be a Eulidean or hyperbolic \emph{polygon}, that is, a bounded intersection of finitely many half-planes in the Euclidean or hyperbolic plane with non-empty interior. We say that $C$ is a \emph{Coxeter polygon} if its interior angles are of the form $\pi/m$ with $m \in \NN, m \ge 2$. Most of the time $C$ will be a triangle. Gluing copies of $C$ along isometries of their edges results in a \emph{polygonal complex} $X$ that is \emph{pure} in the sense that every face (vertex, edge, polygon) is contained in a two-dimensional one. We say that it is \emph{thin}, \emph{firm}, respectively \emph{thick} if every edge is contained in exactly two, at least two, respectively at least three polygons. The complex is \emph{locally finite} if every point is contained in finitely many polygons. We regard each edge of $C$ as having its own type and stipulate that the gluings respect types in order to make $X$ \emph{typed} and let $I$ denote the set of types. The type of a vertex is the two-element subset of $I$ of the types of edges incident with it; we denote by $I_2$ the set of types of vertices. We take the angle in the vertex of type $\{i,j\}$ to be $\pi/m_{ij} = \pi/m_{ji}$. The complex $X$ carries a canonical metric that makes it complete geodesic, see \cite[Theorem~I.7.50]{BridsonHaefliger}. Note that cells may meet in more than one face, for instance there may be doubled edges, but the types ensure that there are no identifications within a closed cell, for instance no loops.

The link of a vertex $v \in X$ is a metric graph whose vertices correspond to edges containing $v$ and whose edges correspond to polygons containing it (the previous remark on types means that we need not talk about half-edges and corners here). The length of the edge corresponding to the polygon $C$ is the angle of $C$ in $v$, i.e.\ $\pi/m_{ij}$ if $v$ is of type $\{i,j\}$.

We say that a locally finite, firm, typed polygonal complex $X$ is a \emph{chamber complex} if it has connected vertex links. This implies that any two polygons that can be connected by a sequence of polygons in which any two consecutive ones share an edge. We call the polygons of a chamber complex \emph{chambers}.

We say that a polygonal complex $X$ is \emph{non-positively curved} if it is locally CAT($0$), cf.\ \cite[Chapter~II.1]{BridsonHaefliger}. This is the case if and only if every vertex link has (metric) girth at least $2\pi$, see \cite[Theorem~II.5.2, Lemma~II.5.6]{BridsonHaefliger}.

Throughout the article we will be concerned with non-positively curved chamber complexes. We add the attribute \emph{metric} when we want to make the distinction with the complexes we will describe next. One important use of non-positive curvature condition is through the Cartan--Hadamard Theorem \cite[Theorem~II.4.1]{BridsonHaefliger}:

\begin{theorem}
  If $X$ is a connected non-positively curved chamber complex then $\tilde{X}$ is CAT($0$). In particular, it is contractible.
\end{theorem}

\smallskip

There is a combinatorial description of this setting that we will be using later on in the situation where $C$ is a simplex. Namely let $Y$ be a locally finite two-dimensional $\Delta$-complex (in the sense of \cite[Section~2.1]{Hatcher02}) that is pure in the sense that every point is contained in a triangle, firm in the sense that every edge is contained in at least two triangles, and typed in the sense that every vertex has a cotype in a three element set $I$ and cotypes are respected by the gluings (the type of a vertex is the complement of its cotype in $I$). The link of every vertex $v$ is a (simplicial) graph which we assume to be connected. It has a (combinatorial) girth $g(v)$ and we define the angle of cotype $i \in I$ to be
\[
  \alpha(i) \defeq \max\left\{\frac{2\pi}{g(u)} \, \big | \, u \text{ is of type } i  \right\}.
\]
We call such $Y$ a \emph{combinatorial chamber complex} and say that it is \emph{non-positively curved} if it satisfies the condition
\begin{equation}\label{eq:comb_npc}
  \alpha(1) + \alpha(2) + \alpha(3) \le \pi.
\end{equation}

The connection is made by the observation that if we take $C$ to be a triangle with angles $\alpha(1),\alpha(2),\alpha(3)$, which is hyperbolic and unique if the sum is less than $\pi$ and is Euclidean and unique up to scaling if the sum is $\pi$, and equip each triangle of $Y$ with the metric on $C$, then we obtain a metric non-positively curved chamber complex. The observation readily generalizes to other polyhedra but we will not have use for that.

Another important use of non-positive curvature is important is to deduce $\pi_1$-injectivity. The following is \cite[Proposition II.4.14]{BridsonHaefliger}

\begin{lemma}\label{lem:locally_convex}
  Let $Y$ be a complete non-positively curved geodesic metric space and let $Y_0 \subseteq Y$ be a locally convex subspace. Let $p \colon X \to Y$ and $p_0 \colon X_0 \to Y_0$ be universal covers and let $x \in X$ and $x_0 \in X$ be such that $p(x) = y = p_0(x_0)$. Then the canonical maps $(X_0,x_0) \to (X,x)$ and $\pi_1(X_0, y) \to \pi_1(X, y)$ are injective.

  This applies in particular if $Y$ is a connected, non-positively curved chamber complex and $Y_0$ is a connected non-positively curved chamber subcomplex.
\end{lemma}    

\subsection{Buildings}

Prime examples of non-positively curved chamber complexes are buildings and their quotients, which we discuss next, keeping the conventions from above. Standard references for buildings are \cite{AbramenkoBrown,Rousseau}.

We say that a graph is thin, firm, respectively thick if every vertex has degree two, at least two, respectively at least three.
A \emph{generalized $m$-gon} (or just \emph{generalized polygon} if we do not want to specify $m$) is a bipartite firm metric graph all of whose edges have length $\pi/m$ such that the diameter is $\pi$ and the girth is $2\pi$. Again one may say that a \emph{combinatorial generalized $m$-gon} is a bipartite firm graph that has diameter $m$ and combinatorial girth $2m$. Generalized triangles ($3$-gons), quadrangles ($4$-gons), respectively hexagons ($6$-gons) are \emph{buildings of type} $A_2$, $C_2$, respectively~$G_2$.

A two-dimensional Euclidean building is a connected, simply connected chamber complex $X$ in which $C$ is a square (type $\tilde{A}_1 \times \tilde{A}_1$), a triangle with all angles $\pi/3$ (type $\tilde{A}_2$), a triangle with angles $\pi/2, \pi/4, \pi/4$ (type $\tilde{C}_2$) or angles $\pi/2, \pi/3, \pi/6$ (type $\tilde{G}_2$) such that any two points are contained in a $2$-flat (isometric copy of the Euclidean plane) called an \emph{apartment} (among other things this implies that $X$ is non-positively curved hence $\CAT(0)$). Note that due to the type-restriction every apartment is necessarily tessellated in such a way that the chambers are the fundamental domains of a reflection group acting on it. A two-dimensional hyperbolic building is a simply connected chamber complex in which $C$ is a hyperbolic polygon such that any two points are contained in an isometric copy of the hyperbolic plane, again called an \emph{apartment}. The type-restriction ensures that $C$ is a Coxeter polygon and every apartment has its chambers as fundamental domains of a reflection group. To make the connection with \cite[Definition~4.1]{AbramenkoBrown} observe that apartments $\Sigma,\Sigma' \subseteq X$ are convex being geodesic subspaces of the CAT($0$)-space $X$. Thus their intersection $\Sigma \cap \Sigma'$ is convex and there is an isometry $\Sigma \to \Sigma'$ extending $\id_{\Sigma \cap \Sigma'}$.

Replacing the assumption that the vertex links have girth $\ge 2\pi$ by the stronger assumption that the links be generalized polygons leads to a stronger conclusion in the Cartan--Hadamard Theorem. The following statement is a rephrasing of \cite{Tits}, the metric characterization in \cite{CharneyLytchak} has much weaker assumptions.

\begin{theorem}\label{thm:building_cartan-hadamard}
  Let $X$ be a thick chamber complex whose chambers have angle $\pi/m_{ij}$ in vertices of type~$\{i,j\}$ and in which every link of a vertex of type $\{i,j\}$ is a generalized $m_{\{i,j\}}$-gon. Then $\tilde{X}$ is a Euclidean or hyperbolic building. It is Euclidean if $C$ is a Euclidean polygon (so $\sum 1/m_{\{i,j\}} = \abs{I}-2$) and hyperbolic if $C$ is a hyperbolic polygon (so $\sum 1/m_{\{i,j\}} < \abs{I}-2$). It is of type $\tilde{A}_1 \times \tilde{A}_1$, $\tilde{A}_2$, $\tilde{C}_2$, respectively $\tilde{G}_2$ if $(m_{\{i,j\}})_{\{i,j\} \in I_2}$ is $(2,2,2,2)$, $(3,3,3)$, $(2,4,4)$, respectively $(2,3,6)$.
\end{theorem}

\begin{proof}
  We want to apply \cite{Tits} and need to translate between chamber complexes whose links are generalized polygons and geometries. We first consider the case when chambers are triangles as it is easier. The complex $X$ corresponds to a geometry of type $M = (m_i)_{i \in I}$ with underlying set $X^{(0)}$, types in $I$ and incidence given by adjacency. The crucial assumption is that rank-$2$-residues, which correspond to $1$-dimensional links, are generalized polygons.

  The assumption that vertex links are connected ensures that a chamber complex is connected if and only if its corresponding geometry is connected. A covering in the sense of \cite{Tits} is a surjective morphism of geometries such that the restriction to any proper residue is an isomorphism and a geometry is simply connected if a covering by a connected geometry is an isomorphism. From this it is clear that coverings of chamber complexes correspond to coverings of geometries and, in particular, the geometry corresponding to $\tilde{X}$ is simply connected. As it is its own only rank-$3$-residue, this makes it residually simply connected and hence \cite[Theorem~1]{Tits} asserts that $\tilde{X}$ is a building. Note also that the $\tilde{A}_1 \times \tilde{A}_1$-case asserts that a simply connected square complex with complete bipartite vertex links is a product of trees which is easy to see. So the Euclidean cases, that we will be interested in, are covered by this.

  In the general case (where $C$ may not be a triangle) put $m_{ij} = \infty$ for distinct $i,j \in I$ with $\{i,j\} \not\in I_2$. Two chambers are $i$-adjacent if their edge of type $i$ is the same and they are $J$-adjacent for $J \subseteq I$ if they can be connected by a sequence of chambers that are $j$-adjacent for some $j \in J$. The residue of type $J$ of a chamber is the subcomplex consisting of all chambers that are $J$-adjacent to it. An important observation is that $J$-residues are locally convex. Indeed, if $v$ is a vertex of type $\{i,j\}$ in the boundary where $j \in J$ and $i \in I \setminus J$ the fact that the angle $\pi/m_{ij}$ is non-obtuse ensures that the residue is locally convex in $v$. It follows that residues of $\tilde{X}$ are contractible and, in particular, simply connected, which for $3$-residues eventually leads to residual simple connectedness, the crucial assumption of \cite[Theorem~1]{Tits}.

\begin{figure}
  \begin{center}
\tdplotsetmaincoords{70}{110}
\begin{tikzpicture}[tdplot_main_coords,scale=4]
  \coordinate (A) at (0,0,0);
  \coordinate (B) at (1,0,0);
  \coordinate (C) at (0.5,0.8660254,0);
  \coordinate (D) at (0.5,0.2886751,0.8164966);
  \draw[thick,blue] (A)--(B)--(C)--cycle;
  \draw[thick,blue] (A)--(D)--(B);
  \draw[thick,blue] (C)--(D);
  \draw[thick] ($(A)!0.5!(B)$) -- ($(A)!0.5!(C)$) -- ($(C)!0.5!(D)$) -- ($(B)!0.5!(D)$) -- cycle;
  \node[fill=white] at (A) {1};
  \node[fill=white] at (B) {2};
  \node[fill=white] at (C) {4};
  \node[fill=white] at (D) {3};
\end{tikzpicture}
\end{center}
\caption{We consider a chamber complex $X$ based on a chamber $C$ that is a square with set of types $I = \{1,2,3,4\}$ and $m_{i,i+1} = 2$ (indices modulo $4$). Shown is a chamber of $X$ together with the corresponding $3$-simplex of $S$. The label $i$ indicates that the vertex has type $I \setminus \{i\}$ so the opposite facet, as well as the edge of $C$ that it contains, have type $i$.}
  \label{fig:tits_correspondence}
\end{figure}

  To explain the translation, first define the simplicial complex $S$ whose simplices of type $I \setminus J$ and dimension $\abs{I \setminus J} - 1$ are the $J$-residues (see Figure~\ref{fig:tits_correspondence} for illustration). Thus a vertex of type $\{i\} \subseteq I$ is a residue of type $I \setminus \{i\}$. (Roughly speaking, this is reversing the Davis realization.) The geometry corresponding to $X$ has as underlying set the vertex set of $S$ and incidence is given by adjacency. It is of type $(m_{ij})_{i,j \in I}$ since the $2$-residues are generalized $m_{ij}$-gons if $\{i,j\} \in I_2$ and are trees, i.e.\ generalized $\infty$-gons if $\{i,j\} \not\in I_2$.

  A $J$-residue corresponds to a component of the full subcomplex of $S$ on edges of type in $J$. If $J$ has $3$ elements, a covering in the sense of \cite{Tits} corresponds to a type-preserving, surjective simplicial map that restricts to an isomorphism on any $J \setminus \{i\}$-residue (component of the corresponding subcomplex). In other words it restricts to isomorphisms on vertex links and therefore is a covering. So in order to show residual simple connectedness in the sense of \cite{Tits} we need to show that subcomplexes of $S$ corresponding to $3$-residues are simply connected. The final observation is that $S$ is the nerve of the cover of $X$ by $(I \setminus \{i\})$-residues which have convex hence contractible intersection if they intersect at all. It follows by the Nerve Cover Lemma \cite[Corollary~4G.3]{Hatcher02} that residues in $C$ are homotopy equivalent to their corresponding subcomplexes in $S$ and, in particular, $3$-residues are simply connected.
\end{proof}

Borrowing from \cite{Kantor86} we call a complex satisfying the assumptions of Theorem~\ref{thm:building_cartan-hadamard} a \emph{geometry that is almost a building} or \emph{GAB} for short.

\subsection{Lattices}

Let $X$ be locally finite Euclidean or hyperbolic building. We let $\Aut(X)$ denote the group of automorphisms of $X$, i.e.\ bijections $X \to X$ that take chambers isometrically to chambers. They are isometries on $X$ but can be described combinatorially as we just did. We equip $\Aut(X)$ with the topology of pointwise convergence (induced by the product topology on $X^X$) which coincides with compact open topology: an open neighborhood of $\varphi \in \Aut(X)$ is defined by coinciding with $\varphi$ on finitely many points.

A \emph{uniform lattice} on $X$ is a group $\Gamma$ acting on $X$ with finite stabilizers such that $\Gamma \backslash X$ is a finite complex. If the action of $\Gamma$ on $X$ is free then $X \to \Gamma \backslash X$ is a (universal) covering and $\Gamma \backslash X$ is a GAB.

\subsection{Residual finiteness, simplicity, and the normal subgroup property}

\begin{definition}
  Let $\Gamma$ be a group. The \emph{finite residual} of $\Gamma$ is
  \[
    \finres\Gamma = \bigcap \{\Lambda < \Gamma \mid [\Gamma:\Lambda] < \infty\} = \bigcap \{\Lambda \lhd \Gamma \mid [\Gamma:\Lambda] < \infty\}.
  \]
  The group $\Gamma$ is \emph{residually finite} if $\finres\Gamma = \{1\}$.

  It is \emph{just infinite} if every non-trivial normal subgroup $\{1\} \ne N \lhd \Gamma$ is finite index, $[\Gamma:N] < \infty$. It is \emph{hereditarily just infinite} if every finite index subgroup is just infinite.

  A uniform lattice $\Gamma$ on a complex $X$ via $\alpha \colon \Gamma \to \Aut(X)$ is said to have the \emph{normal subgroup property} if $\im\alpha \cong \Gamma/\ker\alpha$ is just infinite.
\end{definition}

In the classical setting of arithmetic lattices $\ker \alpha$ is central but in general it could be any finite group. Some authors will take residual finiteness as part of the definition of (hereditarily) just infiniteness but we do not.

\begin{proposition}
  If $\Gamma$ is hereditarily just infinite then it is either residually finite or $\finres\Gamma$ is of finite index and simple.
\end{proposition}

\begin{proof}
  Since $\finres\Gamma$ is normal, if it is non-trivial then it has finite index since $\Gamma$ is just infinite. Then $N \lhd \finres\Gamma$ is non-trivial normal, then it is finite index because $\finres\Gamma$ is just infinite. Then $N = \finres\Gamma$ by the definition of $\finres\Gamma$.
\end{proof}

More generally, in the structure theory of just-infinite groups there is a trichotomy into residually finite virtually hereditarily just infinite, virtually simple, and Branch (residually finite but not virtually hereditarily just infinite).

\section{Lattices on products of trees}\label{sec:tree_products}

The groups in the main theorem are not residually finite because they contain subgroups that are not residually finite. These subgroups are irreducible lattices on products of trees and arise from the celebrated work of Burger--Mozes~\cite{BurgerMozes00} and Wise~\cite{Wise96} and the specific groups we will be working with are called Burger--Mozes--Wise (BMW) groups by Caprace and we follow this convention. A more detailed treatment can be found in \cite[Section~4]{Caprace}.

\begin{definition}
  Let $T_1$ and $T_2$ be regular trees of finite degrees $d_1$ and $d_2$. If $\Gamma < \Aut(T_1) \times \Aut(T_2)$ acts regularly (freely and transitively) on the vertices of $T_1 \times T_2$ then the action of $\Gamma$ on $T_1\times T_2$ is called \emph{BMW-action} and $\Gamma$ is called a \emph{BMW-group} of \emph{degree $(d_1,d_2)$}.

  A BMW-group is \emph{reducible} if it has a finite index subgroup that decomposes as a direct product or, equivalently, its image in $\Aut(T_1)$ or $\Aut(T_2)$ is discrete. Otherwise it is \emph{irreducible}.
\end{definition}

Given a BMW-action $\Gamma \curvearrowright T_1 \times T_2$ of degree $(d_1,d_2)$, we can equip each tree with a bicoloring $T_i \to \{1,2\}$. This yields a type function on the vertices of $T_1 \times T_2$  with values in $\{(1,1), (1,2), (2,1), (2,2)\}$. The subgroup $\Gamma^+$ of $\Gamma$ that preserves these types acts regularly on vertices of each type and the quotient group is the Klein four group $D_2$. The quotient complex $S \colon= \Gamma^+\backslash T_1 \times T_2$ is a square complex with four vertices and each vertex link is the complete bipartite graph of degree $(d_1, d_2)$. Furthermore there exists an action of $D_2$ on $S$ that is regular on the four vertices of $S$.

We call a square complex $S$ with complete bipartite vertex links and a vertex-regular $D_2$-action a \emph{BMW-complex}. The fundamental group $\Gamma^+ \defeq \pi_1(S)$ acts on the universal cover which is product of trees $\tilde{S} = T_1 \times T_2$. The deck transformations of $T_1 \times T_2$ that cover $D_2$ form an extension $\Gamma$ of $\Gamma^+$ by $D_2$ that is a BMW-group.
In what follows we will take both perspectives on BMW-groups: as groups acting on products of trees, and as square complexes with a $D_2$-action.

In order to compute presentations of BMW-groups we use the following special case of \cite[Theorem~1]{Brown}.
Let $X$ be a simply-connected CW-complex and assume that the 1-skeleton of $X$ does not contain loops. 
Assume that $G$ acts on $X$ permuting cells and regularly on vertices. Pick a vertex $x_0 \in X$ and let $E$ be the set of edges of $X$ incident with $x_0$. For $e\in E$ let $t(e,x_0)$ be the vertex of $e$, that is not $x_0$ and let $g_e \in G$ be the unique element that maps $x_0$ to $t(e,x_0)$. Let $F = F(g_e,e \in E)$ be the free group generated by the $g_e$.

Note that $g_e^{-1}$ maps $e \in E$ to an edge $f \in E$, possibly $e=f$. We define
      \[
        R_1 \defeq \{g_e g_f^{-1} \mid g_e^{-1} (e) = f\} \subseteq F.
      \]

If $\gamma = (e_1, \dots, e_k)$ is a combinatorial edge path starting at $x_0$ we define the word $W(\gamma) \in F$ recursively by
      \[
        W(\gamma) \defeq
        \begin{cases}
          g_{e} & \text{if } \gamma = (e), \\
          g_{e} W(g_e^{-1}(\gamma')) & \text{if } \gamma = (e)*\gamma'.\\
        \end{cases}
      \]
For a 2-cell $A$ incident with $x_0$ (pick an orientation for definiteness and) let $\gamma_A \in F$ be the edge loop starting at $x_0$ along the boundary of~$A$.

We define
      \[
        R_2 \defeq \{ W(\gamma_A) \mid A \text{ is a 2-cell incident with } x_0 \} \subseteq F.
      \]

\begin{proposition}\label{prop:brown}
  Let $X$ be a simply-connected CW-complex and assume that the 1-skeleton of $X$ does not contain loops. 
  Assume that $G$ acts on $X$ permuting cells and regularly on vertices and let $E$, $R_1$, and $R_2$ be as above.
  Then the homomorphism $F \to G$ that takes $g_e$ to $g_e$ descends to an isomorphism $\gen{g_e, e \in E \mid R_1 \cup R_2} \cong G$.
\end{proposition}

The following example of a BMW-group was studied by Radu~\cite{Radu_BMW}. It will be our source of non-residual-finiteness in the $q=2$ case of the main theorem.

\begin{figure}
  \centering
  \begin{align*}
    \includestandalone[page=1]{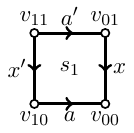}
        &&
        \includestandalone[page=2]{RaduSquareComplex}
        &&
        \includestandalone[page=3]{RaduSquareComplex}
        &&
        \includestandalone[page=4]{RaduSquareComplex}
        &&
        \includestandalone[page=5]{RaduSquareComplex}
        \\
        \includestandalone[page=6]{RaduSquareComplex}
        &&
        \includestandalone[page=7]{RaduSquareComplex}
        &&
        \includestandalone[page=8]{RaduSquareComplex}
        &&
        \includestandalone[page=9]{RaduSquareComplex}
  \end{align*}
  \caption{The square complex $\sradu$.}
  \label{fig:square_complex_R}
\end{figure}

\begin{example}\label{ex:SquareRadu}
  Consider the square complex $\sradu$ indicated in Figure \ref{fig:square_complex_R}. If we equip it with the action by $D_2 = \langle \sigma_A, \sigma_X \rangle$ described below, it becomes a BMW-complex.
  \begin{align*}
    \sigma_A:
        &&
    v_{00} &  \leftrightarrow v_{10},
           &
    v_{01} &  \leftrightarrow v_{11},
           &&&
    \sigma_X: &&
    v_{00} &  \leftrightarrow v_{01},
           &
    v_{10} &  \leftrightarrow v_{11},
    \\
           &&
    a & \leftrightarrow a^{-1},
      &
    b & \leftrightarrow b^{-1},
      &
    c & \leftrightarrow c^{-1},
      &&&
    a & \leftrightarrow a',
      &
    b & \leftrightarrow b',
      &
    c & \leftrightarrow c',
    \\
      &&
    a' & \leftrightarrow (a')^{-1},
       &
    b' & \leftrightarrow (b')^{-1},
       &
    c' & \leftrightarrow (c')^{-1},
       &&&
    x & \leftrightarrow  x^{-1},
      &
    y & \leftrightarrow y^{-1},
      &
    z & \leftrightarrow z^{-1},
    \\
      &&
    x & \leftrightarrow x',
      &
    y & \leftrightarrow y',
      &
    z & \leftrightarrow z'.
      &&&
    x' & \leftrightarrow  (x')^{-1},
       &
    y' & \leftrightarrow (y')^{-1},
       &
    z' & \leftrightarrow (z')^{-1}.
  \end{align*}
\end{example}

Let $\gradu$ denote the BMW-group arising from Example~\ref{ex:SquareRadu}, i.e. the extension of $\pi_1(\sradu)$ by $D_2 = \gen{\sigma_A,\sigma_X}$. Consider the universal cover $\tradu$ as homotopy classes of paths in $\sradu$ starting at $v_{00}$. Define the following automorphisms of $\tradu$ all of which lie in $\Gamma_{\Radu}$.
      \begin{align*}
        a &: [\gamma] \mapsto [a^{-1} * \sigma_A(\gamma)],
          &
        b &: [\gamma] \mapsto [b^{-1} * \sigma_A(\gamma)],
          &
        c &: [\gamma] \mapsto [c^{-1} * \sigma_A(\gamma)],
        \\
        x &: [\gamma] \mapsto [x^{-1} * \sigma_X(\gamma)],
          &
        y &: [\gamma] \mapsto [y^{-1} * \sigma_X(\gamma)],
          &
        z &: [\gamma] \mapsto [z^{-1} * \sigma_X(\gamma)].
      \end{align*}

\begin{proposition}[Radu]\label{prop:gamma33_pres}
  \begin{enumerate}
    \item
      The obvious homomorphism $F(a,b,c,x,y,z) \to \gradu$ descends to an isomorphism
      \begin{align*}
        \langle
        a,  b,  c,  x,  y,  z
        \mid
        a^2,  b^2,  c^2,  x^2,  y^2,  z^2,
        axax, ayay,  azbz,
        bxbx,  bycy,  cxcz
        \rangle \cong \gradu.
      \end{align*}
    \item\label{prop:gamma33_pres_nonrf}
      The BMW-action $(\Gamma_{\Radu}, \tradu)$ is irreducible and the group $\Gamma_{\Radu}$ is not residually finite.
    \item\label{prop:gamma33_pres_profiniteclosure}
      The element $x z$ lies in the profinite closure of $A \defeq \langle a, b, c\rangle$.
    \item\label{prop:gamma33_pres_explicitnrf}
      At least one of the following elements is in the finite residual of $\Gamma_{\Radu}$.
      \begin{align*}
        [ y  (x z)^2 y, x z],
                &&
                [ y  (x z)^2 y, x z b].
      \end{align*}
  \end{enumerate}
\end{proposition}
\begin{proof}
  The presentation and its geometric interpretation can be deduced from Proposition \ref{prop:brown}. For the other statements we refer to \cite[Proposition 5.4]{Radu_BMW}.
\end{proof}
Building on the previous proposition, we can prove the following.
\begin{lemma}\label{lem:shorter_element}
  The element $(x z)^4$ lies in the finite residual of $\Gamma_{\Radu}$.
  In particular it lies in the finite residual of $\pi_1(\sradu, v_{00})$ (represented by the path $(x^{-1}*z)^4$).
\end{lemma}
\begin{proof}
  Let $\delta \defeq xz$ and let $\phi$ be an epimorphism from $\Gamma$ to a finite group.
  We observe that $\delta a = b \delta$, $\delta 
  b = a\delta$ and $\delta c = c \delta^{-1}$. In particular we have $\delta^2 c = c \delta^{-2}$. Now conjugating this equation with $y$ yields
  \[
    y \delta^2 y b = b y \delta^{-2} y. \tag{$*$}
  \]
  Assume that $\phi(y \delta^2 y)$ and $\phi(\delta)$ commute.        
  Then we can conjugate $(*)$ with $\phi(\delta)$ and $\phi(y)$ and we obtain
  \begin{align*}
    \phi(\delta y \delta^2 y \delta^{-1} a) &= \phi(a \delta y \delta^{-2} y \delta^{-1}),\\
    \phi(y \delta^2 y a) &= \phi(a y \delta^{-2} y),\\
    \phi(\delta^2 a) &= \phi(a \delta^{-2}),\\
    \phi(\delta^2) &= \phi(\delta^{-2}).\\
  \end{align*}
  In particular $\phi(\delta^4) = 1$ in that case. If $\phi(y \delta^2 y)$ and $\phi(\delta b)$ commute, then we conjugate $(*)$ with $\phi(\delta b)$ and $\phi(y)$ and we obtain
  \begin{align*}
    \phi(\delta b y \delta^2 y b \delta^{-1} a) &= \phi(a \delta b \delta y \delta^{-2} y b \delta^{-1}),\\
    \phi(y \delta^2 y a) &= \phi(a y \delta^{-2} y),\\
    \phi(\delta^2) &= \phi(\delta^{-2}).\\
  \end{align*}
  Since by Proposition \ref{prop:gamma33_pres}.\ref{prop:gamma33_pres_explicitnrf} at least one of these two cases occurs, we deduce that $\delta^4$ lies in the finite residual of~$\Gamma_{\Radu}$.
\end{proof}

\begin{remark}\label{rem:rep_gamma_r}
  \begin{enumerate}
    \item The previous lemma also implies that $[y (x z)^2 y, x z]$ lies in the finite residual of $\Gamma_{\Radu}$. Indeed since $\delta^4$ vanishes in every finite quotient, we get that $\phi(y (x z)^2 y)$ centralizes $\phi(A)$ for every map $\phi$ to a finite group. Since $xz$ lies in the profinite closure of $A$, we get that $\phi \left([y (x z)^2 y, x z]\right) = 1$.
    \item The elements $\delta^4$ and $\delta^{-4}$ are the shortest elements in the finite residual of $\Gamma_{\Radu}$. This can be seen as follows.
    Let $D$ be the derived subgroup of $\Gamma_{\Radu}$ which has index 8 and is generated by the words $abxyxy$ and $ba$. Let $\alpha \defeq \sqrt{-15}$ and $\beta \defeq \sqrt{17}$ then the following assignment extends to an homomorphism $\Phi:\Gamma \to \SL_2(\QQ(\alpha, \beta))$. 
    \begin{align*}
        abxyxy &\mapsto
        \begin{pmatrix}
            -\frac 1 4 - \frac \alpha 4 & 0\\
            0 & -\frac 1 4 + \frac \alpha 4
        \end{pmatrix}
        &
        ba &\mapsto
        \begin{pmatrix}
            -\frac 1 8 - \frac \alpha {24} - \frac \beta 8 + \frac {\alpha \beta}{40} & 1\\
            - \frac 8 {15} & -\frac 1 8 + \frac \alpha {24} - \frac \beta 8 - \frac {\alpha \beta}{40}
        \end{pmatrix}
    \end{align*}
    The image of $\Phi$ is residually finite, so an element in $\Gamma_{\Radu}^{(\infty)}$ must vanish under $\Phi$.
    We checked with a computer that the only non-trivial elements in $\Gamma_{\Radu}$ of length at most eight that vanish under $\Phi$ are $\delta^4$ and $\delta^{-4}$.  

    We also briefly describe how we construct the homomorphism $\Phi$. Using a computer, we observe that $D$ admits many finite $\SL_2$-quotients in different characteristics, and we suspected that they arise as congruence quotients.
    We choose a particular epimorphism $\Phi_{13}: D \to \SL_2(\FF_{13})$ and Hensel-lift it to sufficiently high precision $\Phi_{13^k}: D \to \SL_2(\ZZ/{13^k}\ZZ)$. 
    We consider the traces of $\Phi_{13^k}(abxyxy)$, $\Phi_{13^k}(ba)$ and $\Phi_{13^k}(xyxy)$ and check if they are roots of relative simple polynomials. They are, so we have algebraic numbers as candidates for the traces.
    Assuming there exists a homomorphism $\Phi: D \to \SL_2(\KK)$, where~$K$ is a number field, we denote the images of $abxyxy$ and $ba$ by $A$ and $B$ and assume that $\Tr(A), \Tr(B)$ and $\Tr(BA)$ are as we suspect.
    We know the characteristic polynomial of~$A$ and assume $A$ to be a diagonal. 
    We have two linear equations in the diagonal entries of~$B$, namely $B_{11} + B_{22} = \Tr(B)$ and
    $A_{11}B_{11} + A_{22}B_{22} = \Tr(BA)$. We solve for $B_{11}$ and $B_{22}$. If $B_{11}B_{22}$ is not 1, then $B_{12}$ and $B_{21}$ are both not 0 and we asssume that this is the case.
    Since conjugating with a diagonal matrix preserves $A$ and the diagonal entries of $B$, we may assume that $B_{12} = 1$. By the determinant equation, the last entry of $B$ is $B_{21} = B_{11}B_{22} -1$.
  \end{enumerate}
\end{remark}

The next example is due to Janzen--Wise \cite{JanzenWise} and will be our source of non-residual-finiteness in the case $q = 3$ of the main theorem.

\begin{example}\label{ex:SquareJanzenWise}
  Consider the square complex $\sjw$ indicated in Figure \ref{fig:square_complex_JW}. If we equip it with the action by $D_4 = \langle \sigma_A, \sigma_X \rangle$ described below, it becomes a BMW-complex.
  \begin{align*}
    \sigma_A\colon
        &&
    v_{00} &  \leftrightarrow v_{10},
           &
    v_{10} &  \leftrightarrow v_{11},
           &
    \sigma_X\colon &&
    v_{00} &  \leftrightarrow v_{01},
           &
    v_{01} &  \leftrightarrow v_{11},
    \\
           &&
    a & \leftrightarrow {\bar a}^{-1},
      &
    b & \leftrightarrow {\bar b}^{-1},
      &&&
    a & \leftrightarrow a',
      &
    b & \leftrightarrow b',\\
      &&
    a' & \leftrightarrow ({\bar a'})^{-1},
       &
    b' & \leftrightarrow ({\bar b'})^{-1},
       &&&
    \bar a & \leftrightarrow \bar a',
           &
    \bar b & \leftrightarrow \bar b',
    \\
           &&
    x & \leftrightarrow x',
      &
    y & \leftrightarrow y',
      &&&
    x & \leftrightarrow \bar x^{-1},
      &
    y & \leftrightarrow \bar y^{-1},\\
      &&
    \bar x & \leftrightarrow \bar x',
           &
    \bar y & \leftrightarrow \bar y'.
           &&&
    x' & \leftrightarrow (\bar x') ^{-1},
       &
    y' & \leftrightarrow (\bar y') ^{-1}.
  \end{align*}

  \begin{figure}
    \centering
    \begin{align*}
      \includestandalone[page=1]{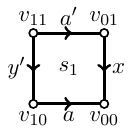}
        &&
        \includestandalone[page=2]{JanzenWiseSquareComplex}
        &&
        \includestandalone[page=3]{JanzenWiseSquareComplex}
        &&
        \includestandalone[page=4]{JanzenWiseSquareComplex}
        &&
        \includestandalone[page=5]{JanzenWiseSquareComplex}
        &&
        \includestandalone[page=6]{JanzenWiseSquareComplex}
        \\
        \includestandalone[page=7]{JanzenWiseSquareComplex}
        &&
        \includestandalone[page=8]{JanzenWiseSquareComplex}
        &&
        \includestandalone[page=9]{JanzenWiseSquareComplex}
        &&
        \includestandalone[page=10]{JanzenWiseSquareComplex}
        &&
        \includestandalone[page=11]{JanzenWiseSquareComplex}
        &&
        \includestandalone[page=12]{JanzenWiseSquareComplex}
        \\
        \includestandalone[page=13]{JanzenWiseSquareComplex}
        &&
        \includestandalone[page=14]{JanzenWiseSquareComplex}
        &&
        \includestandalone[page=15]{JanzenWiseSquareComplex}
        &&
        \includestandalone[page=16]{JanzenWiseSquareComplex}
    \end{align*}
    \caption{The square complex $\sjw$.}
    \label{fig:square_complex_JW}
  \end{figure}

    Let $\tjw$ be the universal cover $\sjw$. Consider the universal cover $\tjw$ as homotopy classes of paths in $\sjw$ starting at $v_{00}$.
    Let $\Gamma_{\JW} \leq \Aut(\tjw)$ be the extension of $\pi_1(\sjw)$ covering $\langle \sigma_A, \sigma_X\rangle$.
  Define the following automorphisms of $\tjw$ in $\gjw$.
        \begin{align*}
          a &: [\gamma] \mapsto [a^{-1} * \sigma_A(\gamma)],
            &
          b &: [\gamma] \mapsto [b^{-1} * \sigma_A(\gamma)],
          \\
          x &: [\gamma] \mapsto [x^{-1} * \sigma_X(\gamma)],
            &
          y &: [\gamma] \mapsto [y^{-1} * \sigma_X(\gamma)].
        \end{align*}
  \begin{proposition}[Janzen--Wise, Caprace]\label{prop:JanzenWise}
    \begin{enumerate}
      \item
      The obvious homomorphism $F(a,b,x,y) \to \Gamma_{\JW}$ descends to an isomorphism
        \begin{align*}
          \Gamma_{\JW} \cong \langle
          a,  b,  x,  y
          \mid
          axay, ax^{-1}by^{-1}, ay^{-1}b^{-1}x^{-1}, bxb^{-1}y^{-1}
          \rangle.
        \end{align*}
      \item
        The BMW-action $(\Gamma_{\JW}, S_{\JW})$ is irreducible and the group $\Gamma_{\JW}$ is not residually finite.
      \item
        The elements $[x^3, y^4]$ and $[y^3, x^4]$ are both in the finite residual of $\Gamma_{\JW}$.
        In particular the homotopy classes of the following loops are contained in the finite residual of $\pi_1(\mathcal S_{\JW}, v_{00})$.
        \begin{align*}
          x^{-1} * \bar x * x^{-1} * \bar y * y^{-1} * \bar y * y^{-1} * x * \bar x^{-1} * x * \bar y^{-1} * y * \bar y^{-1} * y, \\
          y^{-1} * \bar y * y^{-1} * \bar x * x^{-1} * \bar x * x^{-1} * y * \bar y^{-1} * y * \bar x^{-1} * x * \bar x^{-1} * x.
        \end{align*}
    \end{enumerate}
  \end{proposition}

  \begin{proof}
    Again the presentation can be deduced from Proposition \ref{prop:brown}. The irreducibility has been established by Janzen--Wise \cite{JanzenWise}. The computation of the explicit elements in the finite residual is due to Caprace \cite[Remark 4.20]{Caprace}.
  \end{proof}

\end{example}

\section{Non-residually finite $\tilde C_2$-lattices}\label{sec:non-rf_c2tilde}

\begin{figure}
  \centering
  \begin{align*}
    \includestandalone[page=1]{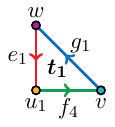}
            &&
            \includestandalone[page=2]{TriangleComplexT}
            &&
            \includestandalone[page=3]{TriangleComplexT}
            &&
            \includestandalone[page=4]{TriangleComplexT}
            &&
            \includestandalone[page=5]{TriangleComplexT}
            &&
            \includestandalone[page=6]{TriangleComplexT}
            \\[-2mm]
            \includestandalone[page=7]{TriangleComplexT}
            &&
            \includestandalone[page=8]{TriangleComplexT}
            &&
            \includestandalone[page=9]{TriangleComplexT}
            &&
            \includestandalone[page=10]{TriangleComplexT}
            &&
            \includestandalone[page=11]{TriangleComplexT}
            &&
            \includestandalone[page=12]{TriangleComplexT}
            \\[-2mm]
            \includestandalone[page=13]{TriangleComplexT}
            &&
            \includestandalone[page=14]{TriangleComplexT}
            &&
            \includestandalone[page=15]{TriangleComplexT}
            &&
            \includestandalone[page=16]{TriangleComplexT}
            &&
            \includestandalone[page=17]{TriangleComplexT}
            &&
            \includestandalone[page=18]{TriangleComplexT}
            \\[-2mm]
            \includestandalone[page=19]{TriangleComplexT}
            &&
            \includestandalone[page=20]{TriangleComplexT}
            &&
            \includestandalone[page=21]{TriangleComplexT}
            &&
            \includestandalone[page=22]{TriangleComplexT}
            &&
            \includestandalone[page=23]{TriangleComplexT}
            &&
            \includestandalone[page=24]{TriangleComplexT}
            \\[-2mm]
            \includestandalone[page=25]{TriangleComplexT}
            &&
            \includestandalone[page=26]{TriangleComplexT}
            &&
            \includestandalone[page=27]{TriangleComplexT}
            &&
            \includestandalone[page=28]{TriangleComplexT}
            &&
            \includestandalone[page=29]{TriangleComplexT}
            &&
            \includestandalone[page=30]{TriangleComplexT}
            \\[-2mm]
            \includestandalone[page=31]{TriangleComplexT}
            &&
            \includestandalone[page=32]{TriangleComplexT}
            &&
            \includestandalone[page=33]{TriangleComplexT}
            &&
            \includestandalone[page=34]{TriangleComplexT}
            &&
            \includestandalone[page=35]{TriangleComplexT}
            &&
            \includestandalone[page=36]{TriangleComplexT}
            \\[-2mm]
            \includestandalone[page=37]{TriangleComplexT}
            &&
            \includestandalone[page=38]{TriangleComplexT}
            &&
            \includestandalone[page=39]{TriangleComplexT}
            &&
            \includestandalone[page=40]{TriangleComplexT}
            &&
            \includestandalone[page=41]{TriangleComplexT}
            &&
            \includestandalone[page=42]{TriangleComplexT}
            \\[-2mm]
            \includestandalone[page=43]{TriangleComplexT}
            &&
            \includestandalone[page=44]{TriangleComplexT}
            &&
            \includestandalone[page=45]{TriangleComplexT}
  \end{align*}
  \caption{The GAB $Y^2_1$, the sub complex generated by the triangles $\boldsymbol{t_1}, \dots, \boldsymbol{t_{18}}$ is isomorphic to a subdivision of $\sradu$.}
  \label{fig:TriangleComplexT}
\end{figure}

\begin{figure}
  \centering
  \begin{align*}
    \includestandalone[page=1]{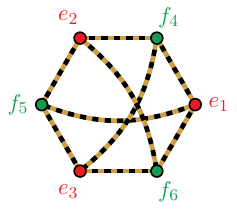}
            &&
            \includestandalone[page=2]{RaduGraph}
            &&
            \includestandalone[page=3]{RaduGraph}
            \\
            \includestandalone[page=4]{RaduGraph}
            &&
            \includestandalone[page=5]{RaduGraph}
  \end{align*}
  \begin{align*}
    \includestandalone[page=6]{RaduGraph}
            &&
            \includestandalone[page=7]{RaduGraph}
  \end{align*}
  \caption{The links in $Y_1^2$ of the vertices $u_1$, $u_2$, $u_3$, $u_4$, $u_5$, $v$, and $w$.}
  \label{fig:perfect_c2tilde}
\end{figure}
In this section we present the complexes $Y_i^q$ from the main theorem. The complexes have been found with computer assistance and we present our approach to find these complexes in Section~\ref{sec:searching}.

We begin with our example for $q=2$.

\begin{theorem}\label{thm:main_2}
  Let $Y_1^2$ be the chamber complex indicated in Figure~\ref{fig:TriangleComplexT}. Its universal cover $X_1^2 = \tilde{Y}_1^2$ is a $\tilde C_2$-building of thickness $3$. Its fundamental group $\Gamma_1^2 = \pi_1(X_1^2)$ is not residually finite. In fact, $\Gamma_1^2$ does not have any finite index subgroups.
\end{theorem}

All assertions are straightforward to verify, except for the last statement. 
The proof for the claim that $\Gamma_1^2$ has no finite-index subgroups is postponed to the end of Section~\ref{sec:normal_forms}; the remaining statements are proved below.

\begin{proof}[Proof of Theorem \ref{thm:main_2}, Part 1]
By considering the links of $Y_1^2$, indicated in Figure~\ref{fig:perfect_c2tilde}, we immediately conclude that it is a GAB of type~$\tilde{C}_2$.
Therefore, $X_1^2$ is a building of type~$\tilde{C}_2$ by Theorem~\ref{thm:building_cartan-hadamard}.
If we subdivide the complex $\sqcomp_{\Radu}$ along the diagonals from~$v_{00}$ to $v_{11}$ and label the new diagonals with $s_i$, we obtain a $\tilde C_2$-GAB that embeds in $Y_k^2$ via the following assignment.
\begin{align*}
    v_{00} &\mapsto v, &
    v_{11} &\mapsto w, &
    v_{10} &\mapsto u_1, &
    v_{01} &\mapsto u_2, &
    a & \mapsto f_4, &
    b & \mapsto f_5, &
    c & \mapsto f_6, \\
    x & \mapsto f_1, &
    y & \mapsto f_2, &
    z & \mapsto f_3, &
    a' & \mapsto e_4, &
    b' & \mapsto e_5, &
    c' & \mapsto e_6, &
    x' & \mapsto e_1, \\
    y' & \mapsto e_2, &
    z' & \mapsto e_3, &
    s_i & \mapsto g_i.
\end{align*}

In particular, $\pi_1(\sqcomp_{\Radu})$ embeds into $\Gamma_1^2$ by Lemma~\ref{lem:locally_convex}, showing that $\Gamma_1^2$ is not residually finite by Proposition~\ref{prop:gamma33_pres}.2. 

\end{proof}

\subsection{Geometric presentations}\label{sec:geometric_presentations}

    Our next goal is to use Proposition~\ref{prop:brown} to derive presentations for $\tilde{C}_2$-lattices that act regularly on the special vertices.

    Let $Y$ be a $\tilde C_2$-GAB equipped with two special vertices $v$ and $w$, one of each special type and assume that $Y$ admits an involutory automorphism $\rho$ that interchanges $v$ and $w$. Let $X = \tilde{Y}$ be the universal cover and let $\Gamma$ be the group of transformations covering $\gen{\rho}$, so that $\Gamma$ acts regularly on special vertices and contains $\pi_1(Y)$ with index $2$.

    We cannot apply Proposition \ref{prop:brown} directly, so we modify $X$ as follows.
    When we remove from $X$ the non-special vertices together with their stars, we are also removing all short edges. Thus we are left with the full subgraph of special vertices and long edges, which we denote~$X^{(s)}$. In $X^{(s)}$ there are many four-cycles, all of which are apartments in links of non-special vertices of~$X$. We call such a four-cycle an \emph{$(A_1 \times A_1)$-boundary} and construct a complex $X'$ by filling in all $(A_1 \times A_1)$-boundaries in $X^{(s)}$.    
    Note that $\Gamma$ acts regularly on the vertices of $X'$. 

    Now let $G$ be the set of edges in $Y$ from $v$ to $w$ (oriented in this way) and call edge paths in $Y$ that lift to an $(A_1 \times A_1)$-boundary also $(A_1 \times A_1)$-boundaries. Consider the following words in $F(G)$
    \begin{align*}
        R_1 &:= \{g h \mid g \in G \text{ and }\rho(g)^{-1} = h\},
        \\
        R_2 &:= \{g_1 g_2^{-1} g_3 g_4^{-1} \mid g_1 * g_2^{-1} * g_3 * g_4^{-1} \text{ is an $(A_1\times A_1)$-boundary} \}.
    \end{align*}

    \begin{lemma}\label{lem:regular_pres}
    In the above situation the space $X'$ is simply connected and $\Gamma \cong \langle G \mid R_1 \cup R_2 \rangle$.
    More precisely, if we consider $X$ as homotopy classes of paths in $Y$ starting  at $v$ then a generator~$g$ of this presentation acts on $X$ as follows
    \[
        [\gamma] \mapsto [g * \rho(\gamma)].
    \]
    In particular, $g \mapsto g.v$ is a $\Gamma$-equivariant isomorphism $\Cay(\Gamma,G) \to X^{(s)}$.
    \end{lemma}
    \begin{figure}
        \centering
            \begin{tikzpicture}[baseline=(current bounding box.center)]
                \def\t{1.5}
                \node[circle, scale = 0.1] (A) at (0,0) {};
                \node[circle, scale = 0.1] (B) at (\t,0) {};
                \node[circle, scale = 0.1] (C) at (0,\t) {};
                \node[circle, scale = 0.1] (D) at (0,-\t) {};
                \node[circle, scale = 0.1] (E) at (-\t,0) {}; 
                \fill[pattern = north east lines, pattern color = lightgray!75] (A.center) -- (B.center) -- (C.center);
                \fill[pattern = north west lines, pattern color = lightgray!75] (A.center) -- (B.center) -- (D.center);
                \fill[pattern = north west lines, pattern color = lightgray!75] (A.center) -- (C.center) -- (E.center);
                \fill[pattern = north east lines, pattern color = lightgray!75] (A.center) -- (D.center) -- (E.center);
                \draw[line width = 0.5mm, color = ForestGreen, ->-] (A) -- (B);
                \draw[line width = 0.5mm, color = RoyalBlue, ->-] (B) -- (C);
                \draw[line width = 0.5mm, color = Red, ->-] (C) -- (A);
                \draw[line width = 0.5mm, color = ForestGreen, ->-] (A) -- (E);
                \draw[line width = 0.5mm, color = Red, ->-] (D) -- (A);
                \draw[line width = 0.5mm, color = RoyalBlue, ->-] (B) -- (D);
                \draw[line width = 0.5mm, color = RoyalBlue, ->-] (E) -- (D);
                \draw[line width = 0.5mm, color = RoyalBlue, ->-] (E) -- (C);
                \draw[fill = Dandelion, thick] (A) circle (0.7mm);
                \draw[fill = BlueGreen, thick] (B) circle (0.7mm);
                \draw[fill = Mulberry, thick] (C) circle (0.7mm);
                \draw[fill = Mulberry, thick] (D) circle (0.7mm);
                \draw[fill = BlueGreen, thick] (E) circle (0.7mm);
                \node[above right] at (A) {$u$};
                \node[right] at (B) {$v$};
                \node[above] at (C) {$w$};
                \node[below] at (D) {$w$};
                \node[left] at (E) {$v$};
                \node[below] at ($(A)!.5!(B)$) {$f_1$};
                \node[right] at ($(A)!.5!(C)$) {$e_1$};
                \node[left] at ($(A)!.5!(D)$) {$e_2$};
                \node[above] at ($(A)!.5!(E)$) {$f_2$};
                \node (F) at (0.69*\t, 0.69*\t) {$g_1$};
                \node (G) at (-0.69*\t, 0.69*\t) {$g_2$};
                \node (H) at (-0.69*\t, -0.69*\t) {$g_3$};
                \node (I) at (0.69*\t, -0.69*\t) {$g_4$};
            \end{tikzpicture}
        \caption{The edge path $(g_1 * g_2^{-1} * g_3 * g_4^{-1})$ is an $(A_1\times A_1)$-boundary.}
        \label{fig:a1a1_boundary}
    \end{figure}
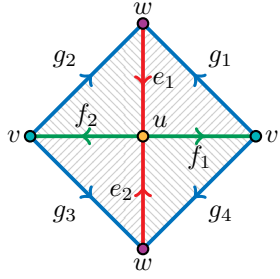
    \begin{proof}
      If $u \in X$ is a non-special vertex and $g$ is an edge in its link (a complete bipartite graph), then $g$ together with edges that meet it in a vertex form a spanning tree $T$ of $\lk u$ and $\pi_1(\lk u) \cong \pi_1(\lk u/T,T)$ is freely generated by edges opposite $g$ (this is a very special case of the Solomon--Tits theorem). It follows that gluing in all four-cycles that contain $g$ produces a simply-connected space. Gluing in all four-cycles rather than only the ones that contain $g$ certainly results in a simply connected space. Since $X$ is simply connected, it follows that removing the star of a non-special vertex and gluing in all four-cycles in its link produces a simply connected space (formally applying Seifert--van Kampen to the diagram
      \[
        \pi_1(\{\text{union of four-cycles}\}) \leftarrow \pi_1(\lk u) \rightarrow \pi_1(X \setminus \st u),
      \]
      where the arrow to the right is an isomorphism and the arrow to the left is trivial).
      Thus $X'$ is simply connected and we can apply Proposition \ref{prop:brown}.
    \end{proof}
    
    Applying Lemma~\ref{lem:regular_pres} to the lattice $\bar{\Gamma}^2_1$ gives:

    \begin{proposition}\label{prop:presentation_bar_2}
        \begin{enumerate}
          \item The lattice $\bar{\Gamma}^2_1$ is presented by generators $g_1,\ldots,g_{15}$ subject to the relations
            \begin{align*}
                &
                g_1 ^2, && g_2^2, && g_3g_6, && g_4^2, && g_5g_8, && g_7g_9, \\
                &
                g_{10}^2 && g_{11}^2 && g_{12}g_{14}, && g_{13}^2, && g_{15}^2,&&
                g_1  g_4     g_5   g_2, &&\\
                &
                g_1  g_4     g_6   g_6, &&
                g_1  g_7     g_7   g_6, && 
                g_1  g_7     g_8   g_2, && 
                g_3  g_{11}  g_{12}  g_{13},&&
                g_3  g_{12}  g_{4}   g_{13},&&
                g_3  g_{12}  g_{13}  g_{15},\\
                &
                g_2  g_{14}  g_{7}   g_{14},&&
                g_7  g_{10}  g_{14}  g_{14},&&
                g_7  g_{14}  g_{14}  g_{15},&&
                g_1  g_{10}  g_{5}   g_{11},&&
                g_5  g_{11}  g_{10}  g_{15},&&
                g_5  g_{13}  g_{11}  g_{10}.
            \end{align*}
            \item If we regard $X^2_1$ as homotopy classes of paths in $Y^2_1$ starting at $v$, then the generator $g_i$ acts on $X^2_1$ via
            \[
                [\gamma] \mapsto [g_i * \rho(\gamma)].
            \]
            In particular, the full subgraph of the 1-skeleton of $X^2_1$ on the set of special vertices is $\bar{\Gamma}^2_1$-equivariantly identified with the Cayley graph of $\bar{\Gamma}^2_1$ with respect to the generating set $\{g_1, \dots, g_{15}\}$.
        \end{enumerate}
    \end{proposition}

    \newcommand{\boundary}[4]{g_{#1} g_{#2}^{-1} g_{#3} g_{#4}^{-1}}

\begin{figure}
        \centering
        \begin{align*}
            \includestandalone[page=1]{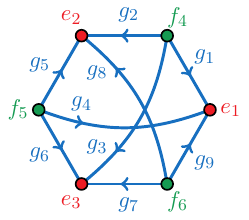}
            &&
            \includestandalone[page=2]{boundaries}
            \\
            \includestandalone[page=3]{boundaries}
            &&
            \includestandalone[page=4]{boundaries}
        \end{align*}
        \caption{The $(A_1\times A_1)$-boundaries around $u_1, u_3, u_4$ and $u_5$ are cycles of length four in the indicated graphs.}
        \label{fig:boundary_graphs}
    \end{figure}

    \begin{proof}
        We apply Lemma \ref{lem:regular_pres} and obtain the following relations over the generating set $(g_i)^{\pm 1}$.
        \begin{align*}
            R_1 &=
            \{g_1^2, \; g_2^2, \; g_3g_6, \; g_4^2, \; g_5g_8, \; g_7g_9,
            \; g_{10}^2, \;g_{11}^2, \; g_{12}g_{14}, \; g_{15}^2\},\\
            R_2 &= \{\text{$(A_1\times A_1)$-boundary relations} \}.
        \end{align*}
        (Up to cyclic permutation) there are 45 relations in $R_2$, nine for each vertex of non-special type in~$\mathcal T$.
        Since $\rho$ swaps $u_1$ and~$u_2$ the relations arising from $(A_1\times A_1)$-boundaries around these vertices are equivalent. 
        We now list the 36 relations relations arising from boundaries around $u_1, u_3, u_4, u_5$. Note that these can be read out easily from Figure \ref{fig:boundary_graphs}.
        \begin{align*}
            r_{1,1} & = \boundary{1}{4}{5}{2},
            &
            r_{1,2} & = \boundary{1}{4}{6}{3},
            &
            r_{1,3} & = \boundary{1}{9}{7}{3},
            &
            r_{1,4} & = \boundary{1}{9}{8}{2},
            \\
            r_{1,5} & = \boundary{2}{5}{6}{3},
            &
            r_{1,6} & = \boundary{2}{8}{7}{3},
            &
            r_{1,7} & = \boundary{4}{9}{8}{5},
            &
            r_{1,8} & = \boundary{4}{9}{7}{6},
            \\
            r_{1,9} & = \boundary{5}{8}{7}{6},
            &
            r_{3,1} &= \boundary{3}{11}{6}{15},
            &
            r_{3,2} &= \boundary{3}{11}{12}{13},
            &
            r_{3,3} &= \boundary{3}{14}{4}{13},
            \\
            r_{3,4} &= \boundary{3}{14}{13}{15},
            &
            r_{3,5} &= \boundary{4}{12}{6}{13},
            &
            r_{3,6} &= \boundary{4}{12}{11}{14},
            &
            r_{3,7} &= \boundary{4}{13}{15}{13},
            \\
            r_{3,8} &= \boundary{6}{13}{14}{11},
            &
            r_{3,9} &= \boundary{6}{15}{13}{12},
            &
            r_{4,1} &= \boundary{2}{12}{7}{12},
            &
            r_{4,2} &= \boundary{2}{12}{15}{14},
            \\
            r_{4,3} &= \boundary{2}{14}{9}{14},
            &
            r_{4,4} &= \boundary{2}{14}{10}{12},
            &
            r_{4,5} &= \boundary{7}{10}{9}{15},
            &
            r_{4,6} &= \boundary{7}{10}{14}{12},
            \\
            r_{4,7} &= \boundary{7}{12}{14}{15},
            &
            r_{4,8} &= \boundary{9}{14}{12}{10},
            &
            r_{4,9} &= \boundary{9}{15}{12}{14},
            &
            r_{5,1} &= \boundary{1}{10}{5}{11},
            \\
            r_{5,2} &= \boundary{1}{10}{15}{10},
            &
            r_{5,3} &= \boundary{1}{11}{8}{10},
            &
            r_{5,4} &= \boundary{1}{11}{13}{11},
            &
            r_{5,5} &= \boundary{5}{11}{10}{15},
            \\
            r_{5,6} &= \boundary{5}{13}{8}{15},
            &
            r_{5,7} &= \boundary{5}{13}{11}{10},
            &
            r_{5,8} &= \boundary{8}{10}{11}{13},
            &
            r_{5,9} &= \boundary{8}{15}{10}{11}
        \end{align*}
        The relations in the presentation in the theorem, which are not in $R_1$ are equivalent to the relations $r_{1,1}, r_{1,2}, r_{1,3}, r_{1,4}, r_{3,2}, r_{3,3}, r_{3,4}, r_{4,1}, r_{4,6}, r_{4,7}, r_{5,1}, r_{5,5}$ and $r_{5,7}$.
        The reaming relations can be deduced as follows.
        \begin{align*}
            r_{1,1} \land r_{1,2} &\leadsto r_{1,5},
            &
            r_{1,3} \land r_{1,4} &\leadsto r_{1,6},
            &
            r_{1,1} \land r_{1,4} &\leadsto r_{1,7},
            &
            r_{1,2} \land r_{1,3} &\leadsto r_{1,8},
            &
            r_{1,7} \land r_{1,8} &\leadsto r_{1,9},
            \\
            r_{3,2} &\leadsto r_{3,8},
            &
            r_{3,3} &\leadsto r_{3,5},
            &
            r_{3,4} &\leadsto r_{3,9},
            &
            r_{3,2} \land r_{3,3} &\leadsto r_{3,6},
            &
            r_{3,3} \land r_{3,4} &\leadsto r_{3,7},
            \\            
            r_{3,8} \land r_{3,9} &\leadsto r_{3,1},
            &
            r_{4,1} &\leadsto r_{4,3},
            &
            r_{4,6} &\leadsto r_{4,8},
            &
            r_{4,7} &\leadsto r_{4,9},
            &
            r_{4,1} \land r_{4,6} &\leadsto r_{4,4},
            \\
            r_{4,1} \land r_{4,7} &\leadsto r_{4,2},
            &
            r_{4,8} \land r_{4,9} &\leadsto r_{4,5},
            &
            r_{5,1} &\leadsto r_{5,3},
            &
            r_{5,5} &\leadsto r_{5,9},
            &
            r_{5,7} &\leadsto r_{5,8},
            \\
            r_{5,1} \land r_{5,5} &\leadsto r_{5,2},
            &
            r_{5,3} \land r_{5,8} &\leadsto r_{5,4},
            &
            r_{5,8} \land r_{5,9} &\leadsto r_{5,6}.
        \end{align*}
        Lemma \ref{lem:regular_pres} also yields, that the generators act on $\tilde {\mathcal T}$ as described in the theorem.
    \end{proof}
    
    Applying the Reidemeister--Schreier procedure to the presentation of $\bar{\Gamma}^2_1$ we obtain a presentation for $\Gamma^2_1$:
    \begin{corollary}\label{cor:presentation_2}
      \begin{enumerate}
        \item
      The lattice $\Gamma^2_1$ is presented by generators $h_2,\ldots,h_{15}$ subject to the relations
      \begin{align*}
                &
                h_4 h_8^{-1} h_2,
                &&
                h_4^{-1} h_5 h_2^{-1},
                &&
                h_4 h_3^{-1} h_6,
                &&
                h_4^{-1} h_6 h_3^{-1},
                \\
                &
                h_7 h_9^{-1} h_6,
                &&
                h_9^{-1} h_7 h_3^{-1},
                && 
                h_7 h_5^{-1} h_2,
                &&
                h_9^{-1} h_7 h_3^{-1},
                \\
                &
                h_3 h_{11}^{-1} h_{12} h_{13}^{-1},
                &&
                h_6^{-1} h_{11} h_{14}^{-1} h_{13},
                &&
                h_3 h_{14}^{-1} h_4 h_{13}^{-1},
                &&
                h_6^{-1} h_{12} h_4^{-1} h_{13},
                \\
                &
                h_3 h_{14}^{-1} h_{13} h_{15}^{-1},
                &&
                h_6^{-1} h_{12} h_{13}^{-1} h_{15},
                &&
                h_2 h_{12}^{-1} h_7 h_{12}^{-1},
                &&
                h_2^{-1} h_{14} h_9^{-1} h_{14},
                \\
                &
                h_7 h_{10}^{-1} h_{14} h_{12}^{-1},
                &&
                h_9^{-1} h_{10} h_{12}^{-1} h_{14},
                &&
                h_7 h_{12}^{-1} h_{14} h_{15}^{-1},
                &&
                h_9^{-1} h_{14} h_{12}^{-1} h_{15},
                \\
                &
                h_{10} h_8^{-1} h_{11},
                &&
                h_{10}^{-1} h_5 h_{11}^{-1},
                &&
                h_5 h_{11}^{-1} h_{10} h_{15}^{-1},
                &&
                h_8^{-1} h_{11} h_{10}^{-1} h_{15},
                \\
                &
                h_5 h_{13}^{-1} h_{11} h_{10}^{-1},
                &&
                h_8^{-1} h_{13} h_{11}^{-1} h_{10}.
                \rangle.
            \end{align*}
          \item The inclusion $\iota:\Gamma^2_1 \to \bar{\Gamma}^2_1$ is given by $h_i\mapsto g_1 g_i$.
          \item Identify $\Gamma_1^2$ as a subgroup of $\bar \Gamma_1^2$ via $\iota$. Then conjugation by $g_1$ in $\bar \Gamma_2^1$ acts on $\Gamma_2^1$ as follows.
            \begin{align*}
                h_2 &\leftrightarrow h_2^{-1},
                &
                h_3 &\leftrightarrow h_6^{-1},
                &
                h_4 &\leftrightarrow h_4^{-1},
                &
                h_5 &\leftrightarrow h_8^{-1},
                &
                h_7 &\leftrightarrow h_9^{-1},
                \\
                h_{10} &\leftrightarrow h_{10}^{-1},
                &
                h_{11} &\leftrightarrow h_{11}^{-1},
                &
                h_{12} &\leftrightarrow h_{14}^{-1},
                &
                h_{13} &\leftrightarrow h_{13}^{-1},
                &
                h_{15} &\leftrightarrow h_{15}^{-1}.
            \end{align*}
        \end{enumerate}
      \end{corollary}

    \subsection{Normal forms}\label{sec:normal_forms}
    Any uniform lattice on a locally finite building is biautomatic by \cite[Theorem~6.7]{Swiatkowski} and the main theorem of \cite{Osjada}. However, extracting explicit automatic structures from this proof is not immediate. In this section we develop an easy algorithm to compute normal forms for presentation arising from Lemma~\ref{lem:regular_pres}.
    We will use it in order to prove that $\Gamma_2^1$ contains no proper finite index subgroups and to compute the full automorphism groups of the buildings $X_i^q$. Moreover, the algorithm allows to perform computations in the group algebras of the lattices $\bar \Gamma_i^q$ which was used to establish property (T) for $\bar \Gamma_i^2$ using Ozawa's method \cite{Ozawa} before \cite{Oppenheim} was available.
    
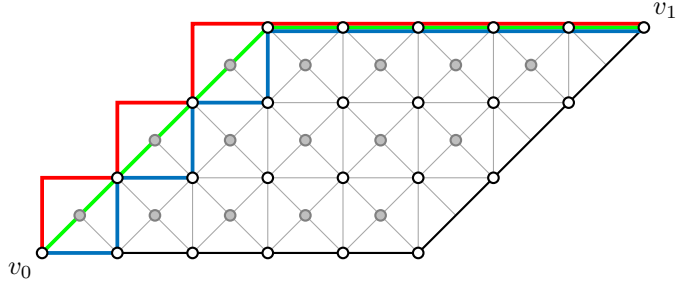
\begin{figure}
      \centering
\begin{tikzpicture}[cm={0,1,1,0,(0,0)}]

\coordinate (A) at (0,0);
\coordinate (B) at (3,3);
\coordinate (C) at (3,8);
\coordinate (D) at (0,5);

\begin{scope}
  \clip (A)--(B)--(C)--(D)--cycle;
  \foreach \x in {0,...,3}{
    \draw[gray!70] (\x,0)--(\x,8);
  }
  \foreach \y in {0,...,8}{
    \draw[gray!70] (0,\y)--(3,\y);
  }
  \foreach \x in {0,...,2}{
    \foreach \y in {0,...,7}{
      \draw[gray!70] (\x,\y)--++(1,1);
      \draw[gray!70] (\x+1,\y)--++(-1,1);
    }
  }
\end{scope}

\draw[thick] (A)--(B)--(C)--(D)--cycle;

\draw[line width=.5mm,green] (A)--($(B)+(0mm,0)$)--($(C)+(0mm,0)$);

\draw[line width=.5mm,RoyalBlue]
  (A)--(0,1)--(1,1)--(1,2)--(2,2)--(2,3)-- ($(B)+(-.5mm,0)$)--($(C)+(-.5mm,0)$);

\draw[line width=.5mm,red]
  (A)--(1,0)--(1,1)--(2,1)--(2,2)--($(3,2)+(.5mm,0)$)-- ($(B)+(+.5mm,0)$)--($(C)+(+.5mm,0)$);

  \foreach \x in {0,...,3}{
  \foreach \y in {0,...,5}{
  \draw[fill = white, thick] (\x,\y+\x) circle (0.7mm);
  }
  }
  \foreach \x in {0,...,2}{
  \foreach \y in {0,...,4}{
  \draw[color = gray, fill = lightgray, thick] (\x+.5,\x+\y+.5) circle (0.7mm);
  }
  }
  \node[below left] at (A) {$v_0$};
  \node[above right] at (C) {$v_1$};
\end{tikzpicture}
\caption{The combinatorial convex hull of two vertices $v$ and $w$ together with a canonical path of edges connecting them (green) and two semi-canonical paths of long edges connecting them (blue, red). The blue path is distinguished among the long-edge paths but we make no use of this.}
\label{fig:normal_forms}
\end{figure}

Throughout this paragraph let $X$ be a $\tilde{C}_2$-building. Let $\Gamma$ be a lattice acting freely on $X$ and regularly on vertices of each special type.
Let $Y = \Gamma \backslash X$ which is a $\tilde C_2$-GAB with two special vertices $v, w$.
Let $\bar \Gamma$ be an extension of $\Gamma$ that acts regular on special vertices of $X$. This determines an involutory automorphism $\rho$ on $Y$ swapping the two special vertices.
Lemma~\ref{lem:regular_pres} gives a presentation of $\bar \Gamma = \langle G \mid R \rangle$ with relations coming from squares as well as edges paired by the action of $\rho$. Note that we take into account all squares, not just a sufficient number to present the group. When we say that $g_{i-1}g_i = g_{i-1}'g_i'$ is a relation, we mean of course that $g_i^{-1}g_{i-1}^{-1}g_{i-1}'g_i'$ is a relator taking cyclic permutations and inverses into account. 
We want to describe normal forms for $\bar{\Gamma}$.

    The geometric starting point is:

    \begin{lemma}
      Let $X$ be a $\tilde{C}_2$-building and let $v_0$ and $v_1$ be special vertices of $X$. There is a unique path of edges $f_1,\ldots,f_{2k},g_1,\ldots,g_m$ (possibly $k = 0$ or $m = 0$) from $v_0$ to $v_1$ such that the $f_i$ are short edges, the $g_j$ are long edges, $f_{i-1}$ and $f_{i}$ as well as $g_{j-1}$ and $g_{j}$ meet in a vertex in which they form an angle of $\pi$ for $1 < i \le 2k$, $1 < j \le m$, $f_{2k}$ and $g_1$ meet in a vertex in which they form an angle of $3\pi/4$.
    \end{lemma}

    \begin{proof}
      The least convex subcomplex of $X$ containing $v_0$ and $v_1$ is a parallelogram with a $\tilde{C}_2$ tiling (see Figure~\ref{fig:normal_forms}): it is the intersection of all apartments containing $v_0$ and $v_1$ and therefore is the (combinatorial) convex hull of the two in any such apartment. For existence, take the path along the boundary of this convex hull. For uniqueness note that if two vertices are connected by such a path, the path runs along the boundary of the convex hull of the two vertices.
    \end{proof}

    \begin{corollary}\label{cor:semi-unique_paths_of_long_edges}
      Let $X$ be a $\tilde{C}_2$-building and let $v_0$ and $v_1$ be special vertices of $X$. There is a path of long edges $g_1,\ldots,g_\ell$ from $v_0$ to $v_1$ and a $k$ such that the edges $g_{j-1}$ and $g_{j}$ meet in a vertex in which they form an angle of $\pi/2$ for $1 < j \le k$ and an angle of $\pi$ for $k < j \le \ell$, and no three consecutive edges lie in the link of a common non-special vertex.
      Every other path of this form is obtained by a finite sequence of replacements of subpaths $g_{2i+1}g_{2i+2}$ by two long edges $g_{2i+1}'g_{2i+2}'$ connecting the same two special vertices. For the resulting path the value $k$ might have changed by~1.
    \end{corollary}

    \begin{remark}
      There is in fact a canonical path even among the ones formed of long edges, namely the one that lies in the convex hull of $v_0$ and $v_1$, which is distinguished as having one more turn than the other ones, see Figure~\ref{fig:normal_forms}. However this is less easy to identify by local conditions.
    \end{remark}

    Let $\omega = g_1 \cdots g_k$ be a sequence of oriented edges in $Y$ that start in $v$ and end in $w$.
    This sequence represents an element in $\bar{\Gamma}$. We say that $g_{i-1}$ and $g_{i}$ \emph{cancel} if $g_i = g_{i-1}^{-1}$ in $\bar{\Gamma}$, meaning that $\rho(g_i) = g_{i-1}^{-1}$, that they form a \emph{turn} if $g_{i-1} g_{i}$ is part of a relation and that they are \emph{straight} otherwise. We say that $g_{i-1}g_ig_{i+1}$ can be \emph{shortened} if $g_{i-1}g_ig_{i+1} = g'$ is a relation.

    From Corollary~\ref{cor:semi-unique_paths_of_long_edges} we get:

    \begin{corollary}\label{cor:semi-normal_form}
      Every element of $\bar{\Gamma}$ is represented by a word $\omega = g_1\ldots g_\ell$ for which there is a $k$ such that $g_{i-1}$ and $g_i$ form a turn for $i \le k$ and are straight for $i > k$ and no shortening is possible. Every other word of this form is obtained from a given one by applying a sequence of rewritings $g_{2i+1}g_{2i+2} \leadsto g_{2i+1}'g_{2i+2}'$ where $g_{2i+1}g_{2i+2} = g_{2i+1}'g_{2i+2}'$ is a relation.
    \end{corollary}

    \begin{definition}
      We say that a word $\omega \in F(G)$ is \emph{reduced} if its length is minimal among all words representing the same element of $\bar{\Gamma}$. We say that it is in \emph{semi-normal form} if it is of the form described in Corollary~\ref{cor:semi-normal_form}. We say that it is in \emph{normal form} if it is minimal in lexicographic order among all words in semi-normal form representing the same element of $\bar{\Gamma}$.
    \end{definition}

    \begin{proposition}
      \begin{enumerate}
        \item Any word in semi-normal form is reduced. Two words in normal form represent the same element of $\bar{\Gamma}$ if and only if they are equal.
        \item Any word of length $\ell$ can be brought into semi-normal form by applying $O(\ell^2)$ rewritings.
        \item Any word in semi-normal form of length $\ell$ can be brought into normal form by applying $O(\ell)$ rewritings.
      \end{enumerate}
    \end{proposition}

    \begin{proof}
    Let $v_0$ be a lift of $v$.
    Every word $\omega$ in $F(G)$ represents an edge path in $X$ along long edges starting in $v_0$ ending in $g.v_0$ where $g \in \bar{\Gamma}$ is the element represented by~$\omega$.
    Note that the graph distance of two special vertices in $X$ is the same whether measured with respect to all edges or with respect to only long edges.
    If $\gamma$ is an edge path from a special vertex $v_0$ to a special vertex $v_1$, let $\Sigma$ be an apartment containing~$v_0$ and $v_1$ and let $\rho_{\Sigma,c} \colon X \to \Sigma$ be the retraction centered at any chamber $c$. Then $\rho_{\Sigma,c}(\gamma)$ is an edge path from $v_0$ to $v_1$ in $\Sigma$. This shows that the minimal edge distance between $v_0$ and $v_1$ is realized inside $\Sigma$. It is easy to see that semi-normal paths are of minimal length among paths in $\Sigma$. Therefore, semi-normal forms are reduced. 
    That words in normal form represent elements uniquely is clear from the definition.

    Let $\omega = g_1, \dots,g_\ell$ be a word and without loss of generality assume that no cancellation or shortening in $\omega$ is possible.
    Let $i$ be the first index such that $g_{i-1}g_i$ is straight, let $a = \ell - i$.
    If no such $i$ exists, let $a \defeq 0$ for definedness.
    If there is a $k > i$ such that $g_{k-1}g_{k}$ is a turn, take the minimal one and let $b = k - i$.
    If no such~$k$ exists, let $b \defeq 0$ for definedness.
    Note that $b\leq a \leq \ell$ and that if $b=0$, then $\omega$ is in semi-normal form.
    We order the set of $(\ell, a,b)$ lexicographically and proceed by induction.
    There is a relation $g_{k-1}g_k \leadsto g_{k-1}'g_k'$ such that $g_{k-2}g_{k-1}'$ is a turn and we apply this rewriting at the positions $k-1$ and $k$ to obtain the word $\omega'$.
    If $k=i+1$ it might happen that $g_{i-2}g_{i-1}g_{k-1}'$ can be shortened to $g_{i-2}'$ and then we apply this shortening.
    The resulting word is of length $\ell-2$ and the index indicating the first straight pair of edges is reduced by two.
    Since the subword $g_1 \ldots g_{i-3} g_{i-2}'$ is reduced, further cancellations or shortenings might decrease the index that indicates the first straight pair, but at by most two per operation. In particular, we can apply cancellations or shortenings as long as possible and $a$ never increases. However, $b$~might increase. 
    If $k=i+1$ but we cannot shorten $g_{i-2}g_{i-1}g_{k-1}'$, then $\omega'$ is of length $\ell$ and the index indicating the first straight pair is at least $i+2$. Furthermore, the subword $g_1, \ldots g_{i+1}$ is reduced.
    Applying further cancellations and shortenings as long as possible, might decrease the index indicating the first straight pair below $i+2$ but only if they reduce the word length by at least the same amount.
    In particular, $\ell$ does not increase, $a$ decreases and $b$ might increase.
    If $k>i+1$, then $\omega'$ is of length $\ell$, the index indicating the first straight pair is $i$, but the distance to the next turn has decreased, i.e.~$b$ has decreased.
    Now performing cancellations and shortenings as long as possible will not increase $a$ by the same argumentation as before.
    In particular, $a$ and $\ell$ never increase. If they do not decrease, $b$ decreases. Therefore, repeating this procedure will eventually yield a word in semi-normal form.
    Since $b\leq a \le \ell$ the number of applications of relations is quadratic in $\ell$.
    
    If $g_1 \ldots g_\ell$ is in semi-normal form, then to obtain normal form we replace the pairs $g_{2i+1} g_{2i+2}$ successively from left to right by $g_{2i+1}' g_{2i+2}'$, choosing $g_{2i+1}'$ minimal among the admissible generators in the given order. This requires only a linear number of replacements, with a constant number of relations available at each step.
    \end{proof}

    We now complete the proof of Theorem~\ref{thm:main_2}.
\begin{proof}[Proof of Theorem \ref{thm:main_2}, Part 2]

    The subcomplex of $Y$ consisting of the triangles $t_1$ to $t_{18}$ (indicated in bold in Figure \ref{fig:TriangleComplexT}) is isomorphic to a subdivision $\dot{S}_R$ of the square complex $S_R$.
    The assignment indicated in Part 1 of this proof, induces a locally convex embedding of $\dot{S}_R$ into $Y_1^2$ by Lemma \ref{lem:locally_convex}. Therefore, the induced homomorphism
    \[
        \iota: \pi_1(S_R, x_{00}) = \pi_1(\dot{S}_R, x_{00}) \to \pi_1(Y_1^2, v)
    \]
    is injective.
    Let $\delta \in \pi_1(S_R, x_{00})$ be as in Lemma~\ref{lem:shorter_element} such that $\delta^4$ lies in the finite residual of $\pi_1(S_R, x_{00})$.
    Then $\bar \delta^4 \defeq \iota(\delta)^4$ lies in the finite residual of $\pi_1(Y_1^2, v)$.
    Tracing the maps we find that the element $\bar\delta$ is represented by the loop $f_1^{-1}*f_3$ in $\pi_1(Y^2_1,v)$. Using the triangles $t_2$ and $t_{12}$ one deduces that this loop is homotopic to $g_1 * g_6^{-1}$, thus in the presentation of Proposition~\ref{prop:presentation_bar_2} we have $\bar\delta = g_1 g_6^{-1}$.
    
    We now show that the normal closure of $\bar\delta^4$ in $\bar \Gamma_2^1$ is $\Gamma_2^1$, which is therefore the finite residual. We claim the following verifications are human-verifiable; however, they are long and tedious and the assistance of a computer is advantageous. To be more precise, we write $g_7g_4$ as a concrete product of elements which lie in the normal closure of $\bar\delta^4$. From there it is easy to see that the normal closure equals $\Gamma_2^1$. To perform computations in the presentation we make excessive use of our algorithm to compute normal forms. In summary, it requires 695 elementary operations, such as cancellations, shortenings or rewritings to deduce that $g_7g_4$ lies in the normal closure of $\bar\delta^4$.

    Let $\nu_0 \defeq \bar\delta^4$, let $\gamma_1 \defeq g_2 g_{11}$, and define $\nu_1 \defeq \nu_0 \gamma_1 \nu_0 \gamma_1^{-1}$.
    Computing the following normal form for $\nu_1$ requires 19 elementary operations.
    \[
        \nu_1 \defeq g_1 g_3 g_4 g_{13} g_4 g_2 g_6 g_{11} g_3 g_1 g_9 g_6 g_{11} g_2.
    \]
    Let $\gamma_2$ and $\gamma_3$ be defined as below. If we define $\nu_2 \defeq \gamma_2 \nu_1 \gamma_2^{-1}$, $\nu_3 \defeq \gamma_3 \nu_1^{-1} \gamma_3^{-1}$ and $\nu_4 \defeq \nu_2 \nu_3$, then one computes their normal forms using 30, 211 and 44 elementary operations, respectively.
    \begin{align*}
        \gamma_2 & \defeq g_3 g_6 g_{12} g_7 g_{11} g_2,\\
        \gamma_3 & \defeq g_{3} g_{2} g_{6} g_{6} g_{14} g_{5} g_{12} g_{6} g_{11} g_{9} g_{3} g_{14} g_{13} g_{2} g_{11} g_{15} g_{14} g_{8} g_{13} g_{5} g_{12} g_{5} g_{14} g_{4} g_{3} g_{1},\\
        \nu_2 &:= g_{7} g_{10} g_{15} g_{7} g_{6} g_{13} g_{3} g_{4} g_{6} g_{11} g_{3} g_{14} g_{13} g_{2},\\
        \nu_3 &:= g_{1} g_{4} g_{14} g_{2} g_{8} g_{8} g_{10} g_{15} g_{8} g_{7} g_{12} g_{5} g_{3} g_{15},\\
        \nu_4 &:= g_{7} g_{6} g_{13} g_{15} g_{7} g_{13} g_{9} g_{5} g_{4} g_{15}.
    \end{align*}
    Define $\gamma_4, \gamma_5$ as below.
    The product of the conjugates $\nu_5 \defeq \gamma_4 \nu_4 \gamma_4^{-1}$, $\nu_6 \defeq \gamma_5 \nu_4 \gamma_5^{-1}$ is is the short word $g_7 g_4 = \nu_5 \nu_6$. 
    One computes the normal forms of $\nu_5$ and $\nu_6$ indicated below using 192 and 188 elementary operations. Then computing the normal form of their product requires 11 elementary operations, which sums up to 695 elementary operations in total.
    \begin{align*}
        \gamma_4 \defeq & g_{1} g_{12} g_{5} g_{15} g_{1} g_{12} g_{5} g_{12} g_{11} g_{2} g_{10} g_{7} g_{2} g_{11} g_{9} g_{14} g_{8} g_{12} g_{8} g_{9} g_{11} g_{15} g_{14} g_{4} g_{2} g_{6} g_{7} g_{1} g_{12} g_{5}\\
        &g_{15} g_{11} g_{9},\\
        \gamma_5 \defeq & g_{1} g_{3} g_{5} g_{13} g_{4} g_{5} g_{11} g_{12} g_{10} g_{11} g_{2} g_{10} g_{7} g_{2} g_{11} g_{9} g_{14} g_{8} g_{12} g_{8} g_{9} g_{11} g_{15} g_{14} g_{4} g_{2} g_{6} g_{7} g_{1} g_{12}\\
        &g_{5} g_{15} g_{11} g_{9},\\
        \nu_5 \defeq & g_{7} g_{14} g_{7} g_{6} g_{9} g_{12} g_{13} g_{5} g_{3} g_{1}, \\
        \nu_6 \defeq & g_{1} g_{11} g_{12} g_{9} g_{1} g_{10} g_{11} g_{3} g_{2} g_{4}.
    \end{align*}
    Let $Q \defeq \bar \Gamma_2^1 / \langle \langle g_7g_4 \rangle \rangle$ and denote its generators by $g_i$. It easily follows that in $Q$ we have $x \defeq g_1 = g_3 = g_4 = g_6 = g_7 = g_9$ and $y \defeq g_2 = g_5 = g_8$. We observe that $g_{12}$ and $g_{13}$ are conjugated by $x$, so $g_{12}^2 = 1$ and therefore $z \defeq g_{12} = g_{14}$. 
    Using the relation $g_7 g_{10} g_{14} g_{14}$ yields $x = g_{10}$. Now the relation $g_1 g_{10} g_5 g_{11}$ implies $y = g_{11}$. We also have that $x = g_{15}$ since $g_7 g_{14} g_{14} g_{15}$ is a relation. Now using the relation $g_3 g_{12} g_{13} g_{15}$ yields $z = g_{13}$.
    Finally $g_3g_{11}g_{12}g_{13}$ yields $x=y$ and $g_5 g_{13}g_{11}g_{10}$ yields $x=z$. In particular, $Q$ is the cyclic group of order 2.
    \end{proof}

    \subsection{The full automorphism group of the building}
    
    Implementing the solution of the word problem of $\bar{\Gamma}^2_1$ in the previous paragraph on a computer, one can reconstruct finite balls in the Cayley graph of $\bar \Gamma_1^2$ and study their automorphism groups. It turns out that the automorphism group of a ball of radius $4$ pointwise fixes the ball of radius $2$, and in particular fixes all special vertices adjacent to it. 
    Since automorphisms of the Cayley graph induce automorphisms of $X_1^2$ and vice versa, it follows by induction that the stabilizer of a special vertex is trivial. From this we conclude:

    \begin{proposition}
      The lattice $\bar{\Gamma}^2_1$ is the full automorphism group and $\Gamma^2_1$ is the type-preserving automorphism group of the building $X^2_1$. In particular the building $X^2_1$ is exotic, since its automorphism group is discrete.
    \end{proposition}

    By combining the rigidity of balls in the Cayley graph of $\bar \Gamma_1^2$ with its non-residual finiteness, we obtain the following graph-theoretic consequence.
    If $Z$ is a locally finite, vertex-transitive graph, we call finite graph $Z_0$ a \textit{perfect finite $r$-local model for $Z$} if for any vertex $v_0 \in Z_0$ the $r$-ball around $v_0$ is isomorphic to the $r$-ball in $Z$ (around any vertex). Following the proof of \cite[Theorem 4]{Thom} we can show:

    \begin{proposition}\label{prop:no_local_model}
        For $r\geq 5$ there exists no perfect finite $r$-local model for the Cayley graph of $\bar \Gamma_1^2$ with respect to its generating set from Proposition \ref{prop:presentation_bar_2}.
    \end{proposition}
    
    \begin{proof}
        Denote the generating set of $\bar \Gamma_1^2$ by $S$ and the Cayley graph by $Z$.
        When considering balls of radius $r$ in $Z$, we will suppress the choice of the center and simply speak of ``the $r$-ball''.
        Assume $Z_0$ is a perfect finite 5-local model for $Z$.
        Let $v_0\in Z_0$ be a vertex.
        Any rooted isomorphism from the 4-ball in~$Z$ to the 4-ball in $Z_0$ around~$v_0$ induces an $S$-labeling on the 4-ball around $v_0$. Since any automorphism of a 4-ball in~$Z$ fixes the centered 2-ball, the induced labeling on the 2-ball around $v_0$ does not depend on the chosen isomorphism.
        In particular, we have a canonical labeling for every 2-ball in $Z_0$.
        
        To see that these labelings are compatible, consider adjacent vertices $v_0,w_0\in Z_0$. 
        Choose a rooted isomorphism from the 5-ball in $Z$ to the 5-ball in $Z_0$ around $v_0$ and use it to label the $5$-ball around $v_0$. 
        In particular, the 4-ball around $w_0$ gets labeled. Since this labeling arises from $Z$ we can find a rooted isomorphism from the 4-ball in $Z$ to the 4-ball around $w_0$ preserving the labels. In particular, the labeling on the 2-ball around $w_0$ induced by the labeling of the 5-ball around $v_0$ is the canonical one.
        
        Therefore, the local labelings are compatible and $Z_0$ admits a canonical $S$-labeling. This defines a vertex-transitive action of the free group over $S$ on $Z_0$. Since any relation of $\Gamma_1^2$ is witnessed in the ball of radius two, the action factors through $\Gamma_1^2$. But its finite residual is of index two, so $Z_0$ would have to be infinite or have at most two vertices, which is impossible.
    \end{proof}

    \subsection{The complexes $Y_k^3$}
    The complexes $Y_k^3$ are similar to the complex $Y_1^2$, but they are bigger. Each of them consists of 160 triangles, 120 edges and twelve vertices. We provide an explicit description of these complexes and a further details in appendix \ref{app:C2tilde_thickness4}. Given the complexes it is easy to verify that their universal covers $X_k^3$ are $\tilde C_2$-buildings of thickness $4$ and that they all contain a subdivision of $S_{\JW}$ as subcomplex. In particular their fundamental groups $\Gamma_k^3$ are non-residually finite $\tilde C_2$-lattices by Lemma \ref{lem:locally_convex} and Theorem \ref{thm:building_cartan-hadamard}. Each of the four complexes admits an involutory automorphism that swaps its two special vertices, which allows us to apply Lemma \ref{lem:regular_pres} and obtain a geometric presentation for the extension $\bar \Gamma_k^3$. Recall that Lemma \ref{prop:JanzenWise} provides explicit elements in $\pi_1(S_{\JW})$ and as in the proof of Theorem \ref{thm:main_2}, we compute the images of these explicit elements in $\bar \Gamma_k^3$. With the help of a computer, we verify that their normal closure is of finite index in~$\bar \Gamma_k^3$ and therefore this normal closure equals the finite residual of $\Gamma_k^3$ and $\bar \Gamma_k^3$. To be more precise, if~$\check \Gamma_k^3$ denotes the finite residual then $[\Gamma_k^3 : \check \Gamma_k^3] = 4$ if $k = 1,2$ and $[\Gamma_k^3 : \check \Gamma_k^3] = 8$ if $k = 3,4$ (we indicate the index in $\Gamma_k^3$ for consistency with the main theorem). As for $X_1^2$ we reconstruct finite balls in the Cayley complex for each case to compute the full automorphism group $\hat \Gamma_k^3$ of the buildings $X_k^3$. It turns out that we always have $[\hat \Gamma_k^3 : \Gamma_k^3] = 8$. In particular the automorphism group is discrete and $\check \Gamma_k^3 \backslash X_k^3$ is the maximal finite quotient of $X_k^3$ and therefore an invariant of the building. One can compute these quotients which are of course finite covers of the complexes $Y_k^3$ and as it turns out they are all different. In particular the buildings $X_k^3$ are pairwise not isomorphic.

    \subsection{Quasi-isometric rigidity}

    We recall quasi-isometric rigidity of Euclidean buildings in order to draw conclusions for lattices. The following is \cite[Theorem~1.1.3]{KleinerLeeb97} where it is asserted that the Moufang condition can be removed using \cite{Leeb00}. An alternative reference is \cite[Theorem~III]{KramerWeiss14} noting that the distance that a quasi-isometry of a building can have from an isometry is bounded in terms of the quasi-isometry constants \cite[Lemma~2.4]{KapovichKleinerLeeb98}.

    \begin{theorem}\label{thm:building_rigidity}
      Let $X$ and $Y$ be thick, irreducible Euclidean buildings of dimension $\ge 2$ and let $f \colon X \to Y$ be an $(A,B)$-quasi isometry. Then there is a unique isomorphism $\bar{f} \colon X \to Y$ at distance $\le D$ from $f$ where $D$ depends only on $A$ and $B$.
    \end{theorem}

    We immediately conclude.

    \begin{corollary}
      The lattices $\Gamma^2_1,\Gamma^3_1,\ldots,\Gamma^3_4$ are pairwise not quasi-isometric.
    \end{corollary}
    
    \begin{proof}
      By Svarc--Milnor $\Gamma^q_i$ is quasi-isometric to $X^q_i$. Theorem~\ref{thm:building_rigidity} translates quasi-isometry of the $\Gamma^q_i$ into isomorphism of the $X^q_i$. As we saw before the buildings $X^2_1,X^3_1,\ldots,X^3_4$ are pairwise non-isomorphic.
    \end{proof}

    Another consequence of Theorem~\ref{thm:building_rigidity} is quasi-isometric rigidity for uniform lattices.

    \begin{proposition}
      Let $X$ be a thick $2$-dimensional Euclidean building with $\Aut(X)$ discrete and let $\Gamma < \Aut(X)$ be a uniform lattice on $X$. If $\Lambda$ is quasi-isometric to $\Gamma$ then there is a finite index subgroup $\Lambda' < \Lambda$ and a homomorphism $\Lambda' \to \Gamma$ with finite kernel and finite-index image.
    \end{proposition}

    \begin{proof}
      Let $X$, $\Gamma$ and $\Lambda$ be as in the statement and pick a basepoint $o \in X$. The orbit map $\tau \colon \Gamma \to X, g \mapsto g.o$ is a quasi-isometry by the Svarc--Milnor lemma. By assumption there is also a quasi-isometry $\iota \colon \Lambda \to \Gamma$. So $\tau \circ \iota \colon \Lambda \to X$ is a quasi-isometry. By slightly perturbing $\tau \circ \iota$ we obtain an injective quasi-isometry $\kappa \colon \Lambda \to X$ and a quasi-inverse $\bar{\kappa} \colon X \to \Lambda$ satisfying $\bar{\kappa} \circ \kappa = \id$.

      For $g \in \Lambda$ let $\lambda_g \colon \Lambda \to \Lambda$ be left-multiplication, which is a quasi-isometry. Then $\tilde{\mu}_g = \kappa \circ \lambda_g \circ \bar{\kappa}$ is a quasi-isometry, and the constants are uniform in $g$. Moreover $g \mapsto \tilde{\mu}_g$ is a quasi action, in fact $\tilde{\mu}_1$ has bounded distance from the identity and $\tilde{\mu}_g \circ \tilde{\mu}_h = \tilde{\mu}_{gh}$.

    By Theorem~\ref{thm:building_rigidity} $\tilde{\mu}_g$ has bounded distance from an automorphism $\mu_g$, which is unique since the building is thick. It follows that the map $\mu \colon \Lambda \to \Aut(X), g \mapsto \mu_g$ is a homomorphism. Let $h = \bar{\kappa}(o)$. Then $\tilde{\mu}_g(o) = \kappa(gh)$, showing that $\Lambda \to X, g \mapsto \tilde{\mu}_g(o)$ is a quasi-isometry. Since the quasi-isometry constants of $\tilde{\mu}_g$ are uniform, the distance from $\tilde{\mu}_g$ to $\mu_g$ is uniformly bounded so $\Lambda \to X, g \mapsto \mu_g(o)$ is a quasi-isometry as well. Since $\Aut(X)$ is discrete, $\Gamma$ has finite index in it. So letting $\Lambda' = \mu^{-1}(\Gamma)$ the restriction of $\mu$ to $\Lambda'$ is the needed homomorphism.
    \end{proof}

    Applying this to our lattices without proper finite index subgroups we get:

    \begin{corollary}
        Let $\Gamma$ be one of $\check \Gamma_1^2, \check \Gamma_1^3, \dots, \check \Gamma_4^3$. If $\Lambda$ is quasi-isometric to $\Gamma$, then there is a finite index subgroup $\Lambda' \leq \Lambda$ and a finite normal subgroup $N\leq \Lambda'$ and such that $\Lambda'/N \cong \Gamma$.
    \end{corollary}

    \subsection{Property (T) using Ozawa's method}

    Lattices on (irreducible) Euclidean enjoy Kazhdan's property (T) \cite{Oppenheim}. This is particularly relevant as it represents half of the proof of the normal subgroup property, the other half being amenability of proper quotients. Previously, the best known bound \cite{Oppenheim15} for lattices on $\tilde{C}_2$-buildings of thickness $q+1$ to have (T) was $q \ge 4$, which does not apply to our lattices.

    We therefore established property (T) for $\bar{\Gamma}^2_1$ using Ozawa's method \cite{Ozawa}. Now the proof, besides showing that Ozawa's method can be applied, has the extra merit of providing quantitative information. Namely, $\bar{\Gamma}^2_1$ with respect to the generating set from Proposition \ref{prop:presentation_bar_2} has Kazhdan radius at most $2$ and Kazhdan constant least $0,4147$.

    We mimic the strategy in \cite{NetzerThom}. In what follows $\Gamma$ denotes a finitely generated group and we fix a symmetric generating set $S$ that does not contain $1$. Several of the upcoming definitions depend on $S$, which is not reflected in our notation. We denote by $\RR\Gamma$ the real group algebra of $\Gamma$, equipped with the anti-automorphism $*$ that takes $\gamma$ to $\gamma^{-1}$. Its fixed elements are called \emph{hermitian}. The (unnormalized) \emph{Laplacian} is
    \[
      \Delta = \abs{S} - \sum_{s \in S} s \in \RR\Gamma.
    \]
    The \emph{augmentation ideal} $I\Gamma$ is the kernel of $\chi \colon \RR\Gamma \to \RR, \gamma \mapsto 1$. The \emph{support} of $x \in \RR\Gamma$ is the set of $\gamma$ with non-zero coefficient. The \emph{square sum} of $q_1,\ldots,q_n \in \RR\Gamma$ is $x = \sum_i q_i^*q_i$ and it is a \emph{sum of squares decomposition} of $x$. It is \emph{supported} on a set $B \subseteq \Gamma$ if all $q_i$ have support in $B$. 
    
    Note that $\Delta \in I\Gamma$, and if $x$ is a sum of squares as above, then
    \[
        \chi(x) = \sum_i \chi(q_i^*)\chi(q_i) = \sum_i \chi(q_i)^2,
    \]
    so $x \in I\Gamma$ if and only if $q_i \in I\Gamma$ for all $i$.

    \begin{theorem}[Ozawa, \cite{Ozawa}]\label{thm:ozawa}
        The following are equivalent.
        \begin{enumerate}
            \item
            $\Gamma$ has property (T).
            \item 
            There exists an $\epsilon>0$ such that $\Delta^2 - \epsilon \Delta$ admits a sum of squares decomposition in $\RR\Gamma$.
        \end{enumerate}
    \end{theorem}

    The \emph{Kazhdan radius} is the infimal $R$ such that there is an $\epsilon > 0$ such that $\Delta^2 - \epsilon \Delta$ admits a sum of squares decomposition supported on the ball of radius $R$.

    We now discuss how the existence of a sum of squares decomposition can be approached via semidefinite programming. Let $b_1,\ldots,b_m \in \RR\Gamma$ be linearly independent and let $c_1,\ldots,c_n \in \Span_{\RR}\{b_i\}$ be written as $c_i = \sum_j q_{ij} b_j$ with $Q = (q_{ij})_{ij} \in \RR^{n \times m}$. If $x \in \RR\Gamma$ is a square sum
    \begin{equation}\label{eq:square_sum}
      x = \sum_i c_i^*c_i
    \end{equation}
    then
    \begin{equation}
      x = \sum_j b_j^* p_{jk} b_k\label{eq:semidefinite_decomposition}
    \end{equation}
    where $P = (p_{jk})_{jk} \in \RR^{m \times m}$ is the symmetric positive semidefinite matrix $P = Q^T Q$. Conversely, it $P \in \RR^{m \times m}$ is symmetric positive semidefinite satisfying \eqref{eq:semidefinite_decomposition}, writing $P = Q^TQ$ and $c_i = \sum_j q_{ij} b_j$ produces a sum of squares decomposition as in \eqref{eq:square_sum}.
    We call such a matrix $P$ a \textit{Gram matrix} for $x$ with respect to~$(b_1, \dots, b_m)$.
    
    Solving \eqref{eq:semidefinite_decomposition} is a semidefinite programming problem, for which solvers exist. The solution will usually involve some numerical error, so it is important to take this into account, which the following lemma does.

    \begin{lemma}[{\cite[Lemma 4.10]{KalubaKielakNowak}}]\label{lem:sufficient_T}
       Assume that $x \in \RR\Gamma$ belongs to the augmentation ideal, admits a sum of squares decomposition and is supported on the ball of radius $R$.
       Let $\epsilon>0$ and let $\nu \defeq \norm{\Delta^2 - \epsilon \Delta - x}_1$ be the 1-norm of this difference. Let $C = 2^{2\lceil \log_2R\rceil}$.
       Then $\Delta^2 - (\epsilon - \omega)\Delta$ admits a sum of squares decomposition supported in the ball of radius~$R$ for every $\omega \ge C\nu$.
    \end{lemma}
    In order to compute a sum of squares decomposition we proceeded as follows. Using the normal forms from Section~\ref{sec:normal_forms} we computed all 12526 elements in the ball of radius $4$ in $\Gamma_1^2$ and the multiplication table on the $166$ elements in ball of radius $2$. Denote these by $\gamma_1, \ldots, \gamma_{166}$. We used the Python-package Cvxpy (\cite{cvxpy}) and the solver Mosek (\cite{mosek}) to obtain an approximation to a Gram matrix $P$ for $\Delta^2-\epsilon\Delta$ with respect to $(\gamma_1, \dots, \gamma_{166})$ and $\epsilon = 1,29$.
    We computed a Cholesky decomposition $Q^T Q = P$. In order to do exact calculations, we approximated $Q$ by $10^{-12}Q'$ where~$Q'$ is an integer matrix. Since $\Delta^2-\epsilon\Delta$ lies in the augmentation ideal, we know that the elements represented by rows of $Q$ lie close to it, meaning that they sum to nearly zero.
    We adjusted $Q'$ slightly to obtain $Q''$ in which the sum of each row is zero.
    We computed the sum of squares $x$, indicated below, which also lies in the augmentation ideal:
    \begin{align*}
      q_i &= \sum_{j=1}^{166} Q''_{ij} \gamma_j & \text{and}& & x = \sum_j q_j^*q_j.
    \end{align*}
    We computed the 1-norm
    \[
      \norm{\Delta^2 - \epsilon \Delta - 10^{-24}x}_1 = 7589138977410503812 \cdot 10^{-24} \approx 7.59 \cdot 10^{-6}.
    \]
    Applying Lemma~\ref{lem:sufficient_T} we obtain:

    \begin{proposition}
    The lattice $\bar \Gamma_1^2$ has Property (T) with Kazhdan radius at most $2$ with respect to the generating set $S$ in Proposition \ref{prop:presentation_bar_2}. If $\Delta \in \RR \Gamma_1^2$ is the Laplacian then we have that
    \[
        \Delta^2 - 1,2899 \Delta
    \]
    admits a sum of squares decomposition. The Kazhdan constant of $\bar \Gamma_1^2$ with respect to $S$ is bounded below by $0,4147 < \sqrt{1,2899 \cdot 2/15} \leq \kappa$.
    \end{proposition}
    \begin{proof}
      The condition is that $\omega \geq 4 \cdot 7589138977410503812 \cdot 10^{-24} \approx 3,036 \cdot 10^{-5}$. Choosing $\omega = 0,0001$ gives the sum of squares decomposition. The Kazhdan constant follows using \cite[Remark 5.4.7]{Bekka}, taking into account that the Laplace operator there is normalized by $1/\abs{S}$.
    \end{proof}

%
%

\section{Searching for non-positively curved chamber complexes}\label{sec:searching}

In this section we describe the algorithm that was used to produce the examples in Section~\ref{sec:non-rf_c2tilde}. 

It can more generally be used to search for finite non-positively combinatorial chamber complexes in the sense of Section~\ref{sec:non-positively_curved_complexes}. The algorithm is a variation of the search algorithm that was used in \cite{Radu_NonDesarguesian} to find a triangle presentation for the Hughes plane.

Recall that a combinatorial chamber complex is a special case of a $\Delta$-complex. We equip vertices with cotypes and edges with types in a type set $I$. If the type set is $I = \{i,j,k\}$, then vertices of cotype $i$ are incident with edges of types $j,k$.

Here is an informal description of the algorithm: let $I= \{i,j,k\}$ be the type set. We decide beforehand on the set of vertices of each cotype $i \in I$ as well as their links. Note that the vertices in the link of a vertex of cotype $i \in I$ naturally carry types in $I \setminus \{i\}$. Now every vertex of cotype~$j$ in the link of a vertex of cotype $i$ corresponds to an edge of type $k$ and thus to a (unique) vertex of cotype $i$ in the link of a vertex of cotype $j$. Thus the first step toward constructing the complex is to pair the vertices of cotype $j$ in the links of vertices of cotype $i$ with the vertices of cotype $i$ in the links of vertices of cotype $j$.

Not every pairing gives rise to a $2$-complex, however, an edge in a link needs to come from a triangle in the $2$-complex, which also induces edges in the vertex links of the other two types. Thus edges $e,f,g$ in links in vertices $u,v,w$ can arise from the same triangle in the $2$-complex if the vertices of $e$ are paired with $v$ and $w$ and so on; and every edge needs to be part of such a triplet in order for the pairing to give rise to a $2$-complex. The search algorithm is essentially a greedy search through possible pairings, taking the triplets that can be formed as a score function. When there is a fixed subcomplex to be embedded into the complex (as there will be in our case) we can fix the according pairings and triangles and not change them during the search.

Rather than using disjoint links and pairing them, the actual algorithm will use the following representation of the situation which is more memory efficient. The local geometries play the roles of links, vertices represent vertices in links, and edges represent at the same time pairings and edges in links.

\begin{definition}
  Let $\Sigma$ be a connected, finite simplicial graph. Assume it is equipped with a type function that assigns to a vertex one of three types in $I$, and maps adjacent vertices to different types. For two types $i,j \in I$, we call the subgraph generated by the vertices of types $i$ and $j$ the \emph{local geometry} of type $\{i,j\}$; we denote it with $\Sigma_{ij}$. The \emph{angle} of the local geometry $\Sigma_{ij}$ is $\theta_{ij} = 2\pi/g_{ij}$, where $g_{ij}$ is the girth of $\Sigma_{ij}$. We call~$\Sigma$ a \emph{Radu graph}~if
  \begin{enumerate}[label=(\alph*)]
    \item 
      every local geometry has the same number of edges,
    \item
      every vertex in a local geometry has valency at least 2,
    \item
      $\theta_{ij} + \theta_{jk} + \theta_{ki} \le \pi$ for $I=\{i,j,k\}$.
  \end{enumerate}
\end{definition}

Given a Radu graph, we define the cotype of an edge $e$ as the unique type $k$, such that the vertices incident with $e$ have types $i,j$ and $I = \{i,j,k\}$.

\begin{lemma}
  Let $Y$ be a finite non-positively curved combinatorial chamber complex with type set $I$. Let $\Sigma$ be the graph defined as follows. The vertex set of $\Sigma$ is the set of edges of $Y$. We connect two vertices in~$\Sigma$ if the corresponding edges are incident with a common triangle. The type function on edges in $\Sigma$ is given by the  type function on edges in $Y$. Then $\Sigma$ is a Radu graph.
  Furthermore, a local geometry $\Sigma$ is the union of the vertex links in $Y$ of a given type.
\end{lemma}
\begin{proof}
  It is clear that the type function on $\Sigma$ is indeed a type function.
  As for a Radu graph, we call the subgraphs of $\Sigma$ generated by two types of vertices local geometries. These are finite bipartite graphs.
  Since every edge in $Y$ is contained in at least two triangles, the valency of a vertex in a local geometry of $\Sigma$ is at least 2. In particular, a local geometry contains a non-trivial cycle and has finite girth.
  Now let~$i$ be a type of $Y$ and let $j,k$ be the other two types. 
  By construction, the local geometry in~$\Sigma$ generated by vertices of types $\{j,k\}$ is isomorphic to the union of the vertex links of vertices of cotype $i$ in~$Y$.
  By non-positive curvature the links must not contain double edges. In particular, $\Sigma$ is simplicial. 
  The connectedness of $\Sigma$ follows from the fact that any two triangles in~$Y$ can be connected by a sequence of triangles that share an edge.
  Every triangle in $Y$ contributes exactly one edge to a vertex link of each cotype of $Y$. In particular, the number of edges in a local geometry is the number of triangles of $Y$.
  It is clear that the angle sum is at most $\pi$.
\end{proof}
\begin{definition}\label{def:perfect}
  Let $\Sigma$ be a Radu graph. A \emph{triangle} in $\Sigma$ is a circle of length 3 in $\Sigma$. Note that the vertices in a triangle are necessarily of three different types. A \emph{partial triangle cover} is a set of triangles such that every edge of $\Sigma$ is contained in at most one triangle of the family. It is an \emph{exact triangle cover} if every edge in $\Sigma$ is contained in exactly one triangle. A Radu graph is \emph{perfect} if it admits an exact triangle cover.
\end{definition}

\begin{proposition}
  A Radu graph is perfect if and only if it is the Radu graph of a finite non-positively curved combinatorial chamber complex.
\end{proposition}

\begin{proof}
  Let $Y$ be a finite non-positively curved combinatorial chamber complex and let $\Sigma$ be the Radu graph associated to it. We obtain a family of triangles in $\Sigma$ as follows. For every triangle in~$Y$ its boundary consists of three edges. These three edges in $Y$ correspond to three vertices in~$\Sigma$ that form a triangle. 
  Since edges in $\Sigma$ arise from boundaries of triangles, this triangle family is an exact cover.
  
  Now assume that we have a perfect Radu graph $\Sigma$.
  We will construct a combinatorial chamber complex complex $Y$ out of it as follows.
  Let $\Sigma_{ij}$ be the local geometry in $\Sigma$ generated by vertices of type $i$ and~$j$. 
  We call a connected component of $\Sigma_{ij}$ a link of cotype $k$, where $\{i,j,k\} = I$.
  For every link of cotype $k$, we have a vertex in~$Y$ of cotype $k$.
  For every vertex~$e$ in~$\Sigma$ of type~$i$ we now glue in an edge in~$Y$ between the pair of vertices in~$Y$, whose corresponding links in~$\Sigma$ contain $e$.
  Note that the two cotypes of these vertices in $Y$ are $j,k$, where $\{i,j,k\} = I$. In particular, the 1-skeleton of~$Y$ does not contain any loops.
  Let $T$ be an exact triangle cover for~$\Sigma$. For a triangle~$t$ in $T$, we glue in a triangle along the edges in $Y$ corresponding to three vertices in~$t$.
  We now verify the defining properties a of combinatorial chamber complex complex. By construction, the vertex links of $Y$ correspond to the links in $\Sigma$. Indeed, for a vertex~$u$ in $Y$, two edges incident with~$u$ share a triangle if and only if the corresponding vertices $e,f$ in $\Sigma$ are contained in a triangle. Since $T$ is an exact cover, this happens exactly when $e$ and $f$ are connected in $\Sigma$. We deduce that every edge is contained in at least two triangles, since every vertex in a vertex link has valency at least 2. Moreover, since $\Gamma$ is connected, any two triangles can be connected by a chain of triangles with consecutive triangles sharing an edge.
  The non-positive curvature follows from the third defining property of a Radu graph.
\end{proof}

\begin{figure}
  \begin{align*}
    \includestandalone[page=1]{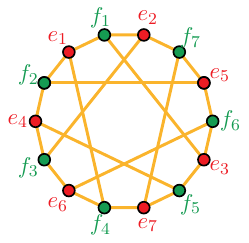}
            &&
            \includestandalone[page=2]{HeawoodRaduGraphPerfect}
            &&
            \includestandalone[page=3]{HeawoodRaduGraphPerfect}
  \end{align*}
  \caption{A perfect Radu graph admitting two exact covers.}
  \label{fig:perfectradugraph}
\end{figure}

\begin{remark}
  Many of the perfect Radu graphs we encountered during our experiments admit a unique triangle cover. In general a perfect Radu graph may admit more than one perfect triangle cover and if it does the associated triangle complexes may be non-isomorphic.
  Figure~\ref{fig:perfectradugraph} shows an example of a Radu graph that admits two triangle covers. For this example the two corresponding complexes are isomorphic.
  An example of a Radu graph arising from triangle complexes that are not isomorphic appears in \cite{Loue}.
\end{remark}
\subsection{Acting on Radu graphs}\label{subsec:acting_on_Radu_graphs}
In this subsection, we define an action on the space of Radu graphs with fixed vertex set and type function, which will later be used to construct perfect Radu graphs.

Let $V$ be a set of vertices with a type function $\type$ onto a three-element set and for $i \in I$ let $V_i$ be the set of vertices of type $i$. Let $G = \prod_i \Sym(V_i)$ be the group of type-preserving permutations of $V$. We will implicitly identify elements in $\Sym(V_i)$ with elements in $G$; in particular, we may omit writing elements as explicit triples.

Let $X$ be the set of Radu graphs with vertex set $V$ and type function $\type$.
We say that $\Sigma, \Sigma' \in X$ are \textit{isomorphic}, if there exists a type-preserving isomorphism between them;
we say that $\Sigma, \Sigma'$ are  \emph{locally isomorphic} if there exist type-preserving isomorphisms $\phi_{i,j} \colon \Sigma_{ij} \to \Sigma_{ij}'$ between the three local geometries.
We will often consider isomorphisms between local geometries as elements in $\Sym(V_i) \times \Sym(V_j)$.

We now define an action of $G$ on $X$.
We fix a cyclic ordering on $I$, represented by a cycle $\rho \in \Sym(I)$.
Let ${\Sigma \in X}$, let $i \in I$ and let $g \in \Sym(V_i)$. Let $j \defeq \rho(i)$ and $k \defeq \rho^2(i)$.  
We define $g.\Sigma$ as the Radu graph in $X$ with the same local geometries of type $\{j, k\}$ and $\{k, i\}$ as $\Sigma$, and with the following edges in the local geometry of type~$\{i,j\}$:  let $e \in V_i$ and $f \in V_j$ then $\{g(e),f\}$ is an edge in $g.\Sigma$ if and only if $\{e,f\}$ is an edge in $\Sigma$. This induces a type-preserving action on graphs with vertex set $V$, which we denote by $(g, \Sigma) \mapsto g.\Sigma$.

Note that isomorphisms between Radu graphs in $X$ are also naturally encoded by elements in $G$. We will be explicit throughout to avoid any ambiguity.

\begin{proposition}\label{prop:transitivity_Radu_graphs}
  The group $G$ acts on $X$. If two Radu graphs in $X$ lie in the same orbit, then they are locally isomorphic.
  If two Radu graphs $\Sigma, \Sigma'$ in $X$ are locally isomorphic, then an isomorphic copy of $\Sigma'$ lies in the orbit of $\Sigma$.
\end{proposition}

\begin{proof}
    It is clear that $G$ acts on $X$ and it is clear that two Radu graphs in the same orbit are locally isomorphic.
    On the other hand let $\Sigma, \Sigma'$ be locally isomorphic Radu graphs. 
    For $i \in I$ there is an isomorphism 
    \[
        \left(\phi^{i\rho(i)}_i, \phi^{i\rho(i)}_{\rho(i)}\right) \in \Sym(V_i) \times \Sym(V_{\rho(i)})
    \]
    between the two local geometries $\Sigma_{i\rho(i)}$ and $(\Sigma')_{i\rho(i)}$.
    Define $g$ as below. We claim that the permutation $h$, defined below, induces a graph isomorphism from~$g.\Sigma$ to~$\Sigma'$.
    \begin{align*}
        g &:= \left(\left(\phi^{\rho^2(i)i}_i\right)^{-1}\phi^{i\rho(i)}_i\right)_i,
        &
        h &:= \left(\phi^{\rho^2(i)i}_i\right)_i.
    \end{align*}
    Indeed, let $i$ be a type, let $j :=\rho(i)$ and $k :=\rho^2(i)$. Let $e \in V_i$ and $f\in V_j$. We have the following equivalences.
    \begin{align}
        & \left\{\left(\phi^{ki}_i\right)^{-1}\phi^{ij}_i(e), f\right\} \text{ is an edge in } (g.\Sigma)_{ij} \label{eq:action1} \\
        \iff & \{e,f\} \text{ is an edge in } \Sigma_{ij} \notag \\
        \iff & \left\{\phi^{ij}_i(e), \phi^{ij}_j(f)\right\} \text{ is an edge in } \Sigma'_{ij}.\label{eq:action2}
    \end{align}
    The permutation $h$ maps the edges in (\ref{eq:action1}) to the edges in (\ref{eq:action2}) and therefore induces a well-defined graph isomorphism from $g.\Sigma$ to $\Sigma'$.
\end{proof}

Thus looking for a triangle complex with prescribed local geometries $\Sigma_{ij}$ amounts to picking an arbitrary Radu graph $\Sigma$ with these local geometries and looking for a perfect one within the orbit~$G.\Sigma$.

  \subsection{The score of a Radu graph}

  Obviously $G$ is too large to do an exhaustive search. Rather we do a greedy search using the score function
        \[
          \score(\Sigma) \defeq \max
          \left\{  
            3 \cdot \lvert T\rvert \mid T \; \text{a partial triangle cover of } \Sigma
          \right\} .
        \]

        Note that a Radu graph is perfect if and only if its score equals its number of edges.

  An elementary but crucial observation is that small variations in $G$ (with respect to the generating set of transpositions) lead to small variations in score:

  \begin{lemma}\label{lem:quasi_cont}
    Let $\Sigma$ be a Radu graph on $X$. Let $\lambda_i$ be a transposition in $\Sym(V_i)$. Let $d_i$ be the maximal degree of a vertex of type $i$ in $\Sigma$. Then we have
    \[
      \abs{\score(\Sigma) - \score(\lambda_i .\Sigma)} \leq 3 d_i.
    \]
  \end{lemma}
  \begin{proof}
    Let $e, e'$ be the two vertices swapped by the transposition. Let~$T$ be a partial triangle cover of $\Sigma$. At most $(d_i/2)$ triangles contain $e$, and at most $(d_i/2)$ triangles contain $e'$. Removing these triangles from $T$ yields a partial triangle cover for $\lambda_i.\Sigma$. The lemma follows from the symmetry of the argument.
  \end{proof}
  Starting with a Radu graph $\Sigma$, we are looking for a finite non-positively cuvrved combinatorial chamber complex $Y$ whose Radu graph $\Sigma_Y$ is locally isomorphic to $\Sigma$. We employ the following strategy: We compute the score of $\Sigma$ and all Radu graphs obtained by acting with transpositions. Out of these we take the one(s) with the highest score and repeat the procedure. To avoid loops, we keep record of the Radu graphs we already checked. This search leads us to local maxima in the $G$-orbit of $\Sigma$. To avoid circulating around local maxima that are not perfect Radu graphs, we restart the search at a different Radu graph after a fixed number of steps.

  In fact, the score is expensive to compute and generally not worth knowing exactly. Rather we follow \cite{Radu_NonDesarguesian} in computing an estimate of the score in a greedy fashion, that has the important property that it equals the score for perfect Radu graphs. It is computed by first looking for an edge that is contained in a unique triangle (which therefore has to cover it in an exact cover) and adding that triangle to the cover. This is done as long as there are such edges. Finally, the exact maximal number of triangles covering the remaining edges is computed. There may be larger partial covers not using all of the initial triangles but if an exact cover exists, it has to use them.

  The above search strategy using the score estimate has been implemented in GAP~\cite{GAP4} and succeeded in producing the examples in Section~\ref{sec:non-rf_c2tilde}.
  The examples in Section~\ref{sec:a2tilde} were found using the estimated score function.

  \subsection{Stabilizers of isomorphism types}

  Certain elements $g \in G$ act on a Radu graph while preserving its isomorphism class. When classifying Radu graphs, it is therefore natural to quotient by the stabilizers of these isomorphism classes. This is the purpose of the following lemma.

    \begin{lemma}\label{lem:isomorphic_Radu}
    Let $\Sigma$ be a Radu graph in $X$ and for each type $i\in I$ let 
    \[
        \left(\phi^{i\rho(i)}_i,\phi^{i\rho(i)}_{\rho(i)}\right) \in \Sym(V_i)\times \Sym(V_{\rho(i)})
    \]be an automorphism of the local geometry of type $\{i, \rho(i)\}$. Let $g = (g_i)_i \in G$ and define $h \in G$ as below. Then $h.\Sigma$ and $g.\Sigma$ are isomorphic Radu graphs via the permutation $k$ defined below.
    \begin{align*}
        h &:= \left(\left(\phi^{\rho^2(i)i}_i \right)^{-1} \circ g_i  \circ \phi^{i\rho(i)}_i \right)_i,
        &
        k &:= \left(\phi^{\rho^2(i)i}_i \right)_i.
    \end{align*}
    \end{lemma}
    \begin{proof}
        Let $i$ be a type, let $j :=\rho(i)$ and $k :=\rho^2(i)$. Let $e \in V_i$ and $f \in V_j$. We have:
        \begin{align*}
            \{e,f\} \text{ is an edge in } \Sigma_{ij}
            \; \iff \;
            \left\{\phi^{ij}_i(e), \phi^{ij}_j(f)\right\} \text{ is an edge in } \Sigma_{ij}. 
        \end{align*}
        We use the above equivalence to characterize edges in $(g.\Sigma)_{ij}$.
        \begin{align}
            & \left\{\left(\phi^{ki}_i \right)^{-1} \circ g_i  \circ \phi^{ij}_i (e), f \right\} \text{ is an edge in } (h.\Sigma)_{ij}\label{eq:action3}\\
            \iff & \{e,f\} \text{ is an edge in } \Sigma_{ij} \notag\\
            \iff &  \left\{g_i \circ \phi^{ij}_i (e), \phi^{ij}_j(f)\right\} \text{ is an edge in } (g.\Sigma)_{ij}.\label{eq:action4}
        \end{align}
        The permutation $k$ maps the edges in (\ref{eq:action3}) to the edges in (\ref{eq:action4}) and therefore induces a well defined graph isomorphism from $h.\Sigma$ to $g.\Sigma'$.
    \end{proof}

\section{$\tilde{A}_2$-lattices}\label{sec:a2tilde}

To find the lattices in Section~\ref{sec:non-rf_c2tilde} we employed the strategy described in Section~\ref{sec:searching} on a search space that was too big to exhaust. In this section we provide the results of an exhaustive search of certain $\tilde{A}_2$-lattices. Before we embark on the actual discussion of these, we briefly discuss important classes of previously known $\tilde{A}_2$-lattices.

\subsection{Chamber-, vertex- and panel-regular lattices}

Lattices on $\tilde{A}_2$-buildings that act regularly on chambers were first studied by Ronan and Tits; Timmesfeld classifies chamber-transitive lattices on $\tilde A_2$-buildings with Desarguesian vertex links, see \cite{Ronan84, Tits85, Timmesfeld89}. 
There exist four chamber-regular $\tilde A_2$-lattices in thickness 3 and $44$ in thickness 9. The examples in thickness 3 are given by the presentations indicated below. 
In each presentation, the cyclic groups generated by $a$, $b$ and $c$ are edge stabilizers, and the groups generated by two of them are vertex stabilizers isomorphic to the Frobenius groups of order~$21$. The building can be recovered from the group as a coset geometry.
The lattice~$\Omega_3$ is an arithmetic subgroup of $\PGL_3(\FF_2((t)))$ \cite{KMW3} and $\Omega_1$ is an arithmetic subgroup of $\PGL_3(\QQ_2)$ \cite[Theorem~8]{KMW2}. The lattices $\Omega_2$ and~$\Omega_4$ act on distinct exotic buildings.
\begin{align*}
  \Omega_1 &= \langle a,b,c \mid a^3 = b^3 = c^3 = 1, (ab)^2 = ba, (bc)^2 = cb, (ca)^2 = ac \rangle, \\
  \Omega_2 &= \langle a,b,c \mid a^3 = b^3 = c^3 = 1, (ab)^2 = ba, (bc)^2 = cb, (ac)^2 = ca \rangle, \\
  \Omega_3 &= \langle a,b,c \mid a^3 = b^3 = c^3 = 1, (ab)^2 = ba, (c^2b)^2 = bc^2, (ac)^2 = ca \rangle, \\
  \Omega_4 &= \langle a,b,c \mid a^3 = b^3 = c^3 = 1, (ab)^2 = ba, (bc^2)^2 = c^2b, (ac)^2 = ca \rangle.
\end{align*}

Lattices that act type-rotatingly and regularly on all vertices were studied by Cartwright, Mantero, Steger, and Zappa \cite{Cartwright}. The examples in thickness 3 are the groups given by the presentations indicated below. 

            In each presentation, the generators take some base vertex to its neighbors as do their inverses, and thus correspond to edges. The relations are products of generators that take the base vertex back to itself and correspond to triangles. The building is the Cayley complex of the presentation. The groups $\Gamma_{A.i}$ are arithmetic lattices in $\PGL_3(\FF_2((t)))$, whereas the groups~$\Lambda_{B.i}$ and $\Lambda_{C.i}$ are arithmetic lattices in $\PGL_3(\QQ_2)$.
            \begin{align*}
                \Lambda_{A.1} &:=
                \langle
                    s_1, \dots, s_7 \mid 
                    s_1 s_2 s_4, \,
                    s_2 s_3 s_5, \,
                    s_3 s_4 s_6, \,
                    s_4 s_5 s_7, \,
                    s_5 s_6 s_1, \,
                    s_6 s_7 s_2, \,
                    s_7 s_1 s_3
                \rangle, \\
                \Lambda_{A.2} &:=
                \langle
                    s_1, \dots, s_7 \mid 
                    s_1 s_4 s_2, \,
                    s_3 s_3 s_1, \,
                    s_5 s_5 s_4, \,
                    s_6 s_6 s_2, \,
                    s_7 s_1 s_6, \,
                    s_7 s_2 s_5, \,
                    s_7 s_4 s_3
                \rangle, \\
                \Lambda_{A.3} &:=
                \langle
                    s_1, \dots, s_7 \mid
                    s_1 s_4 s_3, \,
                    s_2 s_1 s_6, \,
                    s_3 s_3 s_3, \,
                    s_4 s_2 s_5, \,
                    s_5 s_5 s_5, \,
                    s_6 s_6 s_6, \,
                    s_7 s_1 s_3, \,
                    s_7 s_2 s_6, \,
                    s_7 s_4 s_5
                \rangle, \\
                \Lambda_{A.4} &:=
                \langle
                    s_1, \dots, s_7 \mid
                    s_1 s_1 s_1, \,
                    s_1 s_3 s_2, \,
                    s_2 s_5 s_4, \,
                    s_3 s_3 s_5, \,
                    s_4 s_4 s_4, \,
                    s_6 s_6 s_2, \,
                    s_7 s_4 s_1, \,
                    s_7 s_5 s_6, \,
                    s_7 s_7 s_7
                \rangle, \\
                \Lambda_{B.1} &:=
                \langle
                    s_1, \dots, s_7 \mid
                    s_1 s_6 s_3, \,
                    s_2 s_5 s_6, \,
                    s_3 s_6 s_5, \,
                    s_4 s_3 s_5, \,
                    s_7 s_1 s_2, \,
                    s_7 s_2 s_4, \,
                    s_7 s_4 s_1
                \rangle, \\
                \Lambda_{B.2} &:=
                \langle
                    s_1, \dots, s_7 \mid
                    s_3 s_3 s_2, \,
                    s_3 s_5 s_6, \,
                    s_6 s_6 s_4, \,
                    s_5 s_5 s_1, \,
                    s_7 s_1 s_4, \,
                    s_7 s_2 s_1, \,
                    s_7 s_4 s_2
                \rangle, \\
                \Lambda_{B.3} &:=
                \langle
                    s_1, \dots, s_7 \mid
                    s_1 s_1 s_7, \,
                    s_1 s_3 s_6, \,
                    s_2 s_2 s_7, \,
                    s_2 s_6 s_5, \,
                    s_3 s_3 s_3, \,
                    s_4 s_4 s_7, \,
                    s_4 s_5 s_3, \,
                    s_5 s_5 s_5, \,
                    s_6 s_6 s_6
                \rangle, \\
                \Lambda_{C.1} &:=
                \langle
                    s_1, \dots, s_7 \mid
                    s_1 s_2 s_6, \,
                    s_1 s_3 s_5, \,
                    s_1 s_5 s_4, \,
                    s_2 s_4 s_5, \,
                    s_3 s_4 s_6, \,
                    s_7 s_7 s_6, \,
                    s_7 s_2 s_3
                \rangle.
            \end{align*}

            Lattices that preserve types and act regularly on panels of each type were described by Essert and further studied by the second author \cite{Essert,Witzel17}. The examples in thickness 3 are given by the following presentations.
            The group $\Sigma_1$ embeds as an index $3$ subgroup into $\Omega_3$ and is therefore arithmetic in $\PGL_3(\FF_2((t)))$. Similarly, the group $\Sigma_2$ embeds into $\Omega_4$ and is therefore exotic.
            \begin{align*}
                \Sigma_1 &:=
                \langle x,y,z \mid x^7 = y^7 = z^7 = xyz = x^3y^3z^3 =1 \rangle,
                \\
                \Sigma_2 &:=
                \langle x,y,z \mid x^7 = y^7 = z^7 = xyz^3 = x^3y^3z =1 \rangle.
            \end{align*}
          
\subsection{Type-preserving vertex-regular lattices}
    We classify lattices on $\tilde{A}_2$-buildings of thickness 3 that preserve types and act regularly on the vertices of each type. Every vertex regular lattice in \cite{Cartwright} contains such a lattice with index~$3$: the kernel of the action on types. Conversely, not every type-preserving lattice needs to have a type-transitive extension.
    
    \begin{theorem}\label{thm:heawood}
        \begin{enumerate}
            \item
            Consider the $\tilde A_2$-GABs $\Delta_1, \dots, \Delta_{13}$ indicated in Table \ref{tab:deltas} in Appendix \ref{app:a2tildeGABs}.
            Let $\Gamma_i$ be the fundamental group and let $X_i$ be the universal cover of $\Delta_i$ for $1\leq i \leq 13$.
            Then~$\Gamma_i$ is a uniform lattice on $X_i$ that acts regularly on vertices of each type. On the other hand, if~$\Gamma \curvearrowright X$ is such an action then it is equivariant to one of the actions $\Gamma_i\curvearrowright X_i$ for some $1\leq i \leq 13$.
            These 13 lattices are pairwise not isomorphic.
            \item
            The lattices $\Gamma_1, \Gamma_2, \Gamma_3$ embed as arithmetic lattices in $\PGL_3(\FF_2((t)))$ and the lattices $\Gamma_4, \Gamma_5, \Gamma_6$ embed as arithmetic lattices in $\PGL_3(\QQ_2)$. 
            In particular, they act on the associated Bruhat--Tits buildings. 
            The seven lattices $\Gamma_7, \dots, \Gamma_{13}$ act on different exotic $\tilde A_2$-buildings.
            In particular, the 13 lattices $\Gamma_i$ fall into nine quasi-isometry classes.
        \end{enumerate}
    \end{theorem}

    Before we prove the theorem, we provide some more information on the lattices, GABs and buildings in the Tables~\ref{tab:arit}, \ref{tab:exotic}, \ref{tab:twoballs}.
    
    \renewcommand{\arraystretch}{1.5}
    \setlength{\LTcapwidth}{\textwidth}
    \newcolumntype{P}[1]{>{\centering\arraybackslash}p{#1}}
    
    \begin{longtable}{|c|c|c|P{2cm}|P{3cm}|}
            \hline
            Lattice & Ambient $\textbf{G}(\KK)$ & $\Aut(\Delta)$ & $\Aut(\Delta) /\allowbreak \Aut^+(\Delta)$ & Extensions\\
            \hline
            $\Gamma_1$ & $\PGL_3(\FF_2((t)))$ & $(C_7 \rtimes C_3) \times C_3$ & $C_3$ & $\Sigma_1, \Omega_3,\allowbreak\Lambda_{A.1}, \Lambda_{A.2}, \Lambda_{A.3}$\\
            \hline
            $\Gamma_2$ & $\PGL_3(\FF_2((t)))$ & $\Sym(3)$ & $\Sym(3)$ & $\Lambda_{A.4}$ \\
            \hline
            $\Gamma_3$ & $\PGL_3(\FF_2((t)))$ & $(C_7 \rtimes C_6)$ & $C_2$ & $\Sigma_1, \Omega_3$ \\
            \hline
            $\Gamma_4$ & $\PGL_3(\QQ_2)$ & $(C_3 \times C_3) \rtimes C_2$ & $\Sym(3)$ & $\Lambda_{B.1},\Lambda_{B.2},\Lambda_{B.3}$\\
            \hline
            $\Gamma_5$ & $\PGL_3(\QQ_2)$ & $\Sym(3)$ & $\Sym(3)$ & $\Lambda_{C.1}$\\
            \hline
            $\Gamma_6$ & $\PGL_3(\QQ_2)$ & $\Sym(3)$ & $C_2$ & $-$ \\
            \hline
        \caption[Arithmetic lattices $\Gamma_1, \dots, \Gamma_6$.]{
            Arithmetic lattices $\Gamma_1, \dots, \Gamma_6$.
            With $\Delta$ we mean the associated GAB and $\Aut^+$ is the type-preserving automorphism group.
            We indicate extensions that arise as deck transformation groups covering subgroups of $\Aut(\Delta)$.
        }
        \label{tab:arit}
    \end{longtable}

    \begin{longtable}{|c|c|c|P{2cm}|c|}
            \hline
             Lattice & $\Aut(\Delta)$ &  $\abs{\Aut(X) : \hat \Gamma}$ & $\Aut(X) / \allowbreak \Aut^+(X)$ & Comments \\
             \hline
             $\Gamma_7$ & $C_2$ & 21 & $C_2$ &
             \begin{minipage}{5.2cm}
                \begin{tabular}{l}
                    $\Aut(X) = \Omega_2 \rtimes C_2$,\\
                    biggest solvable quotient: $C_6$ 
                \end{tabular}
             \end{minipage}\\
             \hline
             $\Gamma_8$ & $S_3$ & 1 & $C_2$ & 
             \begin{minipage}{5.2cm}
                 \begin{tabular}{l}
                    perfect, \\
                    smallest quotient: $\PSL_2(11)$
                 \end{tabular}
             \end{minipage}\\
             \hline
             $\Gamma_9$ & $S_3$ & 1 & $C_2$ &
             \begin{minipage}{5.2cm}
                 \begin{tabular}{l}
                      perfect, \\
                      smallest quotient: $\Alt(5)$ 
                 \end{tabular}
             \end{minipage}\\
             \hline
             $\Gamma_{10}$ & $C_3$ & 1 & $1$
             &
             \begin{minipage}{5.2cm}
                 \begin{tabular}{l}
                      perfect, \\
                      no quotient of size $<10^6$ 
                 \end{tabular}
             \end{minipage}
             \\
             \hline
             $\Gamma_{11}$ & $C_6$ & 1 & $C_2$ & 
             \begin{minipage}{5.2cm}
                \begin{tabular}{l}
                    biggest solvable quotient: $C_{13}$ 
                \end{tabular}
             \end{minipage}\\
             \hline
             $\Gamma_{12}$ & $C_3$ & 1 & $1$ & 
             \begin{minipage}{5.2cm}
                \begin{tabular}{l}
                    biggest solvable quotient: $C_2$ 
                \end{tabular}
             \end{minipage}\\
             \hline
             $\Gamma_{13}$ & $C_3$ & 1 & $1$ & 
             \begin{minipage}{5.2cm}
                 \begin{tabular}{l}
                      perfect, \\
                      no quotient of size $<10^6$ 
                 \end{tabular}
             \end{minipage}
             \\
             \hline
        \caption[Exotic lattices, automorphisms of their buildings and quotients.]{Exotic lattices, automorphisms of their buildings and quotients. With $X$ we mean the associated building and with $\hat \Gamma$ we mean the deck transformation extension $\Gamma \rtimes \Aut(\Delta)$ of $\Gamma$ induced by $\Aut(\Delta)$.}
        \label{tab:exotic}
    \end{longtable}
    
    \begin{longtable}{|c|c|c|c|c|c|c|c|}
            \hline
            Building & $X_7$ & $X_8$ & $X_9$ & $X_{10}$ & $X_{11}$ & $X_{12}$ & $X_{13}$ \\
            \hline
            2-ball around $\tilde u$ & 0 & 0 & 0 & 0 & 2 & 2 & 2\\
            \hline
            2-ball around $\tilde v$ & 2 & 2 & 2 & 2 & 0 & 0 & 0\\
            \hline
            2-ball around $\tilde w$ & 2 & 2 & 2 & 2 & 0 & 0 & 0\\
            \hline
        \caption[Balls of radius two in the exotic buildings $X_7, \dots, X_{13}$.]
        {Balls of radius two in $X_7, \dots, X_{13}$.
        There are two isomorphism types of such balls in $\tilde A_2$-buildings of thickness 3. Representatives are given by the $2$-balls in the Bruhat--Tits buildings associated with $\PGL_3(\QQ_2)$ and $\PGL_3(\FF_2((t)))$, respectively. In the table, the entries $0$ and $2$ indicate the isomorphism type: $0$ denotes the $2$-ball in the building associated with $\PGL_3(\QQ_2)$, while $2$ denotes the $2$-ball in the building associated with $\PGL_3(\FF_2((t)))$.}
        \label{tab:twoballs}
    \end{longtable}

    \begin{proof}[Proof for Table \ref{tab:arit}]
        In Table \ref{tab:arit}, we claim that the lattices $\Gamma_1, \dots, \Gamma_6$ are arithmetic lattices, we provide some information about cellular automorphisms of the quotient GAB, and we claim that our lattices admit extensions to the $\tilde A_2$-lattices mentioned at the beginning of this chapter. The claims about the cellular automorphism groups can be verified easily. In Table \ref{tab:autos_deltas} in Appendix~\ref{app:a2tildeGABs}, we provide explicit generators.
        The extensions to other $\tilde A_2$-lattices in Table \ref{tab:arit} all arise as deck transformations covering automorphisms of the quotient GAB. We provide the explicit automorphisms that give rise to the extensions in Table \ref{tab:extensions}. It now follows that the lattices $\Gamma_1, \dots, \Gamma_5$ are arithmetic, since they all embed as finite index subgroup in arithmetic groups.
        For arithmeticity of $\Gamma_6$ we refer to Theorem \ref{thm:gamma6}.
    \end{proof}

    \begin{proof}[Proof for Table \ref{tab:twoballs}]
        In Table \ref{tab:twoballs} we present the isomorphism types of the two balls in the buildings $X_7, \dots, X_{13}$.
        These balls have been constructed with the help of a computer and the computability of these balls follows from Section \ref{sec:a2tilde_normal_forms}, where we present an explicit algorithm for computing normal forms of combinatorial paths in certain $\tilde A_2$-buildings.
    \end{proof}

    \begin{proof}[Proof for Table \ref{tab:exotic}]
        In Table \ref{tab:exotic}, we provide information about the lattices $\Gamma_7, \dots, \Gamma_{13}$, their buildings and their quotient GABs. By Tables \ref{tab:twoballs}, these are exotic lattices. Again the explicit automorphism groups of the GABs can be computed easily and in Table \ref{tab:autos_deltas} we provide explicit generators.
        Let $7\leq i \leq 13$ and let $\Gamma \defeq \Gamma_i$. Let $X$ and~$\Delta$ be the associated building and GAB, respectively. 
        Let $\hat \Gamma$ be the extension of $\Gamma$ by $\Aut(\Delta)$.
        Note that $\Aut(\Delta)$ fixes the vertex $u$ so the stabilizer in $\hat \Gamma$ of any lift $\tilde u$ is naturally isomorphic to $\Aut(\Delta)$ and it is a complement to~$\Gamma$. Therefore, $\hat \Gamma = \Gamma \rtimes \Aut(\Delta)$.
        We need to show, that in every case, except for $i = 7$, the the full automorphism group of the building $X$ is $\hat \Gamma$. Fix a lift $\tilde u$ with stabilizer $D$ in~$\hat \Gamma$.
        Let~$B_4$ be the ball of radius four in $X$ around~$\tilde u$. In Table~\ref{tab:fourballs} we indicate the size of certain groups associated to $B_4$. In that table, the group~$A(u)$ denotes the automorphism group of the ball $B_4$.
        If we denote by~$B_1$ and~$B_3$ the balls of radius one and three around $\tilde u$, then $A^1(u)$ denotes the image of the restriction map $\Aut(B_4) \to \Aut(B_1)$ and $A_{(1)}^3(u)$ denotes the image of the fixator of $B_1$ in $\Aut(B_4)$ under the restriction map $\Aut(B_4) \to \Aut(B_3)$. The entries of the table have been calculated with a computer and it turns out that $A_{(1)}^3(u)$ is trivial in every case.
        The size of $A^1(u)$ equals the size of $D$ in every case, except for $i=7$.
        Now assume that $i \neq 7$.
        Let $\sigma$ be an automorphism of $X$ that fixes~$\tilde u$. We can find an element $\delta \in D$ such that $\sigma \circ \delta$ fixes $B_1$ pointwise. Since $A_{(1)}^3(u)$ is trivial, we know that $\sigma \circ \delta$ fixes $B_3$ pointwise. In particular, $\sigma\circ\delta$ fixes all vertices of the same type as $\tilde u$ that are at distance two and furthermore $\sigma\circ\delta$ fixes the 1-ball around these vertices. By induction we deduce that $\sigma \circ \delta$ is trivial and hence $\sigma \in D$. By Table \ref{tab:twoballs}, there cannot exist an automorphism of $X$, that maps $\tilde u$ to a vertex of another type. In particular, $\Aut(X) = \Aut(X)_{\tilde u}\Gamma = D\Gamma = \hat \Gamma$.
        For the case $i=7$ we refer to Proposition \ref{prop:gamma7embedding}. 
        The action on types of the full automorphism group of the building can be deduced by the previous discussion. The comments in the fifth column that concern quotients have been calculated with a computer. The comment regarding the automorphism group of $X_7$ follows again from Proposition \ref{prop:gamma7embedding}.
        \end{proof}

    \begin{proof}[Proof of Theorem \ref{thm:heawood}]
        The arithmeticity of $\Gamma_1,\dots,\Gamma_6$ follows from Table $\ref{tab:arit}$. To see that the groups $\Gamma_1, \dots, \Gamma_6$ are pairwise not isomorphic, we refer to Table \ref{tab:quotients_secondderived}, which shows that the lattices $\Gamma_1, \dots, \Gamma_6$ can be distinguished by their quotients modulo their second derived subgroups.
        We now show that the buildings $X_7, \dots, X_{13}$ are pairwise not isomorphic.
        By the Tables \ref{tab:exotic}, \ref{tab:twoballs}, we know that the only possible isomorphisms between these buildings could be between $X_8$ and $X_9$, or between~$X_{12}$ and $X_{13}$. Now assume that $X_8$ and~$X_9$ are isomorphic. Then the lattices $\Gamma_8$ and~$\Gamma_9$ are both index 3 subgroups in the type-preserving automorphism group of $X_8$. Note that $\Gamma_8$ and $\Gamma_9$ are not isomorphic since they have different smallest simple quotients. In particular, their intersection is a proper subgroup of index at most 3 in $\Gamma_8$ which implies that $\Gamma_8$ contains a normal subgroup of index at most $6 = 3!$. But the smallest quotient of $\Gamma_8$ is $A_5$, so its biggest normal subgroup is of index~60. Therefore, $X_8$ and $X_9$ cannot be isomorphic. The same argument shows that $X_{12}$ and $X_{13}$ are not isomorphic. 

        In summary, we have shown that the lattices $\Gamma_1, \dots, \Gamma_{13}$ are pairwise not isomorphic and that we have nine isomorphism classes of buildings, or equivalently nine quasi-isometry classes of groups.
        
        What remains to show is that any action of a lattice $\Gamma$ with the described properties is equivariant to one of the actions arising from the GABs in our list. 
        So let $\Gamma$ be such a lattice and let~$\Sigma_\Gamma$ be the Radu graph of the quotient GAB. 
        Recall that in Subsection \ref{subsec:acting_on_Radu_graphs} we defined a canonical action on Radu graphs with the same vertex sets. Now Proposition \ref{prop:transitivity_Radu_graphs} implies that an isomorphic copy of $\Sigma_\Gamma$ lies in the orbit of the Radu graph $\Sigma$ indicated in Figure \ref{fig:perfectradugraph}. In particular, we can obtain every lattice with the desired properties by computing the orbit of $\Sigma$, which is of size~$(7!)^3$. Lemma \ref{lem:isomorphic_Radu} indicates, that it is sufficient to consider a much smaller representative system of Radu graphs; a quantification of this is Lemma~\ref{lem:heawood_reps} below. By that lemma we only need to consider 604800 Radu graphs, check for perfectness, and obtain the equivariance classes of the lattices acting on their buildings. We performed this calculation on an ordinary office computer and it took a few minutes.
    \end{proof}

    \begin{lemma}\label{lem:heawood_reps}
        Let $\Sigma$ be the Radu graph indicated in Figure \ref{fig:perfectradugraph}.
        Let $I \defeq \{1,2,3\}$ be the set of types and denote the vertices of type $i$ by $V_i$. 
        Let $\rho \defeq (1 \; 2\;3) \in \Sym(I).$
        For $i\in I$ let $A^i\leq \Sym(V_i) \times \Sym(V_{\rho(i)})$ be the type-preserving automorphism group of the local geometry $\Sigma_{i\rho(i)}$. Denote the restriction of $A^i$ to $V_j$ with $A^i_j$ for $j \in \{ i, \rho(i)\}$. Now let $T_1$ be a double coset representative system for $A^3_1 \backslash \Sym(V_1) / A^1_1$ and let $T_3$ be a right transversal for $A^2_3 \backslash \Sym(V_3)$. Then for every Radu graph $\Sigma'$, that is locally isomorphic to $\Sigma$ there exist $\lambda_1 \in T_1$, $\lambda_2 \in \Sym(V_2)$ and $\lambda_3 \in T_3$ such that $\Sigma'$ is isomorphic to $\left(\lambda_1, \lambda_2, \lambda_3\right).\Sigma$.
    \end{lemma}

    \begin{proof}
        This is a direct consequence of Lemma \ref{lem:isomorphic_Radu}.
    \end{proof}

    \begin{proposition}\label{prop:gamma7embedding}
        The lattice $\Gamma_7$ embeds as an index 21 subgroup in the exotic lattice $\Omega_2$. Furthermore, the automorphism group of the corresponding building is $\Omega_2 \rtimes C_2$, which therefore contains $\Gamma_7$ as an index 42 subgroup.
    \end{proposition}
    
    \begin{proof}
        The lattice $\Omega_2$ admits a quotient $Q$ that is a non-split extension of $C_6$ by $\Alt(7)$. This quotient can be presented as follows.
        \[
        Q :=
        \left\langle
        a,b,c \;\middle|\;
        \begin{aligned}
        & a^3 = b^3 = c^3 = 1,\quad (ab)^2 = ba,\quad (bc)^2 = cb,\\
        & (ac)^2 = ca,\quad c a b^{-1} c a b c^{-1} b a^{-1} b^{-1} c^{-1} b a^{-1} b
        \end{aligned}
        \right\rangle .
        \]
        Let $K$ be the kernel of the corresponding epimorphism.
        If we denote the building that $\Omega_2$ acts on by $X$, then $Y \defeq K\backslash X$ is a finite GAB and $Q$ acts on it by deck transformations. Now the subgroup $R = \langle bca, abcb\rangle \subseteq Q$ is of index 21 and acts regular on vertices of each type.
        The quotient geometry $R \backslash Y$ is again a GAB and is in fact isomorphic to $\Delta_7$. It follows that~$\Gamma_7$ embeds as an index 21 subgroup in~$\Omega_2$, and it turns out, that $\Gamma_7 \cong \langle bca, abcb\rangle \leq \Omega_2$. The fact that $\Omega_2 \rtimes C_2$ is the full automorphism group of $X$ follows from Table \ref{tab:fourballs} with similar arguments as in the proof for Table~\ref{tab:exotic}.
    \end{proof}

    \subsection{The lattice $\Gamma_6$}

        We consider $\Gamma_6$ as the fundamental group $\pi_1(\Delta_6,w)$.
        Being the fundamental group of a 2-complex we obtain a presentation of $\Gamma_6$ by contracting a spanning tree.
        Choosing the tree consisting of the edges $e_4$ and $g_7$ yields the presentation
        \[
            \Gamma_6 = \left\langle
            \;
            \begin{aligned}
            &e_1, e_2, e_3, e_5, e_6, e_7,\\
            &f_1, \dots, f_7, g_1, \dots, g_6    
            \end{aligned}
            \;\middle|\;
            \begin{aligned} 
            &e_1f_1g_2, && e_1f_2g_1, && e_1f_3g_3, && e_2f_1g_1, && e_2f_4g_4, && e_2f_5g_5, && e_3f_1g_4, \\
            &e_3f_6g_2, && e_3f_7g_6, && f_2g_6,    && f_4g_3,    && f_6g_5,    && e_5f_2,    && e_5f_5g_4, \\
            &e_5f_7g_3, && e_6f_3g_5, && e_6f_4,    && e_6f_7g_2, && e_7f_3g_1, && e_7f_5g_6, && e_7f_6
            \end{aligned}
            \;
            \right\rangle .
        \]
        
        According to Table \ref{tab:autos_deltas} the automorphism group of $\Delta_6$ is isomorphic to $\Sym(3)$ and generated by
        \begin{align*}
            r &:=
            (e_1\, e_2\, e_3)(e_5\, e_6\, e_7)(f_2\, f_4\, f_6)(f_3\, f_5\, f_7)(g_1\, g_4\, g_2)(g_3\, g_5\, g_6),\\
            s &:= (e_1\, g_1)(e_2\, g_2)(e_3\, g_4)(e_4\, g_7)(e_5\, g_6)(e_6\, g_5)(e_7\, g_3)(f_4\, f_6)(f_5\, f_7).
        \end{align*}
        Denote the extensions of $\Gamma_6$ consisting of deck transformations covering $\langle r \rangle$ respectively $\Aut(\Delta_6)$ by $\bar \Gamma_6$ respectively by $\hat \Gamma_6$.
        The group $\Aut(\Delta_6)$ preserves the chosen spanning tree and therefore acts on the presentations defining presentations for $\bar \Gamma$ and $\hat \Gamma$. In what follows we consider $\Gamma_6$, $\bar \Gamma_6$ and $\hat \Gamma_6$ as finitely presented groups generated by $\{e_3, g_4\}, \{e_3, g_4, r\}$ and $\{e_3, g_4, r,s\}$. 
        Note that $s$ that swaps $e_3$ and $g_4^{-1}$. The purpose of this section is to describe the arithmetic structure of these groups.

        Let $\alpha$ be a square root of $-23$ and consider the field $K \defeq \QQ(\alpha)$ and its ring of integers $\calO_K \defeq \ZZ[\varphi]$ where $\varphi = (1+\alpha)/2$, equipped with the involution $\overline{\rule[1.5ex]{.1em}{0pt}\cdot\rule[1ex]{.1em}{0pt}}$ that fixes $1$ and interchanges $\alpha$ and $-\alpha$. Let $\omega \defeq (143 + 7\alpha)/2 = 68 + 7\varphi \in \calO_K$ and consider the matrix
      \[
        Q = 
        \left(\begin{array}{rrr}
        83 & \omega  & \omega \\
        \bar \omega  & 83 & \omega \\
        \bar \omega & \bar \omega & 83
        \end{array}\right) \in M_3({\cal O}_{\mathbb K}) \cap \GL_3(K)
      \]
      that is hermitian: $Q^T = \overline{Q}$.
      It defines a hermitian form $f(v,w) = v^T Q \overline{w}$ and general unitary group $\GU(f)$ consisting of similitudes, i.e.\ of matrices that preserve $f$ up to scalars. On an $\calO_K$-algebra $A$ with automorphism extending $\overline{\rule[1.5ex]{.1em}{0pt}\cdot\rule[1ex]{.1em}{0pt}}$ it is given by
      \[
        \GU(f)(A) = \{ (M,\lambda) \in GL_3(A) \times \GL_1(A) \mid M^TQ\overline{M} = \lambda Q\}.
      \]
      It admits the automorphism $\sigma \colon M \mapsto \overline{M}^T$.
      In terms of algebraic groups there is a $\ZZ$-group scheme $\bfG$ such that $\bfG(R) = \GU(f)(\calO_K \otimes_{\ZZ} R)$ for any ring $R$ and $\sigma$ becomes an algebraic automorphism of $\bfG$. Over $K$ we also consider the $K$-group $\PGU(f) = \GU(f)/Z(\GU(f))$.

      Note that $-23$ is a square in $\QQ_2$, so $\calO_K \otimes_{\ZZ} \QQ_2 = \QQ_2[X]/(X^2+23)$ is isomorphic to $\QQ_2 \times \QQ_2$ by the chinese remainder theorem and the involution is given by $\overline{(a,b)} = (b,a)$. Writing $M \in \GL_3(\calO_K \otimes_{\ZZ} \QQ_2) \cong \GL_3(\QQ_2) \times \GL_3(\QQ_2)$ as $(M_1,M_2)$ and $\lambda \in \GL_1(\calO_K \otimes_{\ZZ} \QQ_2) \cong \QQ_2^\times \times \QQ_2^\times$ as $(\lambda_1,\lambda_2)$ the defining relation of $\GU(f)$ can be written as
      \[
        Q^{-1} (M_1^T,M_2^T) Q = (\lambda_1,\lambda_2)(M_2^{-1},M_1^{-1})
      \]
      showing that $\GU(f)(\QQ_2) \cong \GL_3(\QQ_2) \times \GL_1(\QQ_2)$ and $\PGU(f)(\QQ_2) \cong \PGL_3(\QQ_2)$.

      Note that we now have two buildings of type $\tilde{A}_2$ in play: the building $X_6$ that is the universal cover of the GAB $\Delta_6$ used to define $\Gamma_6$, and the Bruhat--Tits building $X$ associated to $\PGU(f)(\QQ_2)$.

      We consider the group $\tilde{\Lambda} = \bfG(\ZZ[1/2]) = \GU(f)(\calO_K[1/2])$ and its image $\Lambda$ in $\PGU(f)(K)$. We also consider the extension $\hat{\Lambda} = \Lambda \rtimes \gen{\sigma}$. We write elements of $\Lambda$ as representatives in $\tilde{\Lambda}$ with square brackets instead of parentheses. Specifically, the following three are elements of $\Lambda$:
       \begin{align*}
        A & \defeq
        \begin{bmatrix}
            (-179 -31 \alpha)/4 & (-109 -17\alpha)/2  &   (-131 -9\alpha)/2 \\
            (-73 +47 \alpha)/4  & (-71+ 55\alpha)/4   &   (7+29\alpha)/2    \\
            (233 -3\alpha)/4    & (135-3\alpha)/2     &   63-6 \alpha
        \end{bmatrix}\\
        B & \defeq
        \begin{bmatrix}
            63 + 6\alpha & (135+3\alpha)/2 & (233+3\alpha)/4 \\
            (7-29\alpha)/2 & (-71-55\alpha)/4 & (-73-47\alpha)/4 \\
            (-131+9\alpha)/2 & (-109+17\alpha)/2 & (-179+31\alpha)/4
        \end{bmatrix}\\
        C & \defeq
        \begin{bmatrix}
            (-927+133\alpha)/2 & (-233+211\alpha)/2 & (-257+91\alpha)/2 \\
            455+67\alpha & (999-29\alpha)/2 & 239+11\alpha \\
            (281-227\alpha)/2 & (-479-187\alpha)/2 & -36-52\alpha
        \end{bmatrix}.
      \end{align*} 

        \begin{theorem}\label{thm:gamma6}
          There is a homomorphism $\Phi \colon \hat{\Gamma}_6 \to \hat{\Lambda}$ defined by
          \begin{align*}
            e_3 &\mapsto A &
            g_4^{-1} &\mapsto B&
        r &\mapsto C&
        s &\mapsto \sigma.
          \end{align*}
          It is injective and has cofinite image. In particular, $\hat{\Gamma}_6$ is a lattice in $\PGU(f)(\QQ_2) \cong \PGL_3(\QQ_2)$. There is an isomorphism $\iota \colon X_6 \to X$ of buildings that is $\hat{\Gamma}_6$-equivariant.
    \end{theorem}

    \begin{proof}
      It is easy to verify that $A$, $B$, $C$, $\sigma$ satifsy the defining relations of $\hat{\Gamma}_6$ and that the image of an infinite order element of $\hat{\Gamma}_6$ has infinite order in $\hat{\Lambda}$. Consequently $\Phi$ is well-defined and has infinite order image.

      It follows from \cite[Theorem~1.1]{BaderCapraceLecureux} that $X_6$ is a Bruhat--Tits building. More precisely, geometry tells us that it is associated to $\PGL_3(D)$ where $D$ is a finite-dimensional division algebra over a local field with residue field $\FF_2$. Moreover by \cite[Theorem~1.2]{BaderCapraceLecureux} a finite index subgroup $\Gamma'$ of $\Gamma_6$ is contained in $\PGL_3(D)$.

      Margulis Superrigidity \cite[Theorem~(2)]{Margulis} implies that $\Phi$ extends uniquely to a homomorphism $\tau \colon \PGL_3(D) \to \PGU(f)(\QQ_2)$. Since $\PGL_3(D)$ is simple, $\tau$ is injective. From the Borel--Tits classification of abstract homomorphisms \cite{BorelTits73} it follows that $D = \QQ_2$ and $\tau$ is an isomorphism. The isomorphism of groups $\tau$ induces an isomorphism $\iota$ of buildings $X_6 \to X$, which is $\Gamma'$-equivariant by construction. Since $\Gamma'$ is a uniform lattice, the action of $\hat{\Gamma}_6$ on $X_6$ and $X$ is uniquely determined by the action on $\Gamma'$, hence $\iota$ is $\hat{\Gamma}_6$-equivariant.
    \end{proof}   

    \begin{remark}
      We expect that $\Phi$ is in fact an isomorphism. Verifying this expectation amounts to showing that $\Lambda \cap K$ has order $3$ where $K$ is the image of $\GL_3(\ZZ_2)$ in $\PGL_3(\QQ_2)$. Indeed, $\Lambda$ contains $\Phi(\bar{\Gamma}_6)$ which is transitive on each type of vertex, so the index can be computed as the index of stabilizers of any special vertex. In $\bar{\Gamma}_6$ such a stabilizer is $\gen{C}$ and in $\PGL_3(\QQ_2)$ it is $K$. So the order computation would show that $\Phi|_{\bar{\Gamma}_6}^\Lambda$ is an isomorphism and it extends to $\Phi$ because $[\hat{\Lambda}:\Lambda] = 2 = [\hat{\Gamma}_6:\bar{\Gamma}_6]$.
    \end{remark}

    \begin{lemma}\label{lem:type_preserving}
        The action of the lattice $\Lambda$ on $X$ preserves types.
    \end{lemma}

    \begin{proof}
      Pick an embedding of $K$ into $\QQ_2$ and an isomorphism $\PGU(f)(\QQ_2) \cong \PGL_3(\QQ_2)$. An element $[g] \in \PGL_3(\QQ_2)$ acts on types by $\nu(\det g) \bmod 3$ where $\nu$ is the valuation of $\QQ_2$ (which is well-defined because the center changes the valuation by a third power). The defining relation $M^TQ\overline{M} = \lambda Q$ asserts that an element $[M]$ of $\Lambda$ has determinant $d \in \calO_{K}[1/2]^\times$ with $d\bar{d} = \lambda^3$, and it suffices to show that $\nu(d) \equiv 0 \bmod 3$.

      The ideal $(2) \lhd \calO_{K}$ is contained in two maximal ideals $\mathfrak p = (2,\varphi)$ and $\overline{\mathfrak p} = (2,\bar{\varphi})$, which are the intersections $2\ZZ_2$ with $\calO_{K}$ under the two possible embeddings $\calO_K \to \QQ_2$ and therefore define the restriction of $\nu$ and its conjugate valuation $\bar{\nu}$. From the last paragraph we have $(d) = \mathfrak{p}^{\nu(d)} \overline{\mathfrak{p}}^{\bar{\nu}(d)}$ (since $d$ is invertible after localizing at $2$) and $\nu(d) + \bar{\nu}(d) = \nu(d\bar{d}) = 3\nu(\lambda) \equiv 0 \bmod 3$.

      Since $\mathfrak{p} \mathfrak{\bar{p}} = (2)$ we have $[\mathfrak{p}] = [\overline{\mathfrak{p}}]^{-1}$ in the ideal class group $\Cl(K)$, so $[(d)] = [\mathfrak{p}]^{\nu(d) - \bar{\nu}(d)}$. Since $[\mathfrak{p}]$ has order $3$ (it is not principal but $\mathfrak{p}^3 ((3+\alpha)/2)$ is and, in fact, the class group $\Cl(K)$ has order $3$ \cite{number_field_database}) we deduce $\nu(d) - \bar{\nu}(d) \equiv 0 \bmod 3$. Together with the last paragraph we conclude $\nu(d) \equiv 0 \bmod 3$.
    \end{proof}

    \begin{remark}
        We computed a representation of $\Gamma_6$ with a similar method as described in Remark \ref{rem:rep_gamma_r}.
        We assume that $\Gamma_6$ admits a representation $\Phi$ to a matrix group in $\SL_3(K)$, where $K$ is a number field.
        We fix a $\Phi_3$ representation to $\SL(3,3)$ and then Hensel-lifted the homomorphism to sufficiently high precision $\Phi_{3^k} : \Gamma_6 \to \SL_3(\ZZ/3^k\ZZ)$.
        As in Remark \ref{rem:rep_gamma_r} we use the approximations to get algebraic numbers as candidates for the characteristic polynomial of $A \defeq \Phi(e_3)$ and the traces of $B \defeq \Phi(g_4^{-1}), \allowbreak AB,\allowbreak A^2B, \allowbreak B^2, \allowbreak AB^2, \allowbreak A^2B^2$.
        We assume that $A$ is diagonal and our candidate for the characteristic polynomial determines the diagonal entries. Since $A$ is diagonal, we can use the candidate values for the traces of $B, AB, A^2B$ to obtain a linear system of equations in the diagonal entries $B_{11}, B_{22}, B_{33}$ of $B$, which we solve.
        
        The diagonal entries of \(B^2\) are given by
        $
        (B^2)_{11}
        =
        B_{11}^2+B_{12}B_{21}+B_{13}B_{31}, \;
        (B^2)_{22}
        =
        B_{22}^2+B_{12}B_{21}+B_{23}B_{32}, \;
        (B^2)_{33}
        =
        B_{33}^2+B_{13}B_{31}+B_{23}B_{32}.
        $
        Define $x_1=B_{12}B_{21}, \; x_2=B_{13}B_{31}, \; x_3=B_{23}B_{32}.$
        Since $A$ is diagonal, the candidate values for the traces of $B^2, AB^2, A^2B^2$ yield a system of three linear equations in $x_1,x_2,x_3$, which we solve.
        
        Define $u \defeq B_{12} B_{23} B_{31}$,  $v \defeq B_{13} B_{21} B_{32}$ and $S \defeq u+v$. Since $\det(B)=1$, we have that $S = 1 - B_{11}B_{22}B_{33} + B_{11}x_3 + B_{22} x_2 + B_{33} x_1$. Since $x_1x_2x_3 = u v$, we conclude that $u,v$ are the roots of $t^2 - St + x_1x_2x_3$. We solve this equation and obtain candidates for $u$ and $v$.
        The candidate values for $u$ and $v$ are not $0$ and in particular, $B_{12}, B_{13}$ are both not $0$.
        After conjugating by a suitable diagonal matrix, we may assume that $B_{12} = B_{13} = 1$, which implies that $B_{21} = x_1$, $B_{31} = x_2$, $B_{23} = u/x_2$, $B_{32} = x_2x_3/u$.
        Using this strategy we successfully constructed a representation of $\Gamma_6$ to $\SL_3(\LL)$, where $\LL=\QQ(\sqrt{-15}, \sqrt{-23})$.

        From there, we used ad hoc methods to modify the representation so that the image lies in $\SL_3(\ZZ[1/2,\sqrt{-23}])$ and preserves the unitary form described above.
        We then computed a matrix inducing the cyclic extension. The extension does not admit a realization in $\SU_3$, but instead can be realized in the projective group $\PGU_3$.
    \end{remark}

\subsection{Normal forms}\label{sec:a2tilde_normal_forms}

\begin{figure}
\centering
\begin{tikzpicture}
\def\gridgray{black!60}
\def\nodefill{black!35}

\tikzset{
  gridline/.style={draw=\gridgray, line width=.25pt},
  boundary/.style={draw=black, line width=.8pt},
  pathBlue/.style={draw=RoyalBlue, line width=.5mm},
  vImportant/.style={circle, draw=black, fill=white, line width=.6pt, inner sep=0pt, minimum size=1.4mm},
  vOther/.style={circle, draw=black, fill=\nodefill, line width=.6pt, inner sep=0pt, minimum size=1.4mm},
  arr/.style={gridline, ->-},
  arrBack/.style={gridline, -<-},
}

\def\a{8}
\def\b{5}

\coordinate (A) at (0,0);
\coordinate (B) at (\a,0);
\coordinate (D) at (60:\b);
\coordinate (C) at ($(B)+(D)$);

\begin{scope}
  \clip (A)--(B)--(C)--(D)--cycle;

  \foreach \x in {0,...,\a}{
    \foreach \y in {0,...,\b}{
      \draw[arr]     ($(0:\x)+(60:\y)$) -- +(60:1);
      \draw[arrBack] ($(0:\x)+(60:\y)$) -- +(0:1);
      \draw[arrBack] ($(0:\x)+(60:\y)$) -- +(120:1);
    }
  }
\end{scope}

\draw[boundary] (A)--(B)--(C)--(D)--cycle;

\draw[pathBlue] (A)--(D)--(C);

\foreach \x in {0,...,\numexpr\b-1\relax}{
  \draw[pathBlue,->-] (60:\x) -- +(60:1);
}
\foreach \x in {0,...,\numexpr\a-1\relax}{
  \draw[pathBlue,->-] ($(60:\b)+(0:\x)$) -- +(0:1);
}

\foreach \x in {0,...,\a}{
  \foreach \y in {0,...,\b}{
    \pgfmathparse{mod(\x - \y,3) == 0 ? int(1) : int(0)}
    \ifnum\pgfmathresult=1
      \node[vImportant] at ($(0:\x)+(60:\y)$) {};
    \else
      \node[vOther] at ($(0:\x)+(60:\y)$) {};
    \fi
  }
}

\node[below left] at (A) {$v_0$};
\node[above right] at (C) {$v_1$};
\end{tikzpicture}

\caption[An edge path in normal form in an $\tilde A_2$-building.]{An edge path in normal form consisting of straight forward edges, a wide turn and straight reverse edges. The gray arrows indicate the preferred orientation, and the blue arrows indicate the orientation in which the edges are traversed.}
\label{fig:a2tilde_normal_forms}
\end{figure}

Throughout this section let $X$ be a $\tilde{A}_2$-building, let $\Gamma$ be a uniform lattice acting freely and type-preservingly on $X$, and let $Y = \Gamma \backslash X$.

This subsection is similar to Section~\ref{sec:normal_forms}, where we computed normal forms for words in $\tilde C_2$-lattices that act regularly on special vertices. Here, however, we take a slightly different approach.
We embrace the fact that $\Gamma \backslash X$ has more than one vertex and give a normal form for edge paths up to homotopy. Conceptually, we work with the fundamental groupoid of $\Gamma \backslash X$ and regard $\Gamma$ as the subgroup of homotopy classes of paths that start and end at the same vertex.

We put a cyclic ordering $\rho \in \Sym(I)$ on the set of types $I$, thus giving each edge of type $\{i,j\}$ a preferred orientation $(i,j)$ such that $\rho(i) = j$. We say that an oriented edge is a \emph{forward} edge if its orientation agrees with the preferred orientation, and is a \emph{reverse} edge if it does not.

\begin{lemma}
  Let $v_0$ and $v_1$ be vertices in $X$ of the same type. There is a unique path of edges $f_1,\ldots,f_{k},g_1,\ldots,g_m$ (possibly $k = 0$ or $m = 0$) from $v_0$ to $v_1$ such that the $f_i$ are forward edges, the $g_j$ are reverse edges, $f_{i-1}$ and $f_{i}$ as well as $g_{j-1}$ and $g_{j}$ meet in a vertex in which they form an angle of $\pi$ for $1 < i \le k$, $1 < j \le m$, $f_{k}$ and $g_1$ meet in a vertex in which they form an angle of $2\pi/3$.
\end{lemma}

\begin{proof}
  The least convex subcomplex of $X$ containing $v_0$ and $v_1$ is a parallelogram with a $\tilde{A}_2$-tiling (see Figure~\ref{fig:a2tilde_normal_forms}): it is the intersection of all apartments containing~$v_0$ and $v_1$ and therefore is the combinatorial convex hull of the two in any such apartment. For existence, take the path along the boundary of this convex hull. For uniqueness note that if two vertices are connected by such a path, the path runs along the boundary of the convex hull of the two vertices.
\end{proof}

Let $\omega = g_1 \cdots g_k$ be a path of oriented edges in $Y$.
We say that $g_{i-1}$ and $g_{i}$ \emph{cancel} if $g_i = g_{i-1}^{-1}$; that they form a \emph{sharp turn} if they form an angle of $\pi/3$; that they form a \emph{wide turn} if they form an angle of $2\pi/3$; and are \emph{straight} if they include an angle of~$\pi$. Note that straight edges and edges forming a sharp turn are either both forward or both reverse while wide turns represent a change from forward to reverse or vice versa. We say that $\omega$ is \emph{reduced} if its length is minimal among all paths with the same endpoints. We say that it is in \emph{normal form} if it contains no cancellations or sharp turns, at most one wide turn, which then has to be from a forward edge to a reverse edge.

\begin{corollary}\label{cor:a2tilde_unique_normal_form}
  Every edge path in $Y$ is homotopic to a unique one in normal form.
\end{corollary}

When manipulating edge paths we use the following rewritings:
free cancellations; rewriting $g_{i-1}g_i \leadsto h$ when $g_{i-1}g_i$ is a sharp turn in the triangle $g_{i-1}g_ih$;
and rewriting $g_{i-1}g_i \leadsto g_{i-1}'g_i'$ when $g_{i-1}g_i$ form a wide turn and $h$ is such that
$(g_{i-1}')^{-1} g_{i-1} h$
and
$g_i' g_i^{-1} h$
form triangles, and $g_{i-1}g_i$ and $g_{i-1}'g_i'$ form a parallelogram.

\begin{proposition}
  \begin{enumerate}
    \item Any edge path in normal form is reduced. Two paths in normal form are homotopic if and only if they are equal.
    \item Any edge path of length $\ell$ can be brought into normal form by applying $O(\ell^2)$ rewritings.
  \end{enumerate}
\end{proposition}

\begin{proof}
    Let $\omega$ be an edge path in $X$ from a vertex $v$ to a vertex $w$. Let $\Sigma$ be an apartment containing both and let $\rho_{\Sigma,c}$ be the retraction centered at some chamber $c \subseteq \Sigma$. Then $\rho_{\Sigma,c}(\omega)$ is an edge path from $v$ to $w$, showing that the minimal edge distance is realized in $\Sigma$. In $\Sigma$ it is easy to see that paths in normal form are reduced. That homotopic paths in normal form are equal follows from Corollary~\ref{cor:a2tilde_unique_normal_form}.
    
    Let $\omega = g_1\ldots g_\ell$ be an edge path.
    Without loss of generality assume that $\omega$ only involves straight pairs of edges and wide turns.
    Let~$i$ be the first index such that~$g_i$ is a reverse edge and define $a \defeq \ell -i$. 
    For definedness set $a \defeq 0$ if the path contains no reverse edge.
    If the word is not in normal form there exists a minimal $k$ greater than~$i$ such that~$g_k$ is a forward edge.
    The path $g_i \ldots g_k$ lies in an apartment and we can apply~${k-i}$ rewritings to obtain a path $g_i'\ldots g_k'$ such that $g_i'$ is forward and $g_{i+1}'\ldots g_k'$ only consists of reversed edges.
    The path $\omega' = g_1 \dots g_{i-1} g_i' \ldots g_k' g_{k+1} \ldots g_\ell$ has its first reverse edge at position $i+1$.
    If $g_{i-1} g_i'$ is a sharp turn, we rewrite it and obtain a path with its first reverse edge at position $i-1$ and length $\ell - 1$. In particular~$\ell$ decreased and $a$ did not increase. If we now apply a cancellation or a rewrite a sharp turn, the position of the first reverse edge may decrease, but at most by the amount by which the word length decreases.
    If $g_{i-1} g_i'$ is straight, then~$\ell$ remains the same and $a$ decreased.
    Again if we apply a cancellation or rewrite a sharp turn, the position of the first reverse edge may decrease but at most by the amount by which the word length decreases.
    In particular, performing these operations as long as possible results in both cases in a path $\omega''$ for which, compared to $\omega$, either $\ell$ decreased and $a$ did not increase, or $\ell$ did not increase and $a$ decreased.
    Repeating this procedure, we obtain a word in normal form after performing at most $(a + \ell)a$ rewritings that do not decrease the length and at most $\ell$ rewritings that decrease the length.
\end{proof}

    \appendix
    
    \section{The complexes $Y_k^3$}\label{app:C2tilde_thickness4}
    In this section we present the $\tilde C_2$-GABs $Y_1^3, \dots,Y_4^3$. They all consist of 12 vertices labeled by $v,w,u_1,\dots u_{10}$, 120 edges labeled by $e_1, \dots, e_{40}$, $f_1 \dots, f_{40}$, $g_1,\dots, g_{40}$, and 160 triangles. The links around the vertices $v,w$ are always the symplectic quadrangle of order 3 and the links around the vertices $u_1, \dots, u_{10}$ are always the complete bipartite graph $K_{4,4}$. The boundary of the edge $g_i$ is always $(v,w)$ and the boundary of the edge $e_i$ is always $(w, u_j)$ where $j = \lceil {{i}/{4}} \rceil$. Let $\tau = (1,2)(3,4)(5,6) \in \Sym(10)$, then the boundary of the edge $f_i \in Y_k^3$ is always $(u_j, v)$ with $j = \tau_k(\lceil {{i}/{4}} \rceil)$. For each complex the assignment $e_i \leftrightarrow f_i$ induces an involutory automorphism $\sigma_k$ of~$Y_k^3$ that swaps the two vertices of special type. In particular, we can apply Lemma \ref{lem:regular_pres} to~$(Y_k^3, \sigma_k)$ and obtain a presentation for the extension $\bar \Gamma_k^3$ of $\Gamma_k^3 = \pi_1(Y_k^3)$ induced by~$\sigma_k$. We also provide these presentations. If we subdivide the complex $S_{\JW}$ along the diagonals from $v_{00}$ to $v_{11}$ and label the new diagonals with $s_i$, we obtain a complex that embeds in each of the complexes $Y_k^3$ via the following assignment.
    \begin{align*}
        a & \mapsto f_5, & \bar a & \mapsto f_6, &
        b &\mapsto f_7, & \bar b &\mapsto f_8, &
        x & \mapsto f_1, & \bar x & \mapsto f_2, &
        y &\mapsto f_3, & \bar y &\mapsto f_4, \\
        a' & \mapsto e_5, &
        \bar a' & \mapsto e_6 &
        b' &\mapsto e_7, & \bar b' &\mapsto e_8, &
        x' & \mapsto e_1, & \bar x' & \mapsto e_2, &
        y' &\mapsto e_3, & \bar y' &\mapsto e_4, \\
        s_i & \mapsto g_i.
    \end{align*}

    The image of these embeddings is always the complex consisting of the first 32 triangles of $Y_k^3$. In particular, the complexes $Y_k^3$ agree on their first 32 triangles and the boundaries of these triangles are the following ones.
    \allowdisplaybreaks
    \begin{align*}
            & (e_{3}, f_{5}, g_{1}), && (e_{5}, f_{1}, g_{1}), && (e_{4}, f_{5}, g_{2}), && (e_{7}, f_{2}, g_{2}), && (e_{1}, f_{5}, g_{3}), && (e_{5}, f_{3}, g_{3}), \\
            & (e_{2}, f_{5}, g_{4}), && (e_{8}, f_{4}, g_{4}), && (e_{3}, f_{6}, g_{5}), && (e_{7}, f_{1}, g_{5}), && (e_{4}, f_{6}, g_{6}), && (e_{6}, f_{2}, g_{6}), \\
            & (e_{1}, f_{6}, g_{7}), && (e_{8}, f_{3}, g_{7}), && (e_{2}, f_{6}, g_{8}), && (e_{6}, f_{4}, g_{8}), && (e_{4}, f_{7}, g_{9}), && (e_{8}, f_{1}, g_{9}), \\
            & (e_{3}, f_{7}, g_{10}), && (e_{8}, f_{2}, g_{10}), && (e_{1}, f_{7}, g_{11}), && (e_{6}, f_{3}, g_{11}), && (e_{2}, f_{7}, g_{12}), && (e_{5}, f_{4}, g_{12}), \\
            & (e_{3}, f_{8}, g_{13}), && (e_{6}, f_{1}, g_{13}), && (e_{4}, f_{8}, g_{14}), && (e_{5}, f_{2}, g_{14}), && (e_{2}, f_{8}, g_{15}), && (e_{7}, f_{3}, g_{15}), \\
            & (e_{1}, f_{8}, g_{16}), && (e_{7}, f_{4}, g_{16}).
    \end{align*}
    Now we indicate the complexes $Y_k^3$ by indicating the boundaries of their remaining 128 triangles. All the properties claimed so far can be checked by hand or with a simple program.
    The boundaries of the remaining triangles in $Y_1^3$ are
    \begin{align*}
        & (e_{9}, f_{13}, g_{1}), && (e_{9}, f_{14}, g_{23}), && (e_{9}, f_{15}, g_{25}), && (e_{9}, f_{16}, g_{36}), && (e_{10}, f_{13}, g_{30}), \\
        & (e_{10}, f_{14}, g_{36}), && (e_{10}, f_{15}, g_{8}), && (e_{10}, f_{16}, g_{19}), && (e_{11}, f_{13}, g_{17}), && (e_{11}, f_{14}, g_{40}), \\
        & (e_{11}, f_{15}, g_{31}), && (e_{11}, f_{16}, g_{6}), && (e_{12}, f_{13}, g_{35}), && (e_{12}, f_{14}, g_{3}), && (e_{12}, f_{15}, g_{33}), \\
        & (e_{12}, f_{16}, g_{37}), && (e_{13}, f_{9}, g_{3}), && (e_{13}, f_{10}, g_{30}), && (e_{13}, f_{11}, g_{17}), && (e_{13}, f_{12}, g_{39}), \\
        & (e_{14}, f_{9}, g_{34}), && (e_{14}, f_{10}, g_{24}), && (e_{14}, f_{11}, g_{40}), && (e_{14}, f_{12}, g_{1}), && (e_{15}, f_{9}, g_{25}), \\
        & (e_{15}, f_{10}, g_{6}), && (e_{15}, f_{11}, g_{22}), && (e_{15}, f_{12}, g_{33}), && (e_{16}, f_{9}, g_{24}), && (e_{16}, f_{10}, g_{26}), \\
        & (e_{16}, f_{11}, g_{8}), && (e_{16}, f_{12}, g_{37}), && (e_{17}, f_{17}, g_{32}), && (e_{17}, f_{18}, g_{11}), && (e_{17}, f_{19}, g_{34}), \\
        & (e_{17}, f_{20}, g_{27}), && (e_{18}, f_{17}, g_{5}), && (e_{18}, f_{18}, g_{33}), && (e_{18}, f_{19}, g_{21}), && (e_{18}, f_{20}, g_{19}), \\
        & (e_{19}, f_{17}, g_{23}), && (e_{19}, f_{18}, g_{21}), && (e_{19}, f_{19}, g_{28}), && (e_{19}, f_{20}, g_{11}), && (e_{20}, f_{17}, g_{27}), \\
        & (e_{20}, f_{18}, g_{26}), && (e_{20}, f_{19}, g_{5}), && (e_{20}, f_{20}, g_{17}), && (e_{21}, f_{21}, g_{18}), && (e_{21}, f_{22}, g_{35}), \\
        & (e_{21}, f_{23}, g_{13}), && (e_{21}, f_{24}, g_{20}), && (e_{22}, f_{21}, g_{39}), && (e_{22}, f_{22}, g_{38}), && (e_{22}, f_{23}, g_{29}), \\
        & (e_{22}, f_{24}, g_{13}), && (e_{23}, f_{21}, g_{7}), && (e_{23}, f_{22}, g_{20}), && (e_{23}, f_{23}, g_{25}), && (e_{23}, f_{24}, g_{26}), \\
        & (e_{24}, f_{21}, g_{29}), && (e_{24}, f_{22}, g_{7}), && (e_{24}, f_{23}, g_{19}), && (e_{24}, f_{24}, g_{40}), && (e_{25}, f_{25}, g_{40}), \\
        & (e_{25}, f_{26}, g_{35}), && (e_{25}, f_{27}, g_{28}), && (e_{25}, f_{28}, g_{15}), && (e_{26}, f_{25}, g_{39}), && (e_{26}, f_{26}, g_{25}), \\
        & (e_{26}, f_{27}, g_{15}), && (e_{26}, f_{28}, g_{32}), && (e_{27}, f_{25}, g_{28}), && (e_{27}, f_{26}, g_{10}), && (e_{27}, f_{27}, g_{30}), \\
        & (e_{27}, f_{28}, g_{22}), && (e_{28}, f_{25}, g_{10}), && (e_{28}, f_{26}, g_{32}), && (e_{28}, f_{27}, g_{31}), && (e_{28}, f_{28}, g_{37}), \\
        & (e_{29}, f_{29}, g_{37}), && (e_{29}, f_{30}, g_{2}), && (e_{29}, f_{31}, g_{29}), && (e_{29}, f_{32}, g_{23}), && (e_{30}, f_{29}, g_{12}), \\
        &(e_{30}, f_{30}, g_{27}), && (e_{30}, f_{31}, g_{22}), && (e_{30}, f_{32}, g_{29}), && (e_{31}, f_{29}, g_{20}), && (e_{31}, f_{30}, g_{31}), \\
        & (e_{31}, f_{31}, g_{21}), && (e_{31}, f_{32}, g_{12}), && (e_{32}, f_{29}, g_{34}), && (e_{32}, f_{30}, g_{20}), && (e_{32}, f_{31}, g_{2}), \\
        & (e_{32}, f_{32}, g_{30}), && (e_{33}, f_{33}, g_{17}), && (e_{33}, f_{34}, g_{23}), && (e_{33}, f_{35}, g_{18}), && (e_{33}, f_{36}, g_{4}), \\
        & (e_{34}, f_{33}, g_{34}), && (e_{34}, f_{34}, g_{33}), && (e_{34}, f_{35}, g_{4}), && (e_{34}, f_{36}, g_{38}), && (e_{35}, f_{33}, g_{18}), \\
        & (e_{35}, f_{34}, g_{14}), && (e_{35}, f_{35}, g_{32}), && (e_{35}, f_{36}, g_{19}), && (e_{36}, f_{33}, g_{14}), && (e_{36}, f_{34}, g_{38}), \\
        & (e_{36}, f_{35}, g_{26}), && (e_{36}, f_{36}, g_{28}), && (e_{37}, f_{37}, g_{21}), && (e_{37}, f_{38}, g_{24}), && (e_{37}, f_{39}, g_{9}), \\
        & (e_{37}, f_{40}, g_{39}), && (e_{38}, f_{37}, g_{36}), && (e_{38}, f_{38}, g_{38}), && (e_{38}, f_{39}, g_{31}), && (e_{38}, f_{40}, g_{16}), \\
        & (e_{39}, f_{37}, g_{16}), && (e_{39}, f_{38}, g_{22}), && (e_{39}, f_{39}, g_{18}), && (e_{39}, f_{40}, g_{24}), && (e_{40}, f_{37}, g_{35}), \\
        & (e_{40}, f_{38}, g_{9}), && (e_{40}, f_{39}, g_{36}), && (e_{40}, f_{40}, g_{27}).
    \end{align*}
    The boundaries of the remaining triangles in $Y_2^3$ are
    \begin{align*}
        & (e_{9}, f_{13}, g_{8}), && (e_{9}, f_{14}, g_{37}), && (e_{9}, f_{15}, g_{24}), && (e_{9}, f_{16}, g_{26}), && (e_{10}, f_{13}, g_{40}), \\
        & (e_{10}, f_{14}, g_{1}), && (e_{10}, f_{15}, g_{34}), && (e_{10}, f_{16}, g_{24}), && (e_{11}, f_{13}, g_{22}), && (e_{11}, f_{14}, g_{33}), \\
        & (e_{11}, f_{15}, g_{25}), && (e_{11}, f_{16}, g_{6}), && (e_{12}, f_{13}, g_{17}), && (e_{12}, f_{14}, g_{39}), && (e_{12}, f_{15}, g_{3}), \\
        & (e_{12}, f_{16}, g_{30}), && (e_{13}, f_{9}, g_{6}), && (e_{13}, f_{10}, g_{40}), && (e_{13}, f_{11}, g_{31}), && (e_{13}, f_{12}, g_{17}), \\
        & (e_{14}, f_{9}, g_{37}), && (e_{14}, f_{10}, g_{3}), && (e_{14}, f_{11}, g_{33}), && (e_{14}, f_{12}, g_{35}), && (e_{15}, f_{9}, g_{36}), \\
        & (e_{15}, f_{10}, g_{23}), && (e_{15}, f_{11}, g_{25}), && (e_{15}, f_{12}, g_{1}), && (e_{16}, f_{9}, g_{19}), && (e_{16}, f_{10}, g_{36}), \\
        & (e_{16}, f_{11}, g_{8}), && (e_{16}, f_{12}, g_{30}), && (e_{17}, f_{17}, g_{30}), && (e_{17}, f_{18}, g_{20}), && (e_{17}, f_{19}, g_{34}), \\
        & (e_{17}, f_{20}, g_{2}), && (e_{18}, f_{17}, g_{20}), && (e_{18}, f_{18}, g_{21}), && (e_{18}, f_{19}, g_{12}), && (e_{18}, f_{20}, g_{31}), \\
        & (e_{19}, f_{17}, g_{23}), && (e_{19}, f_{18}, g_{2}), && (e_{19}, f_{19}, g_{37}), && (e_{19}, f_{20}, g_{29}), && (e_{20}, f_{17}, g_{12}), \\
        & (e_{20}, f_{18}, g_{22}), && (e_{20}, f_{19}, g_{29}), && (e_{20}, f_{20}, g_{27}), && (e_{21}, f_{21}, g_{18}), && (e_{21}, f_{22}, g_{13}), \\
        & (e_{21}, f_{23}, g_{20}), && (e_{21}, f_{24}, g_{35}), && (e_{22}, f_{21}, g_{7}), && (e_{22}, f_{22}, g_{40}), && (e_{22}, f_{23}, g_{19}), \\
        & (e_{22}, f_{24}, g_{29}), && (e_{23}, f_{21}, g_{20}), && (e_{23}, f_{22}, g_{26}), && (e_{23}, f_{23}, g_{25}), && (e_{23}, f_{24}, g_{7}), \\
        & (e_{24}, f_{21}, g_{39}), && (e_{24}, f_{22}, g_{29}), && (e_{24}, f_{23}, g_{13}), && (e_{24}, f_{24}, g_{38}), && (e_{25}, f_{25}, g_{21}), \\
        & (e_{25}, f_{26}, g_{24}), && (e_{25}, f_{27}, g_{9}), && (e_{25}, f_{28}, g_{39}), && (e_{26}, f_{25}, g_{36}), && (e_{26}, f_{26}, g_{38}), \\
        & (e_{26}, f_{27}, g_{31}), && (e_{26}, f_{28}, g_{16}), && (e_{27}, f_{25}, g_{16}), && (e_{27}, f_{26}, g_{22}), && (e_{27}, f_{27}, g_{18}), \\
        & (e_{27}, f_{28}, g_{24}), && (e_{28}, f_{25}, g_{35}), && (e_{28}, f_{26}, g_{9}), && (e_{28}, f_{27}, g_{36}), && (e_{28}, f_{28}, g_{27}), \\
        & (e_{29}, f_{29}, g_{40}), && (e_{29}, f_{30}, g_{35}), && (e_{29}, f_{31}, g_{28}), && (e_{29}, f_{32}, g_{15}), && (e_{30}, f_{29}, g_{39}), \\
        & (e_{30}, f_{30}, g_{25}), && (e_{30}, f_{31}, g_{15}), && (e_{30}, f_{32}, g_{32}), && (e_{31}, f_{29}, g_{28}), && (e_{31}, f_{30}, g_{10}), \\
        & (e_{31}, f_{31}, g_{30}), && (e_{31}, f_{32}, g_{22}), && (e_{32}, f_{29}, g_{10}), && (e_{32}, f_{30}, g_{32}), && (e_{32}, f_{31}, g_{31}), \\
        & (e_{32}, f_{32}, g_{37}), && (e_{33}, f_{33}, g_{28}), && (e_{33}, f_{34}, g_{26}), && (e_{33}, f_{35}, g_{14}), && (e_{33}, f_{36}, g_{38}), \\
        & (e_{34}, f_{33}, g_{19}), && (e_{34}, f_{34}, g_{32}), && (e_{34}, f_{35}, g_{18}), && (e_{34}, f_{36}, g_{14}), && (e_{35}, f_{33}, g_{4}), \\
        & (e_{35}, f_{34}, g_{18}), && (e_{35}, f_{35}, g_{17}), && (e_{35}, f_{36}, g_{23}), && (e_{36}, f_{33}, g_{38}), && (e_{36}, f_{34}, g_{4}), \\
        & (e_{36}, f_{35}, g_{34}), && (e_{36}, f_{36}, g_{33}), && (e_{37}, f_{37}, g_{33}), && (e_{37}, f_{38}, g_{21}), && (e_{37}, f_{39}, g_{5}), \\
        & (e_{37}, f_{40}, g_{19}), && (e_{38}, f_{37}, g_{21}), && (e_{38}, f_{38}, g_{28}), && (e_{38}, f_{39}, g_{23}), && (e_{38}, f_{40}, g_{11}), \\
        & (e_{39}, f_{37}, g_{11}), && (e_{39}, f_{38}, g_{34}), && (e_{39}, f_{39}, g_{32}), && (e_{39}, f_{40}, g_{27}), && (e_{40}, f_{37}, g_{26}), \\
        & (e_{40}, f_{38}, g_{5}), && (e_{40}, f_{39}, g_{27}), && (e_{40}, f_{40}, g_{17}).
    \end{align*}
    The boundaries of the remaining triangles in $Y_3^3$ are
    \begin{align*}
        & (e_{9}, f_{13}, g_{4}), && (e_{9}, f_{14}, g_{33}), && (e_{9}, f_{15}, g_{38}), && (e_{9}, f_{16}, g_{34}), && (e_{10}, f_{13}, g_{28}), \\
        & (e_{10}, f_{14}, g_{11}), && (e_{10}, f_{15}, g_{23}), && (e_{10}, f_{16}, g_{21}), && (e_{11}, f_{13}, g_{26}), && (e_{11}, f_{14}, g_{27}), \\
        & (e_{11}, f_{15}, g_{17}), && (e_{11}, f_{16}, g_{5}), && (e_{12}, f_{13}, g_{18}), && (e_{12}, f_{14}, g_{19}), && (e_{12}, f_{15}, g_{14}), \\
        & (e_{12}, f_{16}, g_{32}), && (e_{13}, f_{9}, g_{14}), && (e_{13}, f_{10}, g_{28}), && (e_{13}, f_{11}, g_{26}), && (e_{13}, f_{12}, g_{38}), \\
        & (e_{14}, f_{9}, g_{33}), && (e_{14}, f_{10}, g_{5}), && (e_{14}, f_{11}, g_{21}), && (e_{14}, f_{12}, g_{19}), && (e_{15}, f_{9}, g_{18}), \\
        & (e_{15}, f_{10}, g_{23}), && (e_{15}, f_{11}, g_{17}), && (e_{15}, f_{12}, g_{4}), && (e_{16}, f_{9}, g_{34}), && (e_{16}, f_{10}, g_{27}), \\
        & (e_{16}, f_{11}, g_{11}), && (e_{16}, f_{12}, g_{32}), && (e_{17}, f_{17}, g_{35}), && (e_{17}, f_{18}, g_{20}), && (e_{17}, f_{19}, g_{13}), \\
        & (e_{17}, f_{20}, g_{18}), && (e_{18}, f_{17}, g_{20}), && (e_{18}, f_{18}, g_{26}), && (e_{18}, f_{19}, g_{25}), && (e_{18}, f_{20}, g_{7}), \\
        & (e_{19}, f_{17}, g_{7}), && (e_{19}, f_{18}, g_{40}), && (e_{19}, f_{19}, g_{19}), && (e_{19}, f_{20}, g_{29}), && (e_{20}, f_{17}, g_{38}), \\
        & (e_{20}, f_{18}, g_{13}), && (e_{20}, f_{19}, g_{29}), && (e_{20}, f_{20}, g_{39}), && (e_{21}, f_{21}, g_{23}), && (e_{21}, f_{22}, g_{29}), \\
        & (e_{21}, f_{23}, g_{37}), && (e_{21}, f_{24}, g_{2}), && (e_{22}, f_{21}, g_{29}), && (e_{22}, f_{22}, g_{22}), && (e_{22}, f_{23}, g_{12}), \\
        & (e_{22}, f_{24}, g_{27}), && (e_{23}, f_{21}, g_{30}), && (e_{23}, f_{22}, g_{2}), && (e_{23}, f_{23}, g_{34}), && (e_{23}, f_{24}, g_{20}), \\
        & (e_{24}, f_{21}, g_{12}), && (e_{24}, f_{22}, g_{21}), && (e_{24}, f_{23}, g_{20}), && (e_{24}, f_{24}, g_{31}), && (e_{25}, f_{25}, g_{39}), \\
        & (e_{25}, f_{26}, g_{15}), && (e_{25}, f_{27}, g_{32}), && (e_{25}, f_{28}, g_{25}), && (e_{26}, f_{25}, g_{10}), && (e_{26}, f_{26}, g_{22}), \\
        & (e_{26}, f_{27}, g_{30}), && (e_{26}, f_{28}, g_{28}), && (e_{27}, f_{25}, g_{32}), && (e_{27}, f_{26}, g_{37}), && (e_{27}, f_{27}, g_{31}), \\
        & (e_{27}, f_{28}, g_{10}), && (e_{28}, f_{25}, g_{40}), && (e_{28}, f_{26}, g_{28}), && (e_{28}, f_{27}, g_{15}), && (e_{28}, f_{28}, g_{35}), \\
        & (e_{29}, f_{29}, g_{39}), && (e_{29}, f_{30}, g_{17}), && (e_{29}, f_{31}, g_{3}), && (e_{29}, f_{32}, g_{30}), && (e_{30}, f_{29}, g_{17}), \\
        & (e_{30}, f_{30}, g_{31}), && (e_{30}, f_{31}, g_{40}), && (e_{30}, f_{32}, g_{6}), && (e_{31}, f_{29}, g_{1}), && (e_{31}, f_{30}, g_{25}), \\
        & (e_{31}, f_{31}, g_{23}), && (e_{31}, f_{32}, g_{36}), && (e_{32}, f_{29}, g_{37}), && (e_{32}, f_{30}, g_{8}), && (e_{32}, f_{31}, g_{24}), \\
        & (e_{32}, f_{32}, g_{26}), && (e_{33}, f_{33}, g_{22}), && (e_{33}, f_{34}, g_{33}), && (e_{33}, f_{35}, g_{6}), && (e_{33}, f_{36}, g_{25}), \\
        & (e_{34}, f_{33}, g_{33}), && (e_{34}, f_{34}, g_{35}), && (e_{34}, f_{35}, g_{37}), && (e_{34}, f_{36}, g_{3}), && (e_{35}, f_{33}, g_{8}),\\
        & (e_{35}, f_{34}, g_{30}), && (e_{35}, f_{35}, g_{19}), && (e_{35}, f_{36}, g_{36}), && (e_{36}, f_{33}, g_{40}), && (e_{36}, f_{34}, g_{1}), \\
        & (e_{36}, f_{35}, g_{24}), && (e_{36}, f_{36}, g_{34}), && (e_{37}, f_{37}, g_{31}), && (e_{37}, f_{38}, g_{16}), && (e_{37}, f_{39}, g_{36}), \\
        & (e_{37}, f_{40}, g_{38}), && (e_{38}, f_{37}, g_{9}), && (e_{38}, f_{38}, g_{35}), && (e_{38}, f_{39}, g_{27}), && (e_{38}, f_{40}, g_{36}), \\
        & (e_{39}, f_{37}, g_{24}), && (e_{39}, f_{38}, g_{21}), && (e_{39}, f_{39}, g_{39}), && (e_{39}, f_{40}, g_{9}), && (e_{40}, f_{37}, g_{18}), \\
        & (e_{40}, f_{38}, g_{24}), && (e_{40}, f_{39}, g_{16}), && (e_{40}, f_{40}, g_{22}).
    \end{align*}
    The boundaries of the remaining triangles in $Y_4^3$ are
    \begin{align*}
        & (e_{9}, f_{13}, g_{12}), && (e_{9}, f_{14}, g_{31}), && (e_{9}, f_{15}, g_{21}), && (e_{9}, f_{16}, g_{20}), && (e_{10}, f_{13}, g_{19}), \\
        & (e_{10}, f_{14}, g_{7}), && (e_{10}, f_{15}, g_{29}), && (e_{10}, f_{16}, g_{40}), && (e_{11}, f_{13}, g_{30}), && (e_{11}, f_{14}, g_{20}), \\
        & (e_{11}, f_{15}, g_{2}), && (e_{11}, f_{16}, g_{34}), && (e_{12}, f_{13}, g_{29}), && (e_{12}, f_{14}, g_{38}), && (e_{12}, f_{15}, g_{39}), \\
        & (e_{12}, f_{16}, g_{13}), && (e_{13}, f_{9}, g_{2}), && (e_{13}, f_{10}, g_{23}), && (e_{13}, f_{11}, g_{37}), && (e_{13}, f_{12}, g_{29}), \\
        & (e_{14}, f_{9}, g_{35}), && (e_{14}, f_{10}, g_{13}), && (e_{14}, f_{11}, g_{20}), && (e_{14}, f_{12}, g_{18}), && (e_{15}, f_{9}, g_{27}), \\
        & (e_{15}, f_{10}, g_{29}), && (e_{15}, f_{11}, g_{12}), && (e_{15}, f_{12}, g_{22}), && (e_{16}, f_{9}, g_{20}), && (e_{16}, f_{10}, g_{25}), \\
        & (e_{16}, f_{11}, g_{26}), && (e_{16}, f_{12}, g_{7}), && (e_{17}, f_{17}, g_{33}), && (e_{17}, f_{18}, g_{3}), && (e_{17}, f_{19}, g_{35}), \\
        & (e_{17}, f_{20}, g_{37}), && (e_{18}, f_{17}, g_{1}), && (e_{18}, f_{18}, g_{36}), && (e_{18}, f_{19}, g_{25}), && (e_{18}, f_{20}, g_{23}), \\
        & (e_{19}, f_{17}, g_{31}), && (e_{19}, f_{18}, g_{40}), && (e_{19}, f_{19}, g_{17}), && (e_{19}, f_{20}, g_{6}), && (e_{20}, f_{17}, g_{30}), \\
        & (e_{20}, f_{18}, g_{19}), && (e_{20}, f_{19}, g_{8}), && (e_{20}, f_{20}, g_{36}), && (e_{21}, f_{21}, g_{24}), && (e_{21}, f_{22}, g_{8}), \\
        & (e_{21}, f_{23}, g_{37}), && (e_{21}, f_{24}, g_{26}), && (e_{22}, f_{21}, g_{6}), && (e_{22}, f_{22}, g_{33}), && (e_{22}, f_{23}, g_{22}),\\
        & (e_{22}, f_{24}, g_{25}), && (e_{23}, f_{21}, g_{30}), && (e_{23}, f_{22}, g_{39}), && (e_{23}, f_{23}, g_{17}), && (e_{23}, f_{24}, g_{3}), \\
        & (e_{24}, f_{21}, g_{34}), && (e_{24}, f_{22}, g_{40}), && (e_{24}, f_{23}, g_{1}), && (e_{24}, f_{24}, g_{24}), && (e_{25}, f_{25}, g_{36}), \\
        & (e_{25}, f_{26}, g_{31}), &&  (e_{25}, f_{27}, g_{38}), && (e_{25}, f_{28}, g_{16}), && (e_{26}, f_{25}, g_{35}), && (e_{26}, f_{26}, g_{36}), \\
        & (e_{26}, f_{27}, g_{9}), && (e_{26}, f_{28}, g_{27}), &&  (e_{27}, f_{25}, g_{18}), && (e_{27}, f_{26}, g_{16}), && (e_{27}, f_{27}, g_{24}), \\
        & (e_{27}, f_{28}, g_{22}), && (e_{28}, f_{25}, g_{9}), && (e_{28}, f_{26}, g_{21}), && (e_{28}, f_{27}, g_{39}), && (e_{28}, f_{28}, g_{24}), \\
        & (e_{29}, f_{29}, g_{17}), && (e_{29}, f_{30}, g_{23}), && (e_{29}, f_{31}, g_{18}), && (e_{29}, f_{32}, g_{4}), &&  (e_{30}, f_{29}, g_{19}), \\
        & (e_{30}, f_{30}, g_{32}), && (e_{30}, f_{31}, g_{14}), && (e_{30}, f_{32}, g_{18}), && (e_{31}, f_{29}, g_{38}), && (e_{31}, f_{30}, g_{4}), \\
        & (e_{31}, f_{31}, g_{33}), && (e_{31}, f_{32}, g_{34}), && (e_{32}, f_{29}, g_{14}), && (e_{32}, f_{30}, g_{38}), && (e_{32}, f_{31}, g_{26}), \\
        & (e_{32}, f_{32}, g_{28}), && (e_{33}, f_{33}, g_{33}), && (e_{33}, f_{34}, g_{19}), && (e_{33}, f_{35}, g_{5}), && (e_{33}, f_{36}, g_{21}), \\
        & (e_{34}, f_{33}, g_{23}), && (e_{34}, f_{34}, g_{28}), &&  (e_{34}, f_{35}, g_{21}), && (e_{34}, f_{36}, g_{11}), && (e_{35}, f_{33}, g_{11}), \\
        & (e_{35}, f_{34}, g_{27}), && (e_{35}, f_{35}, g_{32}), && (e_{35}, f_{36}, g_{34}), && (e_{36}, f_{33}, g_{27}), && (e_{36}, f_{34}, g_{5}), \\
        & (e_{36}, f_{35}, g_{26}), && (e_{36}, f_{36}, g_{17}), && (e_{37}, f_{37}, g_{32}), && (e_{37}, f_{38}, g_{10}), && (e_{37}, f_{39}, g_{37}), \\
        & (e_{37}, f_{40}, g_{31}), && (e_{38}, f_{37}, g_{15}), && (e_{38}, f_{38}, g_{32}), && (e_{38}, f_{39}, g_{39}), && (e_{38}, f_{40}, g_{25}), \\
        & (e_{39}, f_{37}, g_{30}), && (e_{39}, f_{38}, g_{22}), && (e_{39}, f_{39}, g_{28}), && (e_{39}, f_{40}, g_{10}), && (e_{40}, f_{37}, g_{35}), \\
        &(e_{40}, f_{38}, g_{40}), && (e_{40}, f_{39}, g_{15}), && (e_{40}, f_{40}, g_{28}).
    \end{align*}
    With Lemma \ref{lem:regular_pres} we computed presentations for the extensions $\bar \Gamma_k^3$ (and extracted a sufficiently big set of relations to present the group).
    In particular, the Cayley graphs of these are the subgraphs of the 1-skeletons of the buildings $X_k^3 = \tilde Y_k^3$ generated by special vertices.
    The generators of these presentations are $g_1, \dots, g_{40}$ and since the complexes $Y_k^3$ share the first 32 triangles, the four presentation share the following relations.
    \begin{align*}
        & g_1 g_3, && g_2 g_{12}, && g_4 g_{14}, && g_5 g_{11}, && g_6 g_8, && g_7 g_{13}, \\
        & g_9 g_{16}, && g_{10} g_{15}, && g_1^2 g_{15} g_{11}, && g_1 g_7 g_{14} g_{12}, && g_1 g_{15} g_{12} g_{14}, && g_2 g_8^2 g_{14}, \\
        & g_2 g_9 g_{12} g_4, && g_2 g_{11} g_9 g_{15}, && g_5 g_{10} g_7 g_{16}, && g_5 g_{10} g_{11} g_7, && g_5 g_{15} g_9 g_8.
    \end{align*}
The remaining relations of the presentation for $\bar \Gamma_1^3$ are 
\begin{align*}
    & g_{17}^2, && g_{18}^2, && g_{19} g_{26}, && g_{20} g_{29}, && g_{21}^2, \\
    & g_{22} g_{31}, && g_{23} g_{34}, && g_{24} g_{36}, && g_{25}^2, && g_{27}^2, \\
    & g_{28}^2, && g_{30}^2, && g_{32}^2, && g_{33}^2, && g_{35} g_{39}, \\
    & g_{37}^2, && g_{38}^2, && g_{40}^2, && g_1 g_{23} g_{25} g_{33}, && g_1 g_{24} g_6 g_{17}, \\
    & g_1 g_{24} g_{19} g_{30}, && g_1 g_{25} g_8 g_{30}, && g_1 g_{34} g_{36} g_{30}, && g_1 g_{40} g_{17} g_{35}, && g_2 g_{20} g_{22} g_{27}, \\
    & g_2 g_{23} g_{12} g_{31}, && g_2 g_{30} g_{23} g_{20}, && g_2 g_{34} g_{30} g_{29}, && g_2 g_{37} g_{12} g_{27}, && g_4 g_{17} g_{14} g_{28}, \\
    & g_4 g_{23} g_{14} g_{19}, && g_4 g_{33} g_{14} g_{32}, && g_4 g_{34} g_{38} g_{28}, && g_5 g_{19} g_{11} g_{23}, && g_5 g_{21} g_{34} g_{32}, \\
    & g_5 g_{26} g_{17} g_{27}, && g_5 g_{27} g_{23} g_{28}, && g_6 g_{22} g_{33} g_{37}, && g_6 g_{33} g_{37} g_{19}, && g_6 g_{40} g_{36} g_{26}, \\
    & g_7 g_{19} g_{40} g_{20}, && g_7 g_{20} g_{18} g_{39}, && g_7 g_{26} g_{29} g_{38}, && g_9 g_{21} g_{16} g_{18}, && g_9 g_{24} g_{31} g_{38}, \\
    & g_9 g_{27} g_{24} g_{31}, && g_9 g_{35} g_{16} g_{22}, && g_9 g_{39} g_{21} g_{36}, && g_{10} g_{28} g_{39} g_{25}, && g_{10} g_{28} g_{40} g_{39}, \\
    & g_{10} g_{30} g_{31} g_{32}, && g_{10} g_{31} g_{37} g_{32}, && g_{10} g_{37} g_{32} g_{35}, && g_{17} g_{19} g_{33} g_{26}, && g_{17} g_{34} g_{33} g_{23}, \\
    & g_{18} g_{26} g_{38} g_{23}, && g_{18} g_{29} g_{40} g_{20}, && g_{19} g_{20} g_{39} g_{20}, && g_{20} g_{25} g_{29} g_{38}, && g_{21} g_{22} g_{27} g_{31}, \\
    & g_{22} g_{36} g_{39} g_{36}, && g_{23} g_{28} g_{34} g_{32}, && g_{25} g_{32} g_{37} g_{32}.
\end{align*}
    The remaining relations of the presentation for $\bar \Gamma_2^3$ are
\begin{align*}
    & g_{17}^2, && g_{18}^2, && g_{19} g_{26}, && g_{20}^2, && g_{21}^2, \\
    & g_{22} g_{31}, && g_{23} g_{34}, && g_{24} g_{36}, && g_{25}^2, \\
    & g_{27}^2, && g_{28}^2, && g_{29}^2, && g_{30}^2, && g_{32}^2, \\
    & g_{33}^2, && g_{35} g_{39}, && g_{37}^2, && g_{38}^2, && g_{40}^2, \\
    & g_1 g_{25} g_{31} g_{17}, && g_1 g_{40} g_8 g_{37}, && g_2 g_{20} g_{22} g_{27}, && g_2 g_{23} g_{12} g_{22}, && g_2 g_{34} g_{30} g_{20}, \\
    & g_2 g_{37} g_{29} g_{31}, && g_4 g_{17} g_{18} g_{26}, && g_4 g_{17} g_{34} g_{38}, && g_4 g_{23} g_{14} g_{19}, && g_4 g_{34} g_{38} g_{28}, \\
    & g_5 g_{19} g_{11} g_{23}, && g_5 g_{26} g_{17} g_{27}, && g_5 g_{27} g_{23} g_{28}, && g_5 g_{33} g_{11} g_{32}, && g_5 g_{33} g_{26} g_{27}, \\
    & g_6 g_{17} g_{35} g_{37}, && g_6 g_{22} g_8 g_{26}, && g_6 g_{25} g_{34} g_{36}, && g_6 g_{31} g_{17} g_{30}, && g_7 g_{19} g_{13} g_{39}, \\
    & g_7 g_{20} g_{39} g_{38}, && g_7 g_{25} g_{19} g_{29}, && g_7 g_{26} g_{20} g_{18}, && g_7 g_{40} g_{29} g_{35}, && g_9 g_{24} g_{18} g_{31}, \\
    & g_9 g_{24} g_{31} g_{38}, && g_9 g_{27} g_{24} g_{31}, && g_9 g_{35} g_{16} g_{22}, && g_9 g_{36} g_{38} g_{22}, && g_9 g_{39} g_{21} g_{36}, \\
    & g_{10} g_{28} g_{40} g_{39}, && g_{10} g_{30} g_{28} g_{39}, && g_{10} g_{31} g_{32} g_{25}, && g_{10} g_{32} g_{35} g_{40}, && g_{10} g_{37} g_{22} g_{28}, \\
    & g_{10} g_{37} g_{32} g_{35}, && g_{17} g_{18} g_{32} g_{18}, && g_{17} g_{22} g_{33} g_{39}, && g_{17} g_{34} g_{33} g_{23}, && g_{18} g_{20} g_{25} g_{20}, \\
    & g_{19} g_{21} g_{34} g_{27}, && g_{19} g_{24} g_{23} g_{24}, && g_{20} g_{21} g_{20} g_{30}, && g_{22} g_{25} g_{34} g_{40}, && g_{27} g_{29} g_{37} g_{29}.
\end{align*}
    The remaining relations of the presentation for $\bar \Gamma_3^3$ are
    \begin{align*}
        & g_{17}^2, && g_{18} g_{38}, && g_{19}^2, && g_{20}^2, && g_{21} g_{27}, \\
        & g_{22}^2, && g_{23}^2, && g_{24} g_{36}, && g_{25} g_{40}, && g_{26}^2, \\
        & g_{28}^2, && g_{29}^2, && g_{30} g_{37}, && g_{31}^2, && g_{32}^2, \\
        & g_{33}^2, && g_{34}^2, && g_{35}^2, && g_{39}^2, && g_1 g_{24} g_{30} g_{39}, \\
        & g_1 g_{34} g_{25} g_{33}, && g_1 g_{36} g_{37} g_{35}, && g_1 g_{40} g_8 g_{30}, && g_1 g_{40} g_{17} g_{39}, && g_2 g_{23} g_{12} g_{31}, \\
        & g_2 g_{29} g_{21} g_{31}, && g_2 g_{30} g_{34} g_{20}, && g_2 g_{34} g_{37} g_{29}, && g_2 g_{37} g_{29} g_{22}, && g_4 g_{17} g_{21} g_{19}, \\
        & g_4 g_{18} g_{17} g_{26}, && g_4 g_{34} g_{32} g_{38}, && g_5 g_{19} g_{32} g_{21}, && g_5 g_{21} g_{11} g_{27}, && g_5 g_{27} g_{17} g_{23}, \\
        & g_5 g_{27} g_{26} g_{28}, && g_6 g_{17} g_{37} g_{26}, && g_6 g_{25} g_{23} g_{24}, && g_6 g_{31} g_{25} g_{24}, && g_6 g_{40} g_{36} g_{19}, \\
        & g_7 g_{20} g_{38} g_{39}, && g_7 g_{25} g_{26} g_{20}, && g_7 g_{29} g_{18} g_{35}, && g_7 g_{40} g_{19} g_{29}, && g_9 g_{24} g_9 g_{36}, \\
        & g_9 g_{24} g_{22} g_{38}, && g_9 g_{24} g_{38} g_{31}, && g_9 g_{27} g_{24} g_{22}, && g_9 g_{39} g_{27} g_{24}, && g_{10} g_{22} g_{15} g_{39}, \\
        & g_{10} g_{28} g_{35} g_{25}, && g_{10} g_{30} g_{15} g_{40}, && g_{10} g_{31} g_{32} g_{40}, && g_{10} g_{32} g_{39} g_{40}, &&  g_{10} g_{32} g_{40} g_{35}, \\
        & g_{17} g_{21} g_{33} g_{18}, && g_{18} g_{20} g_{40} g_{29}, && g_{18} g_{32} g_{21} g_{28}, && g_{19} g_{24} g_{34} g_{36}, && g_{19} g_{25} g_{26} g_{40}, \\
        & g_{20} g_{27} g_{29} g_{30}, && g_{22} g_{33} g_{35} g_{33}, && g_{22} g_{36} g_{35} g_{24}, && g_{24} g_{25} g_{33} g_{30}.
    \end{align*}
    Then remaining relations of the presentation for $\bar \Gamma_4^3$ are
    \begin{align*}
        & g_{17}^2, && g_{18} g_{38}, && g_{19} g_{23}, && g_{20}^2, && g_{21} g_{27}, \\
        & g_{22} g_{39}, && g_{24}^2, && g_{25} g_{40}, && g_{26} g_{34}, && g_{28}^2, \\
        & g_{29}^2, && g_{30} g_{37}, && g_{31} g_{35}, && g_{32}^2, && g_{33}^2, \\
        & g_{36}^2, && g_1 g_{19} g_{36} g_{37}, && g_1 g_{24} g_{26} g_{30}, && g_1 g_{25} g_{33} g_{39}, && g_1 g_{36} g_3 g_{33}, \\
        & g_1 g_{40} g_8 g_{37}, && g_2 g_{26} g_{20} g_{27}, && g_2 g_{29} g_{22} g_{21}, && g_4 g_{17} g_{38} g_{26}, && g_4 g_{26} g_{18} g_{32}, \\
        & g_4 g_{33} g_{26} g_{18}, && g_4 g_{38} g_{26} g_{28}, && g_5 g_{17} g_{21} g_{23}, && g_5 g_{21} g_{11} g_{21}, && g_5 g_{21} g_{23} g_{28}, \\
        & g_5 g_{23} g_{27} g_{32}, && g_5 g_{27} g_{17} g_{34}, && g_5 g_{33} g_{11} g_{32}, && g_6 g_{17} g_{35} g_{30}, && g_6 g_{25} g_{36} g_{19}, \\
        & g_6 g_{39} g_{37} g_{24}, && g_7 g_{23} g_{30} g_{20}, && g_7 g_{40} g_{23} g_{29}, && g_9 g_{21} g_{16} g_{18}, && g_9 g_{22} g_{38} g_{36}, \\
        & g_9 g_{27} g_{36} g_{31}, && g_9 g_{36} g_{21} g_{22}, && g_{10} g_{28} g_{39} g_{40}, && g_{10} g_{30} g_{39} g_{32}, && g_{10} g_{32} g_{30} g_{39}, \\
        & g_{10} g_{37} g_{32} g_{35}, && g_{17} g_{22} g_{33} g_{39}, && g_{17} g_{25} g_{36} g_{40}, && g_{17} g_{38} g_{33} g_{18}, && g_{18} g_{19} g_{38} g_{34}, \\
        & g_{18} g_{20} g_{37} g_{29}, && g_{18} g_{24} g_{38} g_{36}, && g_{18} g_{31} g_{27} g_{39}, && g_{19} g_{25} g_{34} g_{37}, && g_{19} g_{29} g_{39} g_{29}, \\
        & g_{20} g_{26} g_{20} g_{35}, && g_{21} g_{24} g_{27} g_{36}, && g_{22} g_{37} g_{35} g_{25}, && g_{22} g_{40} g_{26} g_{30}, && g_{24} g_{26} g_{24} g_{34}, \\
        & g_{25} g_{32} g_{40} g_{28}.
    \end{align*}

    Recall that by Proposition \ref{prop:JanzenWise} the homotopy classes of the following loops lie in the finite residual of $\pi_1(S_{\JW})$.
    \begin{align*}
      x^{-1} * \bar x * x^{-1} * (\bar y * y^{-1})^2 * x * \bar x^{-1} * x * (\bar y^{-1} * y)^2, \\
      y^{-1} * \bar y * y^{-1} * (\bar x * x^{-1})^2 * y * \bar y^{-1} * y * (\bar x^{-1} * x)^2.
    \end{align*}
    By tracing the embedding we find that these loops correspond to the following loops in $\Gamma_k^3$.
    \begin{align*}
        f_1^{-1} * f_2 * f_1^{-1} * (f_4 * f_3^{-1})^2 * f_1 * f_2^{-1} * f_1 * (f_4^{-1} * f_3)^2,\\
        f_3^{-1} * f_4 * f_3^{-1} * (f_2 * f_1^{-1})^2 * f_3 * f_4^{-1} * f_3 * (f_2^{-1} * f_1)^2.
    \end{align*}
    And these loops are homotopy equivalent to the following ones.
    \begin{align*}
        \eta_1 \defeq g_1*g_{14}^{-1}*g_1*(g_{12}^{-1}*g_3)^2*g_1^{-1}*g_{14}*g_1^{-1}*(g_{12}*g_3^{-1})^2,\\
        \eta_2 \defeq g_3*g_{12}^{-1}*g_3*(g_{14}^{-1}*g_1)^2*g_3^{-1}*g_{12}*g_3^{-1}*(g_{14}*g_1^{-1})^2.
    \end{align*}
    We can interpret $\eta_1$ and $\eta_2$ as elements in $\bar \Gamma_k^3$ and deduce that they lie in the finite residual. With the help of a computer we compute the indexes of the normal closures $\langle \langle \eta_1, \eta_2\rangle\rangle$ in the $\bar \Gamma_k^3$ and they are 8, 8, 16, 16, respectively. In particular, $\langle \langle \eta_1, \eta_2\rangle\rangle$ is the finite residual in every case.
    Using the algorithm in Section \ref{sec:normal_forms}, we compute balls in the Cayley graphs for the geometric presentation of $\Gamma_k^3$ and study their automorphism group. At the time writing an ordinary office machine is not capable of computing the automorphism group of the 4-ball in a reasonable time.
    However, we can compute the automorphism group of the 3-ball, and apply the following lemma.
    \begin{lemma}\label{lem:discreteness_criterion}
        Let $Z$ be a Cayley graph for a group $\Gamma$ with generating set $S$. For $g \in Z, r\in \NN$ denote the ball of radius $r$ centered at $g$ with $B_r(g)$ and the automorphism group of $B_r(g)$ that fixes the center with $A^r(g)$. We denote the image of the projection of $A^r(g)$ to the automorphism group of the $m$-ball centered at $g$ with $A_m^r(g)$ for $m\leq r$. Define 
        \[
            E = \{g \in B_1(1) \mid \Stab(A_1^3, g) = 1\} \quad \text{and} \quad F = B_1(1) \cup \bigcup_{g\in E} B_1(g^{-1}).
        \] If the pointwise stabilizer of $F$ in $A_2^3(1)$ is trivial,
        then the pointwise stabilizer of $B_1(1)$ in $\Aut(Z)$ is trivial and in particular the index of $\Gamma$ in $\Aut(Z)$ is at most $\abs{A^3_1(1)}$.
    \end{lemma}
    \begin{proof}
        If we have $g \in E$, then any $\sigma \in \Aut(Z)$ that fixes $1$ and $g$ fixes $B_1(1)$.
        Translating this phenomenon to an arbitrary $h \in Z$ yields that if $\sigma \in \Aut(Z)$ fixes $h$ and $hg$, then it fixes $B_1(h)$.
        In particular, every $\sigma \in \Aut(Z)$ that fixes $B_1(1)$ pointwise, fixes $F$ pointwise since it fixes $g^{-1}$ and $g^{-1}g =1$.
        Now assume that the pointwise stabilizer of $F$ in $A_2^3(1)$ is trivial. Then every automorphism of $Z$ that fixes $B_1(1)$ pointwise, fixes $F$ pointwise by the previous discussion, and therefore fixes $B_2(1)$ pointwise. It follows by induction that~$\sigma$ is trivial.
        Therefore, the projection of $\Stab(\Aut(Z), 1)$ to $A^1(1)$ is injective and the image lies in $A^3_1(1)$.
    \end{proof}
    We checked that the lattices $\Gamma_k^3$ satisfy the assumptions of Lemma \ref{lem:discreteness_criterion} and as it turns out the index of $\Gamma_k^3$ in the automorphism group of the Cayley graph (which equals $\Aut(X_k^3)$) is at most 4 in every case. On the other hand, the automorphism group of the complexes $Y_k^3$ is always 8. Therefore, the group of deck transformations covering $\Aut(Y_k^3)$ is $\Aut(X_k^3)$.

    \section{Small $\tilde A_2$-GABs and related data}\label{app:a2tildeGABs}
    
    In this section we list the 13 GABs from Theorem \ref{thm:heawood} and provide some more data on these GABs, their lattices and their buildings. The GABs all share the following combinatorial properties: they consists of three vertices $u,v,w$ and 21 edges $e_1,\dots,e_7, f_1, \dots, f_7, g_1, \dots, g_7$ and 21 triangles. The edges $e_i$ point from $w$ to $u$, the edges $f_i$ point from $u$ to $v$ and the edges $g_i$ point from $v$ to $w$. The GABs can be defined via the boundaries of their triangles and we do so in Table \ref{tab:deltas}.

    \begin{longtable}{ |c|c| }
    \hline
    GAB & Triangles \\
    \hline
    $\Delta_1$ & 
    \hspace{-3.5mm}\begin{minipage}{12.83cm}
        \begin{tabular}{c}
        $(e_1, f_1, g_1),  (e_1, f_2, g_3),  (e_1, f_4, g_2),  (e_2, f_1, g_4),  (e_2, f_3, g_1),  (e_2, f_7, g_5),  (e_3, f_1, g_2),$ \\
        $(e_3, f_5, g_6),  (e_3, f_6, g_4),  (e_4, f_2, g_1),  (e_4, f_3, g_6),  (e_4, f_5, g_7),  (e_5, f_2, g_7),  (e_5, f_6, g_3),$ \\
        $(e_5, f_7, g_4),  (e_6, f_3, g_5),  (e_6, f_4, g_3),  (e_6, f_6, g_6),  (e_7, f_4, g_5),  (e_7, f_5, g_2),  (e_7, f_7, g_7)$
        \end{tabular}
    \end{minipage}\\
    \hline
    $\Delta_2$ & 
    \hspace{-3.5mm}\begin{minipage}{12.83cm}
        \begin{tabular}{c}
            $(e_1, f_1, g_3),  (e_1, f_2, g_1),  (e_1, f_3, g_2),  (e_2, f_1, g_4),  (e_2, f_4, g_1),  (e_2, f_5, g_5),  (e_3, f_1, g_7),$ \\
            $(e_3, f_6, g_1),  (e_3, f_7, g_6),  (e_4, f_2, g_6),  (e_4, f_4, g_4),  (e_4, f_7, g_2),  (e_5, f_2, g_7),  (e_5, f_5, g_4),$ \\
            $(e_5, f_6, g_3),  (e_6, f_3, g_7),  (e_6, f_4, g_2),  (e_6, f_6, g_5),  (e_7, f_3, g_5),  (e_7, f_5, g_6),  (e_7, f_7, g_3)$
        \end{tabular}
    \end{minipage}\\
    \hline
    $\Delta_3$ & 
    \hspace{-3.5mm}\begin{minipage}{12.83cm}
        \begin{tabular}{c}
        $(e_1, f_1, g_2),  (e_1, f_2, g_3),  (e_1, f_3, g_1),  (e_2, f_1, g_4),  (e_2, f_4, g_5),  (e_2, f_5, g_1),  (e_3, f_1, g_7),$ \\
        $(e_3, f_6, g_1),  (e_3, f_7, g_6),  (e_4, f_2, g_6),  (e_4, f_4, g_2),  (e_4, f_6, g_4),  (e_5, f_2, g_7),  (e_5, f_5, g_2),$ \\
        $(e_5, f_7, g_5),  (e_6, f_3, g_7),  (e_6, f_4, g_3),  (e_6, f_7, g_4),  (e_7, f_3, g_5),  (e_7, f_5, g_6),  (e_7, f_6, g_3)$
        \end{tabular}
    \end{minipage}\\
    \hline
    $\Delta_4$ & 
    \hspace{-3.5mm}\begin{minipage}{12.83cm}
        \begin{tabular}{c}
            $(e_1, f_1, g_1),  (e_1, f_2, g_2),  (e_1, f_3, g_3),  (e_2, f_1, g_5),  (e_2, f_4, g_1),  (e_2, f_7, g_4),  (e_3, f_1, g_6),$ \\
            $(e_3, f_5, g_4),  (e_3, f_6, g_2),  (e_4, f_2, g_5),  (e_4, f_4, g_2),  (e_4, f_5, g_7),  (e_5, f_2, g_4),  (e_5, f_6, g_7),$ \\
            $(e_5, f_7, g_3),  (e_6, f_3, g_5),  (e_6, f_4, g_3),  (e_6, f_6, g_6),  (e_7, f_3, g_7),  (e_7, f_5, g_1),  (e_7, f_7, g_6)$
        \end{tabular}
    \end{minipage}\\
    \hline
    $\Delta_5$ & 
    \hspace{-3.5mm}\begin{minipage}{12.83cm}
        \begin{tabular}{c}
        $(e_1, f_1, g_2),  (e_1, f_2, g_3),  (e_1, f_3, g_1),  (e_2, f_1, g_1),  (e_2, f_4, g_5),  (e_2, f_7, g_4),  (e_3, f_1, g_4),$ \\
        $(e_3, f_5, g_2),  (e_3, f_6, g_6),  (e_4, f_2, g_1),  (e_4, f_4, g_7),  (e_4, f_5, g_6),  (e_5, f_2, g_7),  (e_5, f_6, g_4),$ \\
        $(e_5, f_7, g_3),  (e_6, f_3, g_5),  (e_6, f_4, g_2),  (e_6, f_6, g_7),  (e_7, f_3, g_6),  (e_7, f_5, g_3),  (e_7, f_7, g_5)$
        \end{tabular}
    \end{minipage}\\
    \hline
    $\Delta_6$ & 
    \hspace{-3.5mm}\begin{minipage}{12.83cm}
        \begin{tabular}{c}
        $(e_1, f_1, g_2),  (e_1, f_2, g_1),  (e_1, f_3, g_3),  (e_2, f_1, g_1),  (e_2, f_4, g_4),  (e_2, f_5, g_5),  (e_3, f_1, g_4),$ \\
        $(e_3, f_6, g_2),  (e_3, f_7, g_6),  (e_4, f_2, g_6),  (e_4, f_4, g_3),  (e_4, f_6, g_5),  (e_5, f_2, g_7),  (e_5, f_5, g_4),$ \\
        $(e_5, f_7, g_3),  (e_6, f_3, g_5),  (e_6, f_4, g_7),  (e_6, f_7, g_2),  (e_7, f_3, g_1),  (e_7, f_5, g_6),  (e_7, f_6, g_7)$
        \end{tabular}
    \end{minipage}\\
    \hline
    $\Delta_7$ & 
    \hspace{-3.5mm}\begin{minipage}{12.83cm}
        \begin{tabular}{c}
        $(e_1, f_1, g_2),  (e_1, f_2, g_7),  (e_1, f_3, g_5),  (e_2, f_1, g_1),  (e_2, f_4, g_4),  (e_2, f_7, g_5),  (e_3, f_1, g_3),$ \\
        $(e_3, f_5, g_7),  (e_3, f_6, g_4),  (e_4, f_2, g_4),  (e_4, f_4, g_2),  (e_4, f_5, g_6),  (e_5, f_2, g_1),  (e_5, f_6, g_6),$ \\
        $(e_5, f_7, g_7),  (e_6, f_3, g_6),  (e_6, f_4, g_5),  (e_6, f_6, g_3),  (e_7, f_3, g_1),  (e_7, f_5, g_2),  (e_7, f_7, g_3)$
        \end{tabular}
    \end{minipage}\\
    \hline
    $\Delta_8$ & 
    \hspace{-3.5mm}\begin{minipage}{12.83cm}
        \begin{tabular}{c}
        $(e_1, f_1, g_2),  (e_1, f_2, g_1),  (e_1, f_7, g_3),  (e_2, f_1, g_3),  (e_2, f_3, g_7),  (e_2, f_6, g_4),  (e_3, f_1, g_1),$ \\
        $(e_3, f_4, g_4),  (e_3, f_5, g_5),  (e_4, f_2, g_4),  (e_4, f_3, g_6),  (e_4, f_4, g_2),  (e_5, f_2, g_5),  (e_5, f_5, g_2),$ \\
        $(e_5, f_6, g_7),  (e_6, f_3, g_1),  (e_6, f_5, g_7),  (e_6, f_7, g_6),  (e_7, f_4, g_6),  (e_7, f_6, g_3),  (e_7, f_7, g_5)$
        \end{tabular}
    \end{minipage}\\
    \hline
    $\Delta_9$ & 
    \hspace{-3.5mm}\begin{minipage}{12.83cm}
        \begin{tabular}{c}
            $(e_1, f_2, g_1),  (e_1, f_4, g_2),  (e_1, f_7, g_3),  (e_2, f_1, g_1),  (e_2, f_6, g_4),  (e_2, f_7, g_5),  (e_3, f_3, g_1),$ \\
            $(e_3, f_5, g_6),  (e_3, f_7, g_7),  (e_4, f_2, g_4),  (e_4, f_5, g_2),  (e_4, f_6, g_6),  (e_5, f_1, g_2),  (e_5, f_2, g_7),$ \\
            $(e_5, f_3, g_5),  (e_6, f_3, g_6),  (e_6, f_4, g_5),  (e_6, f_6, g_3),  (e_7, f_1, g_3),  (e_7, f_4, g_4),  (e_7, f_5, g_7)$
        \end{tabular}
    \end{minipage}\\
    \hline
    $\Delta_{10}$ & 
    \hspace{-3.5mm}\begin{minipage}{12.83cm}
        \begin{tabular}{c}
            $(e_1, f_2, g_1),  (e_1, f_5, g_2),  (e_1, f_6, g_3),  (e_2, f_3, g_1),  (e_2, f_4, g_5),  (e_2, f_5, g_4),  (e_3, f_1, g_1),$ \\
            $(e_3, f_5, g_7),  (e_3, f_7, g_6),  (e_4, f_1, g_6),  (e_4, f_2, g_4),  (e_4, f_4, g_2),  (e_5, f_2, g_3),  (e_5, f_3, g_7),$ \\
            $(e_5, f_7, g_4),  (e_6, f_4, g_3),  (e_6, f_6, g_6),  (e_6, f_7, g_5),  (e_7, f_1, g_2),  (e_7, f_3, g_5),  (e_7, f_6, g_7)$
        \end{tabular}
    \end{minipage}\\
    \hline
    $\Delta_{11}$ & 
    \hspace{-3.5mm}\begin{minipage}{12.83cm}
        \begin{tabular}{c}
            $(e_1, f_1, g_3),  (e_1, f_2, g_5),  (e_1, f_3, g_6),  (e_2, f_1, g_1),  (e_2, f_4, g_2),  (e_2, f_7, g_3),  (e_3, f_1, g_2),$ \\
            $(e_3, f_5, g_6),  (e_3, f_6, g_4),  (e_4, f_2, g_1),  (e_4, f_4, g_4),  (e_4, f_5, g_5),  (e_5, f_2, g_4),  (e_5, f_6, g_3),$ \\
            $(e_5, f_7, g_7),  (e_6, f_3, g_1),  (e_6, f_4, g_7),  (e_6, f_6, g_6),  (e_7, f_3, g_7),  (e_7, f_5, g_2),  (e_7, f_7, g_5)$
        \end{tabular}
    \end{minipage}\\
    \hline
    $\Delta_{12}$ & 
    \hspace{-3.5mm}\begin{minipage}{12.83cm}
        \begin{tabular}{c}
            $(e_1, f_1, g_2),  (e_1, f_2, g_1),  (e_1, f_4, g_3),  (e_2, f_1, g_1),  (e_2, f_5, g_4),  (e_2, f_7, g_5),  (e_3, f_1, g_6),$ \\
            $(e_3, f_3, g_1),  (e_3, f_6, g_7),  (e_4, f_4, g_2),  (e_4, f_6, g_6),  (e_4, f_7, g_4),  (e_5, f_3, g_7),  (e_5, f_4, g_5),$ \\
            $(e_5, f_5, g_2),  (e_6, f_2, g_3),  (e_6, f_3, g_5),  (e_6, f_7, g_6),  (e_7, f_2, g_4),  (e_7, f_5, g_7),  (e_7, f_6, g_3)$
        \end{tabular}
    \end{minipage}\\
    \hline
    $\Delta_{13}$ & 
    \hspace{-3.5mm}\begin{minipage}{12.83cm}
        \begin{tabular}{c}
            $(e_1, f_1, g_2),  (e_1, f_3, g_1),  (e_1, f_6, g_3),  (e_2, f_1, g_1),  (e_2, f_4, g_5),  (e_2, f_7, g_4),  (e_3, f_1, g_7),$ \\
            $(e_3, f_2, g_1),  (e_3, f_5, g_6),  (e_4, f_2, g_4),  (e_4, f_3, g_6),  (e_4, f_4, g_2),  (e_5, f_3, g_5),  (e_5, f_5, g_2),$ \\
            $(e_5, f_7, g_7),  (e_6, f_4, g_3),  (e_6, f_5, g_4),  (e_6, f_6, g_7),  (e_7, f_2, g_3),  (e_7, f_6, g_6),  (e_7, f_7, g_5)$
        \end{tabular}
    \end{minipage}\\
    \hline
    \caption{The GABs $\Delta_1, \dots, \Delta_{13}$.}
    \label{tab:deltas}
    \end{longtable}

    \vspace{5mm}

    \begin{longtable}{|c|c|}
            \hline
            GAB & Generators of $\Aut(\Delta)$ \\
            \hline
            $\Delta_1$
            &
            \begin{minipage}{11cm}
                {
                \vspace{-2.5mm}
                \begin{gather*}
                    (e_2\, e_7\, e_5) (e_3\, e_6\, e_4)
                    (f_1\, f_4\, f_2)(f_3\, f_5\, f_6)
                    (g_1\, g_2\, g_3)(g_4\, g_5\, g_7),
                    \\
                    (e_1\, g_3\, f_6\, e_3\, g_2\, f_4\, e_6\, g_6\, f_5\, e_7\, g_5\, f_3\, e_4\, g_7\, f_7\, e_2\, g_1\, f_2\, e_5\, g_4\, f_1)
                \end{gather*}
                \vspace{-5mm}
                }
            \end{minipage}\\
            \hline
            $\Delta_2$
            &
            \begin{minipage}{11cm}
                {
                \vspace{-2.5mm}
                \begin{gather*}
                    (e_1\, e_5)(e_2\, e_6)(e_4\, e_7)
                    (f_1\, g_7)(f_2\, g_3)(f_3\, g_4)(f_4\, g_5)
                    (f_5\, g_2)(f_6\, g_1)(f_7\, g_6),
                    \\
                    (e_1\, g_1\, f_2)
                    (e_2\, g_6\, f_3)
                    (e_3\, g_7\, f_1)
                    (e_4\, g_2\, f_4)
                    (e_5\, g_3\, f_6)
                    (e_6\, g_4\, f_7)
                    (e_7\, g_5\, f_5)
                \end{gather*}
                \vspace{-5mm}
                }
            \end{minipage}\\
            \hline
            $\Delta_3$
            &
            \begin{minipage}{11cm}
                {
                \vspace{-2.5mm}
                \begin{gather*}
                    (e_1\, e_6\, e_5\, e_7\, e_3\, e_2\, e_4)
                    (f_1\, f_4\, f_2\, f_3\, f_7\, f_5\, f_6)
                    (g_1\, g_4\, g_2\, g_3\, g_7\, g_5\, g_6),
                    \\
                    (e_1\, f_7)(e_2\, f_5\, e_7\, f_6\, e_4\, f_4)
                    (e_3\, f_2\, e_6\, f_1\, e_5\, f_3)
                    (g_1\, g_6\, g_3\, g_4\, g_2\, g_5)
                \end{gather*}
                \vspace{-5mm}
                }
            \end{minipage}\\
            \hline
            $\Delta_4$
            &
            \begin{minipage}{11cm}
            {
            \vspace{-2.5mm}
            \begin{gather*}
            (e_3\, e_7)(e_4\, e_6)
            (f_1\, g_1)(f_2\, g_3)(f_3\, g_2)
            (f_4\, g_5)(f_5\, g_6)(f_6\, g_7)(f_7\, g_4),
            \\
            (e_1\, g_5\, f_4)(e_2\, g_1\, f_1)(e_3\, g_4\, f_5)
            (e_4\, g_2\, f_2)(e_5\, g_7\, f_6)(e_6\, g_3\, f_3)(e_7\, g_6\, f_7),
            \\
            (e_1\, g_5\, f_4)(e_2\, g_3\, f_2)(e_3\, g_6\, f_6)
            (e_4\, g_1\, f_3)(e_5\, g_4\, f_7)(e_6\, g_2\, f_1)(e_7\, g_7\, f_5)
            \end{gather*}
            \vspace{-5mm}
            }
            \end{minipage}\\
            \hline
            $\Delta_5$
            &
            \begin{minipage}{11cm}
            {
            \vspace{-2.5mm}
            \begin{gather*}
                (e_1\, e_6)(e_3\, e_7)(f_1\, g_5)(f_2\, g_7)(f_3\, g_2)(f_4\, g_1)(f_5\, g_6)(f_6\, g_3)(f_7\, g_4),
                \\
                (e_1\, g_6\, f_4)
                (e_2\, g_3\, f_6)
                (e_3\, g_5\, f_2)
                (e_4\, g_2\, f_3)
                (e_5\, g_4\, f_7)
                (e_6\, g_1\, f_5)
                (e_7\, g_7\, f_1)
            \end{gather*}
            \vspace{-5mm}
            }
            \end{minipage}\\
            \hline
            $\Delta_6$
            &
            \begin{minipage}{11cm}
            {
            \vspace{-2.5mm}
            \begin{gather*}
                (e_1\, e_2\, e_3)(e_5\, e_6\, e_7)(f_2\, f_4\, f_6)(f_3\, f_5\, f_7)(g_1\, g_4\, g_2)(g_3\, g_5\, g_6),
                \\
                (e_1\, g_1)(e_2\, g_2)(e_3\, g_4)(e_4\, g_7)(e_5\, g_6)(e_6\, g_5)(e_7\, g_3)(f_4\, f_6)(f_5\, f_7)
            \end{gather*}
            \vspace{-5mm}
            }
            \end{minipage}\\
            \hline
            $\Delta_7$
            &
            \begin{minipage}{11cm}
            {
            \vspace{-2.5mm}
            \begin{gather*}
            (e_1\, f_7)(e_2\, f_3)(e_3\, f_5)(e_4\, f_6)(e_5\, f_2)(e_6\, f_4)(e_7\, f_1)(g_2\, g_3)(g_4\, g_6)
            \end{gather*}
            \vspace{-5mm}
            }
            \end{minipage}\\
            \hline
            $\Delta_8$
            &
            \begin{minipage}{11cm}
            {
            \vspace{-2.5mm}
            \begin{gather*}
                (e_1\, e_4\, e_5)(e_2\, e_7\, e_6)(f_1\, f_4\, f_5)(f_3\, f_6\, f_7)(g_1\, g_4\, g_5)(g_3\, g_6\, g_7),
                \\
                (e_1\, f_1)(e_2\, f_7)(e_3\, f_2)(e_4\, f_5)(e_5\, f_4)(e_6\, f_3)(e_7\, f_6)(g_4\, g_5)(g_6\, g_7)
            \end{gather*}
            \vspace{-5mm}
            }
            \end{minipage}\\
            \hline
            $\Delta_9$
            &
            \begin{minipage}{11cm}
            {
            \vspace{-2.5mm}
            \begin{gather*}
                (e_2\, e_5\, e_7)(e_3\, e_4\, e_6)(f_2\, f_4\, f_7)(f_3\, f_5\, f_6)(g_1\, g_2\, g_3)(g_4\, g_5\, g_7),
                \\
                (e_1\, f_1)(e_2\, f_2)(e_3\, f_3)(e_4\, f_6)(e_5\, f_7)(e_6\, f_5)(e_7\, f_4)(g_2\, g_3)(g_5\, g_7)
            \end{gather*}
            \vspace{-5mm}
            }
            \end{minipage}\\
            \hline
            $\Delta_{10}$
            &
            \begin{minipage}{11cm}
            {
            \vspace{-2.5mm}
            \begin{gather*}
            (e_1\, e_2\, e_3)(e_4\, e_5\, e_7)(f_1\, f_2\, f_3)(f_4\, f_7\, f_6)(g_2\, g_4\, g_7)(g_3\, g_5\, g_6)
            \end{gather*}
            \vspace{-5mm}
            }
            \end{minipage}\\
            \hline
            $\Delta_{11}$
            &
            \begin{minipage}{11cm}
            {
            \vspace{-2.5mm}
            \begin{gather*}
            (e_1\, f_4)(e_2\, f_2\, e_6\, f_1\, e_4\, f_3)(e_3\, f_5\, e_7\, f_7\, e_5\, f_6)(g_2\, g_5\, g_7\, g_3\, g_4\, g_6)
            \end{gather*}
            \vspace{-5mm}
            }
            \end{minipage}\\
            \hline
            $\Delta_{12}$
            &
            \begin{minipage}{11cm}
            {
            \vspace{-2.5mm}
            \begin{gather*}
            (e_2\, e_4\, e_6)(e_3\, e_5\, e_7)(f_1\, f_4\, f_2)(f_3\, f_5\, f_6)(g_1\, g_2\, g_3)(g_4\, g_6\, g_5)
            \end{gather*}
            \vspace{-5mm}
            }
            \end{minipage}\\
            \hline
            $\Delta_{13}$
            &
            \begin{minipage}{11cm}
            {
            \vspace{-2.5mm}
            \begin{gather*}
            (e_1\, e_2\, e_7)(e_3\, e_5\, e_6)(f_1\, f_7\, f_6)(f_2\, f_3\, f_4)(g_1\, g_5\, g_3)(g_2\, g_4\, g_6)
            \end{gather*}
            \vspace{-5mm}
            }
            \end{minipage}\\
            \hline
        \caption{Automorphism groups of the GABs $\Delta_1, \dots, \Delta_{13}$.}
        \label{tab:autos_deltas}
    \end{longtable}

    \vspace{5mm}

    \begin{longtable}{|c|c|c|}
            \hline
            GAB & Corresponding automorphisms in $\Aut(\Delta)$ & Extension \\
            \hline
            $\Delta_1$
            &
            \hspace{-2mm}\begin{minipage}{10.5cm}
            {
            \vspace{-2.5mm}
                \begin{gather*}
                    \big\langle
                    (e_1\, e_6\, e_4\, e_5\, e_3\, e_7\, e_2)
                    (f_1\, f_4\, f_3\, f_2\, f_6\, f_5\, f_7)
                    (g_1\, g_3\, g_6\, g_7\, g_4\, g_2\, g_5)
                    \big \rangle
                \end{gather*}
            \vspace{-5mm}
            }
            \end{minipage}
            &
            $\Sigma_1$\\
            \hline
            $\Delta_1$
            &
            \hspace{-2mm}\begin{minipage}{10.5cm}
            {
            \vspace{-2.5mm}
                \begin{gather*}
                    \big\langle
                    (e_2\, e_7\, e_5)(e_3\, e_6\, e_4)
                    (f_1\, f_4\, f_2)(f_3\, f_5\, f_6)
                    (g_1\, g_2\, g_3)(g_4\, g_5\, g_7),
                    \\
                    (e_1\, e_6\, e_4\, e_5\, e_3\, e_7\, e_2)
                    (f_1\, f_4\, f_3\, f_2\, f_6\, f_5\, f_7)
                    (g_1\, g_3\, g_6\, g_7\, g_4\, g_2\, g_5)
                    \big \rangle
                \end{gather*}
            \vspace{-5mm}
            }
            \end{minipage}
            &
            $\Omega_3$\\
            \hline
            $\Delta_1$
            &
            \hspace{-2mm}\begin{minipage}{10.5cm}
            {
            \vspace{-2.5mm}
                \begin{gather*}
                    \big\langle
                    (e_1\, g_6\, f_7)
                    (e_2\, g_3\, f_5)
                    (e_3\, g_5\, f_2)
                    (e_4\, g_4\, f_4)
                    (e_5\, g_2\, f_3)
                    (e_6\, g_7\, f_1)
                    (e_7\, g_1\, f_6)
                    \big \rangle
                \end{gather*}
            \vspace{-5mm}
            }
            \end{minipage}
            &
            $\Lambda_{A.1}$\\
            \hline
            $\Delta_1$
            &
            \hspace{-2mm}\begin{minipage}{10.5cm}
            {
            \vspace{-2.5mm}
                \begin{gather*}
                    \big\langle
                    (e_1\, g_2\, f_1)
                    (e_2\, g_3\, f_5)
                    (e_3\, g_1\, f_4)
                    (e_4\, g_5\, f_6)
                    (e_5\, g_7\, f_7)
                    (e_6\, g_6\, f_3)
                    (e_7\, g_4\, f_2)
                    \big \rangle
                \end{gather*}
            \vspace{-5mm}
            }
            \end{minipage}
            &
            $\Lambda_{A.2}$\\
            \hline
            $\Delta_1$
            &
            \hspace{-2mm}\begin{minipage}{10.5cm}
            {
            \vspace{-2.5mm}
                \begin{gather*}
                    \big\langle
                    (e_1\, g_1\, f_1)
                    (e_2\, g_2\, f_2)
                    (e_3\, g_3\, f_3)
                    (e_4\, g_4\, f_4)
                    (e_5\, g_5\, f_5)
                    (e_6\, g_6\, f_6)
                    (e_7\, g_7\, f_7)
                    \big \rangle
                \end{gather*}
            \vspace{-5mm}
            }
            \end{minipage}
            &
            $\Lambda_{A.3}$\\
            \hline
            $\Delta_2$
            &
            \hspace{-2mm}\begin{minipage}{10.5cm}
            {
            \vspace{-2.5mm}
                \begin{gather*}
                    \big\langle
                    (e_1 \, f_2 \, g_1) (e_2 \, f_3 \, g_6) (e_3 \, f_1 \, g_7) (e_4 \, f_4 \, g_2)
                    (e_5 \, f_6 \, g_3) (e_6 \, f_7 \, g_4) (e_7 \, f_5 \, g_5)
                    \big\rangle
                \end{gather*}
            \vspace{-5mm}
            }
            \end{minipage}
            &
            $\Lambda_{A.4}$\\
            \hline
            $\Delta_3$
            &
            \hspace{-2mm}\begin{minipage}{10.5cm}
            {
            \vspace{-2.5mm}
                \begin{gather*}
                    \big\langle
                    (e_1 \, e_6 \, e_5 \, e_7 \, e_3 \, e_2 \, e_4)
                    (f_1 \, f_4 \, f_2 \, f_3 \, f_7 \, f_5 \, f_6)
                    (g_1 \, g_4 \, g_2 \, g_3 \, g_7 \, g_5 \, g_6)
                    \big\rangle
                \end{gather*}
            \vspace{-5mm}
            }
            \end{minipage}
            &
            $\Sigma_1$\\
            \hline
            $\Delta_3$
            &
            \hspace{-2mm}\begin{minipage}{10.5cm}
            {
            \vspace{-2.5mm}
            \begin{gather*}
            \big\langle
            (e_2\, e_7\, e_4)(e_3\, e_6\, e_5)
            (f_1\, f_3\, f_2)(f_4\, f_5\, f_6)
            (g_1\, g_3\, g_2)(g_4\, g_5\, g_6),
            \\
            (e_1\, e_6\, e_5\, e_7\, e_3\, e_2\, e_4)
            (f_1\, f_4\, f_2\, f_3\, f_7\, f_5\, f_6)
            (g_1\, g_4\, g_2\, g_3\, g_7\, g_5\, g_6)
            \big\rangle
            \end{gather*}
            \vspace{-5mm}
            }
            \end{minipage}
            &
            $\Omega_3$\\
            \hline
            $\Delta_4$
            &
            \hspace{-2mm}\begin{minipage}{10.5cm}
            {
            \vspace{-2.5mm}
                \begin{gather*}
                    \big\langle
                    (e_1\, g_5\, f_4)
                    (e_2\, g_2\, f_3)
                    (e_3\, g_7\, f_7)
                    (e_4\, g_3\, f_1)
                    (e_5\, g_6\, f_5)
                    (e_6\, g_1\, f_2)
                    (e_7\, g_4\, f_6)
                    \big\rangle
                \end{gather*}
            \vspace{-5mm}
            }
            \end{minipage}
            &
            $\Lambda_{B.1}$\\
            \hline
            $\Delta_4$
            &
            \hspace{-2mm}\begin{minipage}{10.5cm}
            {
            \vspace{-2.5mm}
                \begin{gather*}
                    \big\langle
                    (e_1 \, f_4 \, g_5) (e_2 \, f_2 \, g_3) (e_3 \, f_6 \, g_6) (e_4 \, f_3 \, g_1)
                (e_5 \, f_7 \, g_4) (e_6 \, f_1 \, g_2) (e_7 \, f_5 \, g_7)
                    \big\rangle
                \end{gather*}
            \vspace{-5mm}
            }
            \end{minipage}
            &
            $\Lambda_{B.2}$\\
            \hline
            $\Delta_4$
            &
            \hspace{-2mm}\begin{minipage}{10.5cm}
            {
            \vspace{-2.5mm}
                \begin{gather*}
                    \big\langle
                    (e_1 \, f_4 \, g_5) (e_2 \, f_1 \, g_1) (e_3 \, f_5 \, g_4) (e_4 \, f_2 \, g_2)
                (e_5 \, f_6 \, g_7) (e_6 \, f_3 \, g_3) (e_7 \, f_7 \, g_6)
                    \big\rangle
                \end{gather*}
            \vspace{-5mm}
            }
            \end{minipage}
            &
            $\Lambda_{B.3}$\\
            \hline
            $\Delta_5$
            &
            \hspace{-2mm}\begin{minipage}{10.5cm}
            {
            \vspace{-2.5mm}
                \begin{gather*}
                    \big\langle
                    (e_1 \, f_4 \, g_6) (e_2 \, f_6 \, g_3) (e_3 \, f_2 \, g_5) (e_4 \, f_3 \, g_2)
                (e_5 \, f_7 \, g_4) (e_6 \, f_5 \, g_1) (e_7 \, f_1 \, g_7)
                    \big\rangle
                \end{gather*}
            \vspace{-5mm}
            }
            \end{minipage}
            &
            $\Lambda_{C.1}$\\
            \hline
        \caption{Extensions of the lattices $\Gamma_1, \dots, \Gamma_5$. The extension is the group of deck transformations covering the automorphisms indicated.}
        \label{tab:extensions}
    \end{longtable}

    \vspace{5mm}

    \begin{longtable}{|c|c|c|c|c|c|c|c|}
            \hline
            Building & $X_7$ & $X_8$ & $X_9$ & $X_{10}$ & $X_{11}$ & $X_{12}$ & $X_{13}$ \\
            \hline
            $|A(u))|$ & 10752 & 1536 & 1536 & 1536 & 3072 & 1536 & 1536\\
            \hline
            $|A^1(u)|$ & 42 & 6 & 6 & 3 & 6 & 3 & 3\\ 
            \hline
            $|A_{(1)}^3(u)|$ & 1 & 1 & 1 & 1 & 1 & 1 & 1\\
            \hline
        \caption[Balls of radius four in the exotic buildings.]{Balls of radius four in the exotic buildings $X_7,\dots, X_{13}$.
        We denote by $A(u)$ the automorphism group of a ball of radius four in the building around a lift $\tilde u$ of $u$.
        We denote by $A^1(u)$ the group of automorphisms on the 1-ball around $\tilde u$ induced by $A(u)$. We denote by $A^3_{(1)}(u)$ the groups of automorphism of the 3-ball induced by the pointwise stabilizer of the 1-ball in $A(u)$.} 
        \label{tab:fourballs}
    \end{longtable}

    \vspace{5mm}

    \begin{longtable}{|c|c|c|c|c|c|c|}
            \hline
            Lattice & $\Gamma_1$ & $\Gamma_2$ & $\Gamma_3$ & $\Gamma_4$ & $\Gamma_5$ & $\Gamma_6$ \\
            \hline
            $|\Gamma : \Gamma^{(2)}|$ & 14336 & 441 & 448 & 8 & 16 & 1\\
            \hline
        \caption[The lattices $\Gamma_1, \dots, \Gamma_6$ and their quotients modulo their second derived subgroups.]{The lattices $\Gamma_1, \dots, \Gamma_6$ can be distinguished by their quotients modulo their second derived subgroups.}
        \label{tab:quotients_secondderived}
    \end{longtable}

    \printbibliography 
    
\end{document}